\documentclass{article} 
\usepackage{latexsym,amsmath,amscd,amssymb,graphics}
\usepackage[backend=bibtex,sorting=none]{biblatex}
\usepackage{enumerate}
\usepackage[]{mcode}
\usepackage[toc,page]{appendix}
\usepackage[a4paper,margin=.75in]{geometry}
\usepackage{array}
\usepackage[section]{placeins}
\usepackage{afterpage}
\usepackage{verbatim}

\usepackage{graphicx}
\usepackage{subfig}
\usepackage{framed}
\usepackage[colorlinks]{hyperref}
\usepackage{url}
\usepackage{mathrsfs}
\usepackage{epstopdf}

\usepackage[all]{xy}
\usepackage{indentfirst}
\usepackage[mathscr]{euscript}
\usepackage{algorithm}
\usepackage{algpseudocode}
\usepackage{authblk}

\makeatletter
\renewcommand*{\@fnsymbol}[1]{\ensuremath{\ifcase#1\or \dagger\or *\or \ddagger\or
		\mathsection\or \mathparagraph\or \|\or **\or \dagger\dagger
		\or \ddagger\ddagger \else\@ctrerr\fi}}
\makeatother

\makeatletter

\@addtoreset{figure}{section}
\def\thefigure{\thesection.\@arabic\c@figure}
\def\fps@figure{h, t}
\@addtoreset{table}{bsection}
\def\thetable{\thesection.\@arabic\c@table}
\def\fps@table{h, t}
\@addtoreset{equation}{section}

\makeatother


\makeatletter
\AtBeginDocument{%
	\expandafter\renewcommand\expandafter\subsection\expandafter{%
		\expandafter\@fb@secFB\subsection
	}%
}
\makeatother

\newtheorem{theorem}{Theorem}

\newtheorem{remark}[theorem]{Remark}

\numberwithin{theorem}{subsection}

\def\bea{\begin{eqnarray}}
\def\eea{\end{eqnarray}}
\def\ba{\begin{array}}
\def\ea{\end{array}}

\def\bomega{\boldsymbol{\omega}}
\def\bOm{\boldsymbol{\Omega}}

\def\diag{\mathrm{ \textbf{diag}}}

\def\bx{{\boldsymbol {x} }}

\let\<\langle
\let \>\rangle

\newcommand{\rem}[1]{}

\newcommand{\de}{\delta}

\newcommand{\bA}{\boldsymbol{A}}

\newcommand{\bu}{\boldsymbol{u}}
\newcommand{\bPsi}{\boldsymbol{\Psi}}
\newcommand{\bDelta}{\boldsymbol{\Delta}}
\newcommand{\bv}{\boldsymbol{v}}
\newcommand{\bz}{\boldsymbol{z}}

\newcommand{\br}{\boldsymbol{r}}

\newcommand{\bGam}{\boldsymbol{\Gamma}}
\newcommand{\bGamma}{\boldsymbol{\Gamma}}
\newcommand{\bzeta}{\boldsymbol{\zeta}}
\newcommand{\bom}{\boldsymbol{\omega}}

\newcommand{\bsigma}{\boldsymbol{\sigma}}
\newcommand{\bSigma}{\boldsymbol{\Sigma}}
\newcommand{\bpsi}{\boldsymbol{\psi}}
\newcommand{\balpha}{\boldsymbol{\alpha}}

\newcommand{\bxi}{\boldsymbol{\xi}}

\newcommand{\bs}{\mathbf{s}}

\newcommand{\bY}{\mathbf{Y}}

\newcommand{\bPi}{\boldsymbol{\Pi}}

\newcommand{\bkappa}{\boldsymbol{\kappa}}
\newcommand{\bpi}{\boldsymbol{\pi}}
\newcommand{\blam}{\boldsymbol{\lambda}}
\newcommand{\bnu}{\boldsymbol{\nu}}
\newcommand{\bk}{\mathbf{k}}

\newcommand{\btheta}{\boldsymbol{\theta}}
\newcommand{\inertia}{\mathbb{I}}

\newcommand{\pp}[2]{\frac{\partial #1}{\partial #2}}
\newcommand{\dd}[2]{\frac{\mathrm{d} #1}{\mathrm{d} #2}}

\newcommand{\dt}{\mathrm{d}t}

\newcommand{\Om}{\Omega}

\newcommand{\mse}{\mathfrak{se}}

\newcommand{\dprime}{\prime \prime}
\newcommand{\thprime}{\prime \prime \prime}
\newcommand{\Sym}{\textbf{Sym}}
\DeclareMathOperator{\ReLU}{ReLU}

\newcommand{\todo}[1]{\vspace{5 mm}\par \noindent
\framebox{\begin{minipage}[c]{0.95 \textwidth}
\tt #1 \end{minipage}}\vspace{5 mm}\par}

\newcommand{\revision}[2]{#2} 

\newcommand{\revisionS}[2]{#2} 

\newcommand{\revisionT}[2]{#2}

\usepackage[normalem]{ulem}

\graphicspath{ {./gen/} {./bsim10/}  }

\title{On the Optimal Control of a Rolling Ball Robot Actuated by Internal Point Masses}
\author[1,2]{Vakhtang Putkaradze\thanks{Email address: \texttt{putkarad@ualberta.ca}}}
\author[3]{Stuart Rogers\thanks{Email address: \texttt{srogers@umn.edu}}}
\affil[1]{Department of Mathematical and Statistical Sciences, Faculty of Science, University of Alberta, CAB 632, Edmonton, AB T6G 2G1, Canada}
\affil[2]{ATCO SpaceLab, 5302 Forand ST SW, Calgary, AB T3E 8B4, Canada} 
\affil[3]{Institute for Mathematics and its Applications, College of Science and Engineering, University of Minnesota, 207 Church ST SE, 306 Lind Hall, Minneapolis, MN 55455, USA}
\date{\today}

\providecommand{\keywords}[1]{\textbf{\textit{Keywords:}} #1}
\begin{document}

\maketitle

\abstract{\noindent 
The controlled motion of a rolling ball actuated by internal point masses that move along arbitrarily-shaped rails fixed within the ball is considered. Application of the  variational  Pontryagin's minimum principle yields the ball's controlled equations of motion, a solution of which obeys the ball's uncontrolled equations of motion, satisfies prescribed initial and final conditions, and minimizes a prescribed performance index.} 
\\
\\
\keywords{optimal control, rolling ball robots}
\tableofcontents 

\section{Introduction} \label{sec_introduction}

\subsection{Motivation and Methodology}
\begin{figure}[h] 
	\centering
	\subfloat[Sphero's toy incarnation of \textit{Star Wars'} BB-9E, the evil nemesis of BB-8  \cite{BB_9E}.]{\includegraphics[scale=.1025]{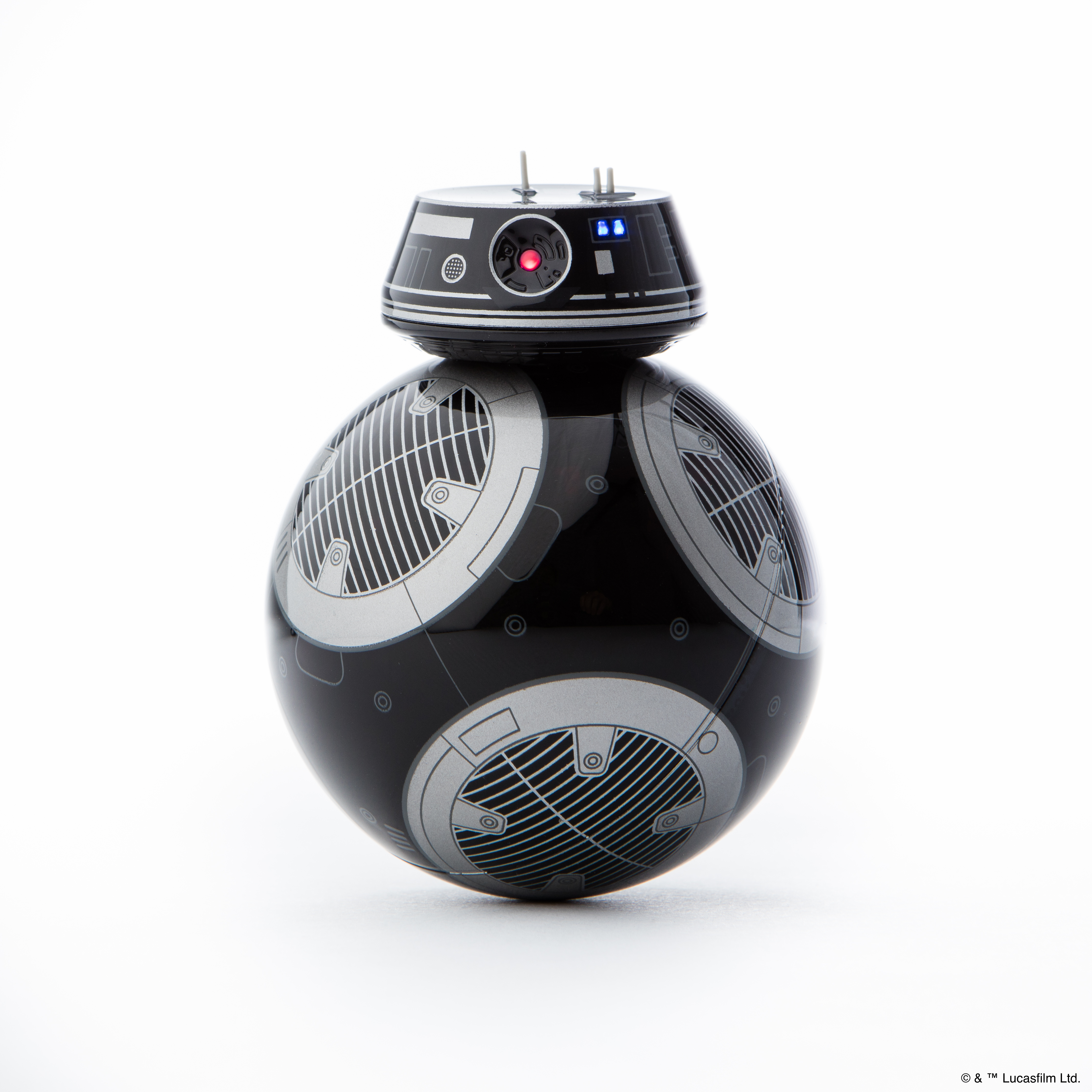}\label{fig_bb9e}}
	\hspace{5mm}
	\subfloat[Rosphere can be used in agriculture for monitoring crops, \copyright\ 2013 Emerald \cite{hernandez2013moisture}.]{\includegraphics[scale=.151]{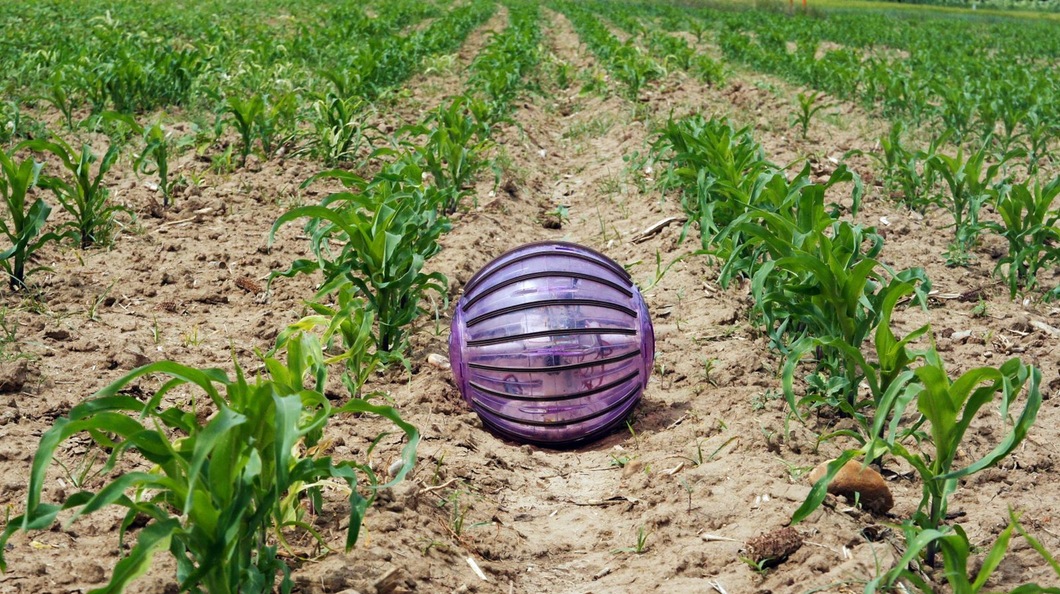}\label{fig_rosphere}}
	\caption{Examples of real rolling ball robots.}
\end{figure}

The first six films in the famous \textit{Star Wars} space saga starred the sidekick robot R2-D2, which looks like a dome-topped cylindrical trash can attached to a three-wheeled tripod.  A new generation of sidekick robots, the BB-series astromech droids 2BB-2, BB-4, BB-8, and BB-9E, were introduced in the seventh and eighth films in that saga, \textit{The Force Awakens} and \textit{The Last Jedi}. Unlike R2-D2, an astromech BB-unit consists of a rolling ball topped by a free-moving sensor and control unit as depicted in Figure~\ref{fig_bb9e}. 
While these rolling ball robots were mostly computer-generated in the films, the toy company Sphero sells working toy models of BB-8 and BB-9E which can be remotely controlled via a smartphone app.  However, rolling ball robots are not just gimmicks used by the entertainment and toy industries. The law enforcement, security, defense, energy, and agricultural industries are also interested in exploiting sensor-equipped rolling ball robots, such as Rosphere shown in Figure~\ref{fig_rosphere}, for such tasks as surveillance, observation, and environmental monitoring. 
More specifically, suppose a rolling ball like an astromech BB-unit or Rosphere is actuated by some internal mechanism which may be controlled, for example, by spinning internal rotors, by swinging an internal pendulum, or by moving internal point masses along arbitrarily-shaped rails fixed within the ball. In addition, suppose initial and final conditions, like the ball's initial and final positions and velocities, algebraic (i.e. non-differential) path inequality constraints, like engineering limitations on the internal mechanism's acceleration, and a performance index, such as the mean error between the ball's actual and prescribed trajectory, are prescribed. How can the ball's internal mechanism be controlled to minimize the prescribed performance index while satisfying the prescribed initial and final conditions and the prescribed  algebraic  path inequality constraints?

This paper attempts to answer  this question using a theoretical approach  assuming that the ball is actuated by internal point masses traveling along fixed trajectories in the ball's coordinate frame, such as arbitrarily-shaped rails. This problem, as we discuss in the next subsection, has not received sufficient attention in the literature due to the relative complexity of the theoretical considerations. We believe that  internal control mechanisms of this type for rolling ball robots may be beneficial both for practical applications and  for theoretical studies of geometric control theory.  Recent papers \cite{Putkaradze2018dynamicsP,putkaradze2018normalpub} by the authors deal exclusively with the derivation and analysis of the uncontrolled equations of motion and the contact point forces  (a.k.a the dynamics) for rolling ball robots  actuated by internal point masses. The goal of this paper is to apply the theory of optimal control to formulate the controlled equations of motion of  such rolling ball robots, to locomote over a prescribed trajectory, avoid obstacles, and/or perform some other maneuver by minimizing a prescribed  performance index. A subsequent paper \cite{RBOCnumerics} numerically solves particular cases of these controlled equations of motion.

Before optimal control can be applied to the rolling ball,  or, for that matter, to any other dynamical system, the equations of motion must be derived first.  Henceforth, the ordinary differential equations of motion of a dynamical system will be referred to as  the dynamics, the equations of motion, or the uncontrolled equations of motion to distinguish them from the controlled equations of motion which may be obtained by the indirect method of optimal control. The  dynamics of a rolling ball  actuated by internal  point  masses moving  along fixed, arbitrarily-shaped trajectories inside the ball is nontrivial, and, as far as we are aware, was only solved in this  generality in recent papers by the authors \cite{Putkaradze2018dynamicsP,putkaradze2018normalpub}. For  conciseness,  we shall not rederive  the uncontrolled  equations of motion  here,  but instead refer the  interested reader to those papers \cite{Putkaradze2018dynamicsP,putkaradze2018normalpub}. 

\subsection{Background}

\begin{figure}[h]
	\centering
	\includegraphics[width=0.5\linewidth]{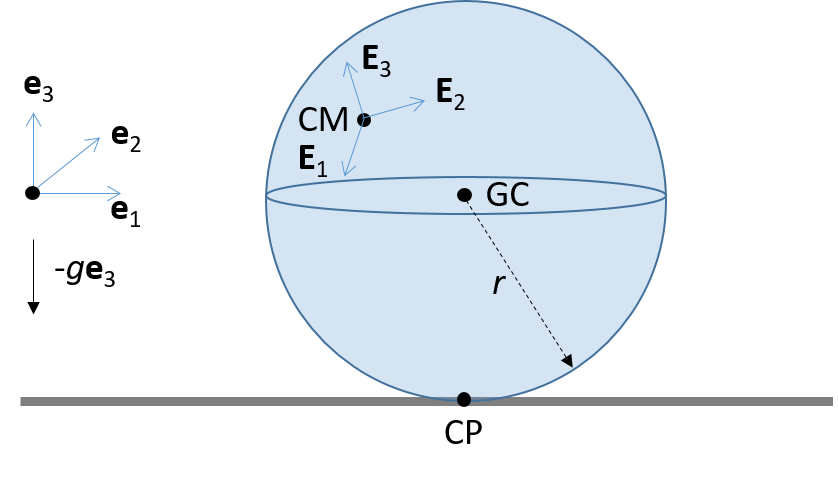}
	\caption{A ball of radius $r$ rolls without slipping on a horizontal surface in the presence of a uniform gravitational field of magnitude $g$. The ball's geometric center, center of mass, and contact point with the horizontal surface are denoted by GC, CM, and CP, respectively.}
	\label{fig:simple_rolling_ball}
\end{figure}

Consider a ball rolling without slipping on a horizontal surface in the presence of a uniform gravitational field. Figure~\ref{fig:simple_rolling_ball} shows a ball of radius $r$ rolling without slipping on a horizontal surface in the presence of a uniform gravitational field of magnitude $g$.

There are several terminologies in the literature to describe a ball rolling without slipping on a horizontal surface in the presence of a uniform gravitational field, depending on its mass distribution and the location of its center of mass \footnote{ A historical note on the terminology is warranted here to avoid confusion.  A Chaplygin sphere is a ball with an inhomogeneous mass distribution, but with its center of mass located at the ball's geometric center \cite{shen2008controllability}. A Chaplygin top is a ball with an inhomogeneous mass distribution, but with its center of mass not located at the ball's geometric center \cite{shen2008controllability}. Reference \cite{Ho2011_pII} does not distinguish between these two cases, calling a Chaplygin ball a ball with an inhomogeneous mass distribution, regardless of the location of its center of mass; as a special case of a Chaplygin ball, \cite{Ho2011_pII} calls a Chaplygin concentric sphere a ball with an inhomogeneous mass distribution with its center of mass coinciding with the ball's geometric center. Thus, the Chaplygin concentric sphere (used by \cite{Ho2011_pII}) and the Chaplygin sphere (used by \cite{shen2008controllability}) are different terms for the same mechanical system. Note that a ball with a homogeneous mass distribution (in a uniform gravitational field) necessarily has its center of mass at the ball's geometric center, and is therefore not very interesting. In this paper, these terminologies are not used, rather the mechanical system is referred to simply as a ball or a rolling ball, regardless of its mass distribution (homogeneous vs inhomogeneous) and regardless of the location of its center of mass (at the ball's geometric center vs not at the ball's geometric center). 
}. 
In this paper, the motion of the rolling ball is investigated assuming dynamic internal structure. The dynamics of the rolling ball with static internal structure was first solved analytically by Chaplygin for the cylindrically symmetric rolling ball, i.e. a ball such that the line joining the ball's center of mass and geometric center forms an axis of symmetry, in 1897 \cite{chaplygin2002motion} and for the Chaplygin sphere in 1903 \cite{chaplygin2002ball}, though dynamical properties of the cylindrically symmetric rolling ball were previously investigated by Routh \cite{routh1884advanced} and Jellet \cite{jellett1872treatise}. More recently, \cite{borisov2013problem} provides a detailed analysis of the trajectory of the Chaplygin sphere's contact point, and it has been shown that the dynamics of the Chaplygin top exhibit a strange attractor \cite{borisov2016spiral} and the phenomenon of reversal \cite{borisov2014reversal}. The dynamics of the rolling ball with dynamic internal structure is also an active topic in the nonholonomic mechanics literature \cite{das2001design,mojabi2002introducing,shen2008controllability,borisov2012control,bolotin2012problem,gajbhiye2016geometric,GaBa:RNC3457,kilin2015spherical,burkhardt2016reduced}. \rem{For example, the recent paper \cite{GaBa:RNC3457} explores \emph{local} controllability of a rolling ball with dynamic internal structure, in which an internal pendulum and yoke serve as the control mechanisms.}

\begin{figure}[h] 
	\centering
	\subfloat[A ball actuated by $3$ rotors, studied in \cite{borisov2012control,bolotin2012problem,gajbhiye2016geometric}, \copyright\ 2016 IFAC   \cite{gajbhiye2016geometric}.]{\includegraphics[width=0.3\textwidth]{ball_rotors}\label{fig_ball_rotors}}
	\hspace{5mm}
	\subfloat[A ball actuated by $6$ magnets, each in its own linear, solenoidal tube, studied in \cite{burkhardt2014energy,davoodi2014moball,asama2015design,davoodi2015moball,bowkett2016combined,burkhardt2016reduced}.]{
		\includegraphics[width=0.3\textwidth]{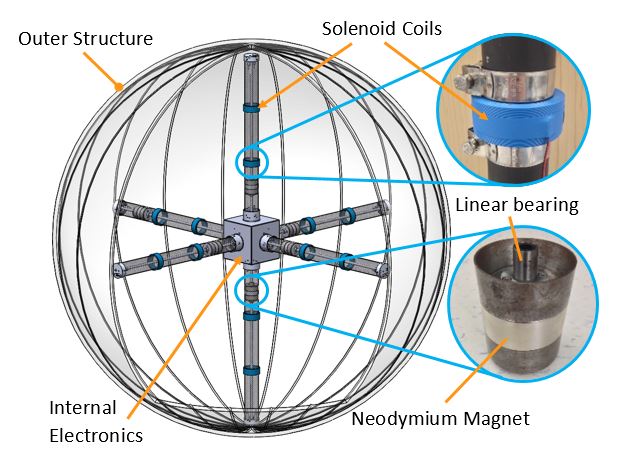}\label{fig_ball_magnets}}
	\hspace{5mm}
	\subfloat[A ball actuated by $3$ point masses, each on its own circular rail, studied in \cite{Putkaradze2018dynamicsP,putkaradze2018normalpub,RBOCnumerics}.]{
		\includegraphics[width=0.3\textwidth]{bsim2_control_masses_rails_png}\label{fig_intro_bsim2_control_masses_rails}}
	\hspace{5mm}
	\\
	\subfloat[A ball actuated by a gyroscopic pendulum, studied in \cite{kilin2015spherical}.]{
		\includegraphics[width=0.3\textwidth]{ball_gyro_pend_custom}\label{fig_ball_gyro_pend}}
	\hspace{5mm}
	\subfloat[A ball actuated by a spherical pendulum, studied in \cite{bolotin2013motion}.]{
		\includegraphics[width=0.3\textwidth]{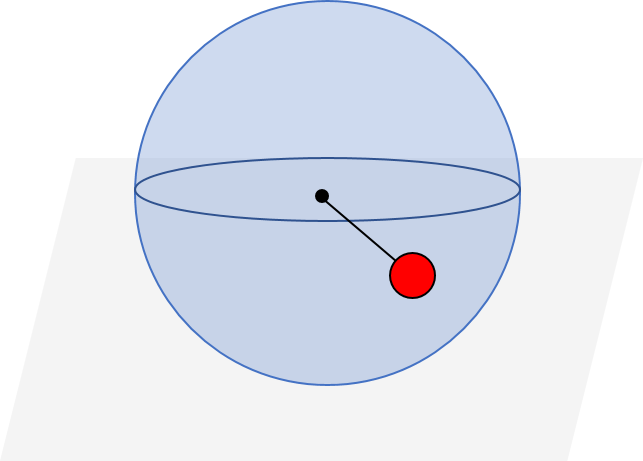}\label{fig_ball_sph_pend}}
	\hspace{5mm}
	\subfloat[Sphero has $4$ wheels wedged inside the spherical shell, but only the lower $2$ are spun by the motor \cite{Sphero_cut}.]{
		\includegraphics[width=0.3\textwidth]{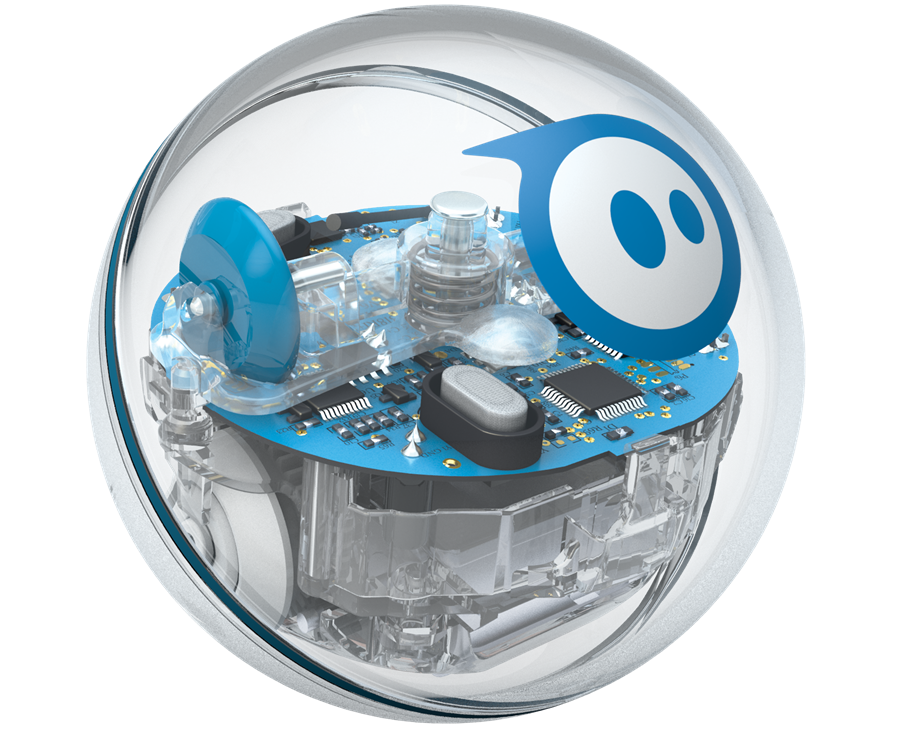}\label{fig_ball_sphero_internal}}	
	\caption{Different methods to actuate a rolling ball.} \label{fig_ball_actuation}
\end{figure}

Many methods have been proposed (and some realized) to actuate a rolling ball, such as illustrated in Figure~\ref{fig_ball_actuation}. References \cite{borisov2012control,bolotin2012problem,gajbhiye2016geometric} actuate the rolling ball by internal rotors such as shown in Figure~\ref{fig_ball_rotors}, while \cite{burkhardt2014energy,davoodi2014moball,asama2015design,davoodi2015moball,bowkett2016combined,burkhardt2016reduced} actuate the rolling ball via $6$ internal magnets, each of which slides inside its own linear, solenoidal tube, i.e. a straight tube embedded within a solenoid that generates a magnetic field along the tube's longitudinal axis as illustrated in Figure~\ref{fig_ball_magnets}. 
References  \cite{das2001design,mojabi2002introducing} study the locomotion and trajectory-tracking of a ball with masses moving along straight rails inside the ball, as well as practical realizations of such a device. In particular, \cite{mojabi2002introducing} actuates the rolling ball by internal masses which reciprocate along spokes. Reference \cite{shen2008controllability} actuates the rolling ball by a combination of internal rotors and sliders, \cite{kilin2015spherical} actuates the rolling ball by an internal gyroscopic pendulum as shown in Figure~\ref{fig_ball_gyro_pend}, \cite{bolotin2013motion,pivovarova2012stability,ivanova2013dynamics,ivanova2018controlled} actuate the rolling ball by an internal spherical pendulum as shown in Figure~\ref{fig_ball_sph_pend}, and \cite{GaBa:RNC3457} actuates the rolling ball by an internal pendulum and yoke. 
References \cite{karavaev2016nonholonomic,kilin2017experimental} study the theory of and experiment with a rolling ball actuated by an internal omniwheel platform, which is a more sophisticated version of Sphero's driving mechanism depicted in Figure~\ref{fig_ball_sphero_internal}. This paper considers a rolling ball actuated by internal point masses that move along arbitrarily-shaped rails fixed within the ball, such as depicted in Figure~\ref{fig_intro_bsim2_control_masses_rails}. Actuating the rolling ball by moving internal point masses along general rails has not been considered yet in the literature. References \cite{burkhardt2016reduced,das2001design,mojabi2002introducing,shen2008controllability} actuate the rolling ball by moving internal masses with inertias along linear trajectories (e.g. spokes or hollow tubes) in the ball's frame. 
\revisionT{EQ4}{Reference \cite{tomik2012design} designs a transfer function between the inputs and outputs of such a device and shows an example of a robot navigating a path within a building while being controlled by human input specifying the desired values of momenta in time. }  The very recent work \cite{ilin2017dynamics} derives and simulates the dynamics of a simplified model of a beaver ball, in which a point mass moves with constant angular velocity along a circular trajectory fixed inside the ball. Prior to this paper, controlling the motion of a nonholonomic mechanical system by moving internal point masses has been studied in \cite{osborne2005steering}, which investigates the controlled motion of the Chaplygin sleigh actuated by a single internal point mass. 

More generally,  the optimal control of other nonholonomic mechanical systems has been investigated recently; the continuous variable transmission and the \revision{R1Q6}{Chaplygin} sleigh are investigated in \cite{Bl-etal-2015}, while Suslov's problem is investigated in \cite{putkaradze2018constraint}.
In a comprehensive review of nonholonomic optimal control, \cite{Bloch2003} briefly discusses the controllability and optimal control (in the sense of the variational Pontryagin's minimum principle) of a rolling ball, where an external control force pushes the ball's geometric center. Reference \cite{grong2016submersions} uses the language of differential geometry to abstractly study the optimal control of a ball rolling on a horizontal surface, but does not specify a particular actuation mechanism. A recent survey \cite{borisov2017dynamical} provides an overview of the theory of mechanical systems having different types of nonintegrable velocity constraints, such as vakonomic, which is useful for control theory, and nonholonomic. Reference~\cite{shen2008controllability} investigates the controllability of a rolling ball actuated by internal rotors and sliders, while \cite{GaBa:RNC3457} studies the controllability of a rolling ball actuated by an internal pendulum and yoke. For a rolling ball actuated by internal rotors, reference~\cite{borisov2012control} performs an analysis of controllability and trajectory tracking control, \cite{bolotin2012problem} analyzes the optimal control (in the sense of the variational Pontryagin's minimum principle) for some very specific  performance indexes, and \cite{gajbhiye2016geometric} investigates orientation and contact point trajectory tracking control by constructing feedback control laws.
\revisionT{EQ1\\EQ2\\EQ3}{ Appendix~\ref{app_intrinsic} describes several methods for kinematic control of the rolling ball. The work in that field started with the celebrated papers by Jurdjevic providing an exact solution to the optimal control problem \cite{Jurdjevic1993,Jurdjevic1999a,Jurdjevic1999b,ohsawa2020geometric} and subsequent application of geometric theory to the trajectory planning of bodies rolling with and  without twisting at the contact point \cite{AlChLo2010}. Other methods of kinematic control have also been considered. For example, paper \cite{mukherjee2002motion} proposes two ways of controlling a rolling ball robot, one using a trajectory formed by a combination of circular arcs and straight line segments and another using the parallel transport theorem. The paper \cite{das2004exponential} derives a nonsmooth control that provides  exponential stabilization of the rolling ball and motion along straight lines and circular arcs. Finally, \cite{das2006reconfiguration} provides a geometric viewpoint of the problem using a mapping of a three-dimensional kinematic problem to a planar problem, namely a problem of wrapping and unwrapping a taut rope on a planar curve.   }
Closer to the topic of this paper, reference \cite{mojabi2002introducing} investigates the trajectory tracking control of a rolling ball, actuated by internal masses which reciprocate along spokes, by minimizing the ball's energy expenditure during the motion by discretizing the trajectory and finding a global minimum of a certain function, i.e. by using the direct method. Reference~\cite{kilin2015spherical} investigates the trajectory tracking control of a rolling ball actuated by an internal gyroscopic pendulum by piecing together special gaits, each of whose motion is determined analytically. To contrast with previous work, this paper investigates the optimal control (in the sense of the variational Pontryagin's minimum principle) of a rolling ball actuated by internal point masses for very general  performance indexes and mass trajectories. 

\rem{
\begin{figure}[h] 
	\centering
	\subfloat[A ball actuated by $3$ rotors, studied in \cite{borisov2012control,bolotin2012problem,gajbhiye2016geometric}, \copyright\ 2016 IFAC   \cite{gajbhiye2016geometric}.]{\includegraphics[width=0.45\textwidth]{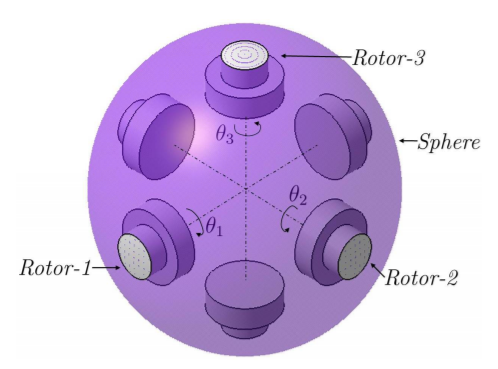}\label{fig_ball_rotors}}
	\hspace{5mm}
	\subfloat[A ball actuated by $3$ point masses, each on its own circular control rail, studied in this paper.]{
		\includegraphics[width=0.45\textwidth]{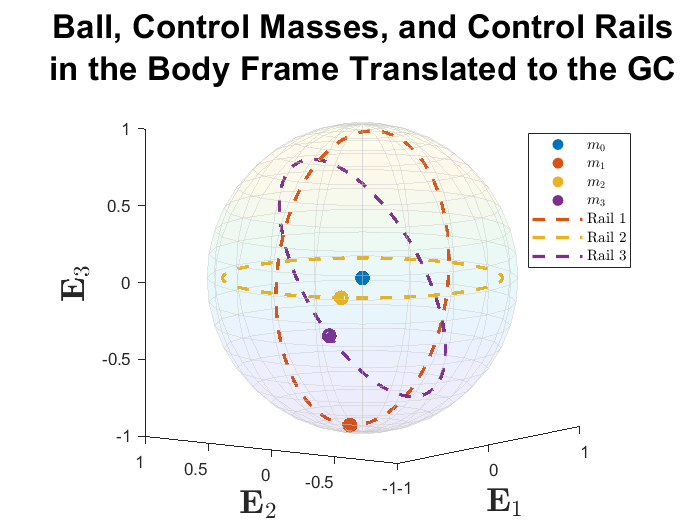}\label{fig_intro_bsim2_control_masses_rails}}
	\\
	\subfloat[A ball actuated by a gyroscopic pendulum, studied in \cite{kilin2015spherical}.]{
		\includegraphics[width=0.45\textwidth]{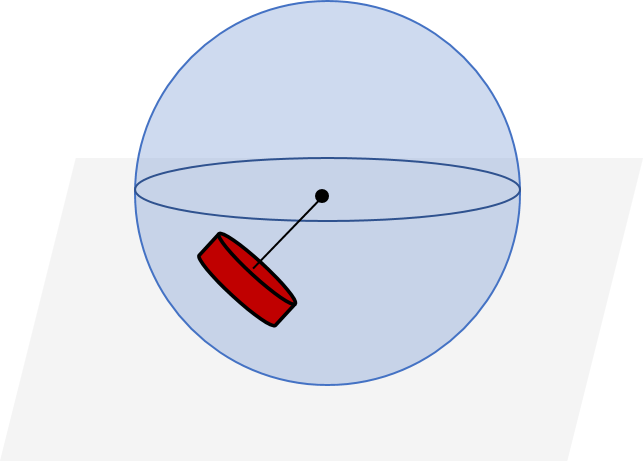}\label{fig_ball_gyro_pend}}
	\hspace{5mm}
	\subfloat[Sphero has $4$ wheels wedged inside the spherical shell, but only the lower $2$ are spun by the motor \cite{Sphero_cut}.]{
		\includegraphics[width=0.45\textwidth]{sphero_cutaway}\label{fig_ball_sphero_internal}}	
	\caption{Different methods to actuate a rolling ball.} \label{fig_ball_actuation}
\end{figure} }

\subsection{Contributions and Outline}
The key contribution of this paper is the construction of the controlled equations of motion for a rolling ball actuated by internal point masses that move along arbitrarily-shaped rails fixed within the ball, a solution of which obeys the uncontrolled equations of motion, satisfies prescribed initial and final conditions, and minimizes a prescribed performance index. The paper is organized as follows.  Section~\ref{sec_uncontrolled} discusses the specific type of rolling disk and ball considered, defines coordinate systems and notation used to describe this rolling disk and ball, and presents the uncontrolled equations of motion for this rolling disk and ball derived earlier in \cite{Putkaradze2018dynamicsP,putkaradze2018normalpub}. Section~\ref{sec_ball_controlled} formulates the regular Hamiltonian and endpoint function required to construct the controlled equations of motion for the rolling disk and ball. Section~\ref{sec_conclusions} summarizes the results of the paper and discusses topics for future work. Appendix~\ref{app_ccem} derives formulas needed to construct the controlled equations of motion for the rolling disk and ball. Appendix~\ref{app_intrinsic} derives the uncontrolled and controlled equations of motion for the rolling ball in intrinsic coordinates.

\section{Uncontrolled Equations of Motion \label{sec_uncontrolled} } 
This section presents the notation, coordinate systems, and assumptions used to describe the rolling ball and its special case, the rolling disk. The uncontrolled equations of motion for these mechanical systems are recalled from \cite{Putkaradze2018dynamicsP,putkaradze2018normalpub}. We refer  the interested reader to those papers for the  detailed derivations.  

\subsection{Uncontrolled Equations of Motion for the Rolling Ball} \label{ssec_ball_uncontrolled}
\begin{figure}[h]
	\centering
	\includegraphics[width=0.5\linewidth]{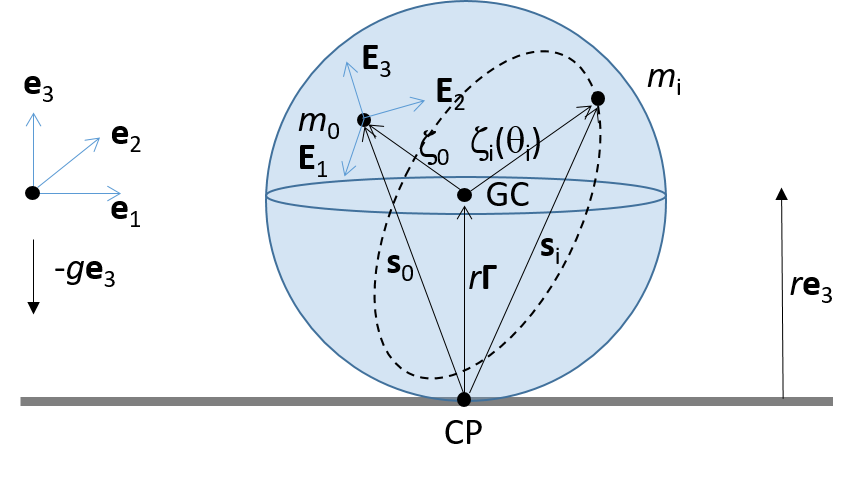}
	\caption{A ball of radius $r$ and mass $m_0$ rolls without slipping on a horizontal surface in the presence of a uniform gravitational field of magnitude $g$. The ball's geometric center, center of mass, and contact point with the horizontal surface are denoted by GC, $m_0$, and CP, respectively. The spatial frame has origin located at height $r$ above the horizontal surface and orthonormal axes $\mathbf{e}_1$, $\mathbf{e}_2$, and $\mathbf{e}_3$. The body frame has origin located at the ball's center of mass (denoted by $m_0$) and orthonormal axes $\mathbf{E}_1$, $\mathbf{E}_2$, and $\mathbf{E}_3$. The ball's motion is actuated by $n$ point masses, each of mass $m_i$, $1 \le i \le n$, and each moving along its own rail fixed inside the ball. The $i^\mathrm{th}$ rail is depicted here by the dashed hoop. The trajectory of the $i^\mathrm{th}$ rail, with respect to the body frame translated to the GC, is denoted by $\bzeta_i$ and is parameterized by $\theta_i$. All vectors inside the ball are expressed with respect to the body frame, while all vectors outside the ball are expressed with respect to the spatial frame. }
	\label{fig:detailed_1dparam_rolling_ball}
\end{figure}
\paragraph{Physical assumptions on the ball's structure and motion.}
Consider a rigid ball of radius $r$ containing some static internal structure as well as $n$ point masses, \revision{R1Q7}{where either $n$ is a positive integer denoting the number of moving masses or $n=0$ if no moving masses are used and the structure of the ball is static}. 
This ball rolls without slipping on a horizontal  surface in the presence of a uniform gravitational field of magnitude $g$, as illustrated in Figure~\ref{fig:detailed_1dparam_rolling_ball}. The ball with its static internal structure has mass $m_0$ and the $i^\mathrm{th}$ point mass has mass $m_i$ for $1 \le i \le n$. Let $M = \sum_{i=0}^n m_i$ denote the mass of the total system. The total mechanical system consisting of the ball with its static internal structure and the $n$ point masses is referred to as the ball or the rolling ball, the ball with its static internal structure but without the $n$  point masses may also be referred to as $m_0$, and the $i^\mathrm{th}$  point mass may also be referred to as $m_i$ for $1 \le i \le n$.

\paragraph{Coordinate systems: spatial frame, body  frame, and notation.}
Two coordinate systems, or frames of reference, will be used to describe the motion of the rolling ball, an inertial spatial coordinate system and a body coordinate system in which each particle within the ball is always fixed. For brevity, the spatial coordinate system will be referred to as the spatial frame and the body coordinate system will be referred to as the body frame. These two frames are depicted in Figure~\ref{fig:detailed_1dparam_rolling_ball}. The spatial frame has orthonormal axes $\mathbf{e}_1$, $\mathbf{e}_2$, $\mathbf{e}_3$, such that the $\mathbf{e}_1$-$\mathbf{e}_2$ plane is parallel to the horizontal surface and passes through the ball's geometric center (i.e. the $\mathbf{e}_1$-$\mathbf{e}_2$ plane is a height $r$ above the horizontal surface), such that $\mathbf{e}_3$ is vertical (i.e. $\mathbf{e}_3$ is perpendicular to the horizontal surface) and points ``upward" and away from the horizontal surface, and such that $\left(\mathbf{e}_1, \mathbf{e}_2, \mathbf{e}_3 \right)$ forms a right-handed coordinate system. For simplicity, the spatial frame axes are chosen to be
\begin{equation}
\mathbf{e}_1 = \begin{bmatrix} 1 & 0 & 0 \end{bmatrix}^\mathsf{T}, \quad \mathbf{e}_2 = \begin{bmatrix} 0 & 1 & 0 \end{bmatrix}^\mathsf{T}, \quad \mathrm{and} \quad \mathbf{e}_3 = \begin{bmatrix} 0 & 0 & 1 \end{bmatrix}^\mathsf{T}.
\end{equation}
The acceleration due to gravity in the uniform gravitational field is $\mathfrak{g} = -g \mathbf{e}_3  = \begin{bmatrix} 0 & 0 & -g  \end{bmatrix}^\mathsf{T}$ in the spatial frame.

The body frame's origin is chosen to coincide with the position of $m_0$'s center of mass. The body frame has orthonormal axes $\mathbf{E}_1$, $\mathbf{E}_2$, and $\mathbf{E}_3$, chosen to coincide with $m_0$'s principal axes, in which $m_0$'s inertia tensor $\inertia$ is diagonal, with corresponding \revision{R1Q9}{principal} moments of inertia $d_1$, $d_2$, and $d_3$. That is, in this body frame the inertia tensor is the diagonal matrix $\inertia = \diag \left( \begin{bmatrix} d_1 & d_2 & d_3 \end{bmatrix} \right)$.
Moreover, $\mathbf{E}_1$, $\mathbf{E}_2$, and $\mathbf{E}_3$ are chosen so that $\left(\mathbf{E}_1, \mathbf{E}_2, \mathbf{E}_3 \right)$ forms a right-handed coordinate system. For simplicity, the body frame axes are chosen to be
\begin{equation}
\mathbf{E}_1 = \begin{bmatrix} 1 & 0 & 0 \end{bmatrix}^\mathsf{T}, \quad \mathbf{E}_2 = \begin{bmatrix} 0 & 1 & 0 \end{bmatrix}^\mathsf{T}, \quad \mathrm{and} \quad \mathbf{E}_3 = \begin{bmatrix} 0 & 0 & 1 \end{bmatrix}^\mathsf{T}.
\end{equation}
In the spatial frame, the body frame is the moving frame $\left(\Lambda \left(t\right) \mathbf{E}_1, \Lambda \left(t\right) \mathbf{E}_2, \Lambda \left(t\right) \mathbf{E}_3  \right)$, where $\Lambda \left(t\right) \in SO(3)$ defines the orientation (or attitude) of the ball at time $t$ relative to its reference configuration, for example at some initial time.   

\paragraph{Assumptions on the motion of the masses inside the ball.}
For $1 \le i \le n$, it is assumed that $m_i$ moves along its own 1-d rail. It is further assumed that the $i^\mathrm{th}$  rail is parameterized by a 1-d parameter $\theta_i$, so that the trajectory $\bzeta_i$ of the $i^\mathrm{th}$  rail, in the body frame translated to the ball's geometric center, as a function of $\theta_i$ is $\bzeta_i(\theta_i)$. Refer to Figure~\ref{fig:detailed_1dparam_rolling_ball} for an illustration.
\revisionS{R1Q3}{We further assume that $\theta_i$ is a prescribed function of time. In general, for internal masses actuated by motors with given torques, the motion of the masses as a function of time cannot be prescribed, but instead must be solved for in conjunction with the ball's motion \cite{BaSvYa2018}. 
We envision a different driving mechanism based on stepper motors which are rigidly attached to the internal frame of the rolling ball. Unlike a regular electric motor which generates a given torque based on the input voltage/current, a stepper motor is a device which rotates the motor's shaft a given amount measured in a discrete number of steps, where each step is typically 1-2 degrees depending on the motor's design, with a typical error of 1-2\% of the step angle. Thus, as long as appropriately powerful stepper motors are used, it is reasonable to assume the prescribed motion of the internal masses. We refer the reader to  \cite{putkaradze2018normalpub} for details of this discussion.}

Therefore, the body frame vector from the ball's geometric center to $m_i$'s center of mass is denoted by $\bzeta_i(\theta_i (t))$. Since $m_0$ is stationary in the body frame and to be consistent with the positional notation for $m_i$ for $1 \le i \le n$, $\bzeta_0 \equiv \bzeta_0(\theta_0) \equiv \bzeta_0(\theta_0 (t))$ is the constant  (time-independent) vector from the ball's geometric center to $m_0$'s center of mass for any scalar-valued, time-varying function $\theta_0(t)$. In addition, suppose a time-varying external force $\mathbf{F}_\mathrm{e}(t)$ acts at the ball's geometric center. \revisionS{R1Q6}{$\mathbf{F}_\mathrm{e}(t)$ does not involve the  static  friction  induced by  the surface to enforce   the no-slip constraint. Instead, it involves forces due to other environmental factors such as  air resistance (i.e. drag) and wind force. We assume that $\mathbf{F}_\mathrm{e}(t)$ is a prescribed function of time. In reality, $\mathbf{F}_\mathrm{e}$ may depend  on the ball's position, orientation, angular velocity, and linear velocity. While this dependence can be introduced in the model, it is hard to predict the analytical form of that dependence. Thus, we shall only consider the simplest case where $\mathbf{F}_\mathrm{e}(t)$ depends only on time. } 

Let $\mathbf{z}_i(t)$ denote the position of $m_i$'s center of mass in the spatial frame so that the position of $m_i$'s center of mass in the spatial frame is $\mathbf{z}_i(t)=\mathbf{z}_0(t)+\Lambda(t) \left[\bzeta_i(\theta(t))-\bzeta_0\right]$. In general, a particle with position $\mathbf{w}(t)$ in the body frame has position $\mathbf{z}(t) = \mathbf{z}_0(t)+\Lambda(t) \mathbf{w}(t)$ in the spatial frame and has position $\mathbf{w}(t)+\bzeta_0$ in the body frame translated to the ball's geometric center.

For conciseness, the ball's geometric center is often denoted GC, $m_0$'s center of mass is often denoted CM, and the ball's contact point with the surface is often denoted CP. The GC is located at $\mathbf{z}_\mathrm{GC}(t) = \mathbf{z}_0(t)-\Lambda(t) \bzeta_0$ in the spatial frame, at $-\bzeta_0$ in the body frame, and at $\mathbf{0} = \begin{bmatrix} 0 & 0 & 0 \end{bmatrix}^\mathsf{T}$ in the body frame translated to the GC. The CM is located at $\mathbf{z}_0(t)$ in the spatial frame, at $\mathbf{0}$ in the body frame, and at $\bzeta_0$ in the body frame translated to the GC. The CP is located at $\mathbf{z}_\mathrm{CP}(t) = \mathbf{z}_0(t)-\Lambda(t) \left[r\bGam(t)+\bzeta_0 \right]$ in the spatial frame, at $-\left[r\bGam(t)+\bzeta_0 \right]$ in the body frame, and at $-r\bGam(t)$ in the body frame translated to the GC, where $\bGamma(t) \equiv \Lambda^{-1}(t) \mathbf{e}_3$. Since the third spatial coordinate of the ball's GC is always $0$ and of the ball's CP is always $-r$, only the first two spatial coordinates of the ball's GC and CP, denoted by $\bz(t)$, are needed to determine the spatial location of the ball's GC and CP. 

For succintness, the explicit time dependence of variables is often dropped. That is, the orientation of the ball at time $t$ is denoted simply $\Lambda$ rather than $\Lambda(t)$, the position of $m_i$'s center of mass in the spatial frame at time $t$ is denoted $\mathbf{z}_i$ rather than $\mathbf{z}_i(t)$,  the position of $m_i$'s center of mass in the body frame translated to the GC at time $t$ is denoted $\bzeta_i$ or $\bzeta_i(\theta_i)$ rather than $\bzeta_i(\theta_i(t))$, the spatial $\mathbf{e}_1$-$\mathbf{e}_2$ position of the ball's GC at time $t$ is denoted $\bz$ rather than $\bz(t)$, and the external force is denoted $\mathbf{F}_\mathrm{e}$ rather than $\mathbf{F}_\mathrm{e}(t)$. 

\paragraph{Uncontrolled equations of motion.} 
\revisionS{R2Q2}{The ball's Lagrangian is the kinetic energy of translation and rotation  minus the gravitational potential energy. The uncontrolled equations of motion are derived by using  Lagrange-d'Alembert's principle for nonholonomic systems  or by using Newton's laws. These equations are computed in \cite{Putkaradze2018dynamicsP,putkaradze2018normalpub}  to have the following form: }
\begin{equation}
\begin{split} \label{uncon_ball_eqns_explicit_1d}
\dot \bOm &= \left[\sum_{i=0}^n m_i \widehat{\mathbf{s}_i}^2  -\inertia \right]^{-1}  \Bigg[\bOm \times \inertia \bOm+r \tilde \bGamma \times \bGamma 
\\ &\hphantom{=} + \sum_{i=0}^n m_i \mathbf{s}_i \times  \left\{ g \bGamma+ \bOm \times \left(\bOm \times \bzeta_i +2 \dot \theta_i \bzeta_i^{\prime} \right) + \dot \theta_i^2 \bzeta_i^{\dprime} + \ddot \theta_i  \bzeta_i^{\prime} \right\}  \Bigg], \\
\dot \Lambda &= \Lambda \widehat{\bOm}, \\
\dot \bz &= \left( \Lambda \bOm \times r \mathbf{e}_3 \right)_{12},
\end{split}
\end{equation}
where $\mathbf{s}_i \equiv r \bGamma +\bzeta_i$ is the body frame vector from the CP to $m_i$ for $0\le i\le n$, $\bOm \equiv \left( \Lambda^{-1} \dot \Lambda \right)^\vee$ is the ball's body angular velocity, $\bGamma \equiv \Lambda^{-1} \mathbf{e}_3$ is the spatial unit normal expressed in the body frame, and $\tilde \bGamma \equiv \Lambda^{-1} \mathbf{F}_\mathrm{e}$ is the external force expressed in the body frame. For $\mathbf{v} = \begin{bmatrix} v_1 & v_2 & v_3  \end{bmatrix}^\mathsf{T} \in \mathbb{R}^3$, $\widehat{\mathbf{v}}^2=\widehat{\mathbf{v}}\widehat{\mathbf{v}}$ is the symmetric matrix given by
\begin{equation}
\widehat{\mathbf{v}}^2 = \begin{bmatrix}
-(v_2^2+v_3^2) & v_1 v_2 & v_1 v_3 \\
v_1 v_2 & -(v_1^2+v_3^2)  & v_2 v_3 \\
v_1 v_3 & v_2 v_3 & -(v_1^2+v_2^2) 
\end{bmatrix}
\end{equation}
and $\mathbf{v}_{12}$ is the projected vector consisting of the first two components of $\mathbf{v}$ so that
\begin{equation}
\mathbf{v}_{12} = \begin{bmatrix}v_1 & v_2 \end{bmatrix}^\mathsf{T} \in \mathbb{R}^2.
\end{equation}
 
\subsection{Uncontrolled Equation of Motion for the Rolling Disk} \label{ssec_disk_uncontrolled}
\begin{figure}[h]
	\centering
	\includegraphics[width=0.6\linewidth]{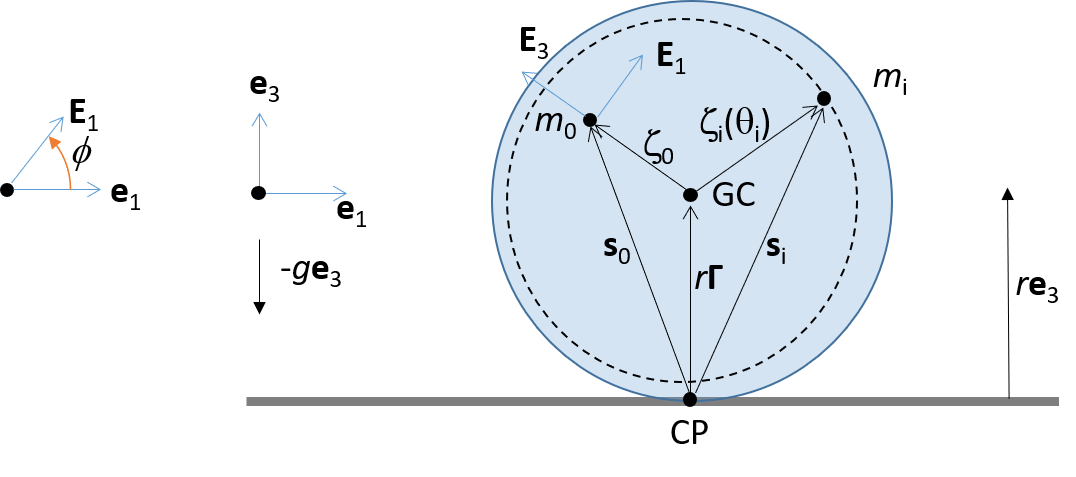}
	\caption{A disk of radius $r$ and mass $m_0$ rolls without slipping in the $\mathbf{e}_1$-$\mathbf{e}_3$ plane. $\mathbf{e}_2$ and $\mathbf{E}_2$ are directed into the page and are omitted from the figure. The disk's center of mass is denoted by $m_0$. The disk's motion is actuated by $n$ point masses, each of mass $m_i$, $1 \le i \le n$, and each moving along its own rail fixed inside the disk. The point mass depicted here by $m_i$ moves along a circular hoop in the disk that is not centered on the disk's geometric center (GC). The disk's orientation is determined by $\phi$, the angle measured counterclockwise from $\mathbf{e}_1$ to $\mathbf{E}_1$.  }
	\label{fig:detailed_rolling_disk}
\end{figure}
Now suppose that  the ball's inertia is such that one of  the ball's principal axes, say the one labeled $\mathbf{E}_2$, is orthogonal to the plane containing the GC and CM. 
Also assume that all the  point masses move along 1-d  rails which lie in the plane containing the GC and CM. Moreover, suppose that the ball is oriented initially so that the plane containing the GC and CM coincides with the $\mathbf{e}_1$-$\mathbf{e}_3$ plane and that the external force $\mathbf{F}_\mathrm{e}$ acts in the $\mathbf{e}_1$-$\mathbf{e}_3$ plane. Then for all time, the ball will remain oriented so that the plane containing the GC and CM coincides with the $\mathbf{e}_1$-$\mathbf{e}_3$ plane and the ball will only move in the $\mathbf{e}_1$-$\mathbf{e}_3$ plane, with the ball's rotation axis always parallel to $\mathbf{e}_2$. Note that the dynamics of this system are equivalent to that of the Chaplygin disk \cite{Ho2011_pII}, equipped with point masses, rolling in the $\mathbf{e}_1$-$\mathbf{e}_3$ plane, and where the Chaplygin disk (minus the  point masses) has polar moment of inertia $d_2$. Therefore, henceforth, this particular ball with this special inertia, orientation,  and placement of the rails and point masses, may be referred to as the disk or the rolling disk. Figure~\ref{fig:detailed_rolling_disk} depicts the rolling disk. Let $\phi$ denote the angle between $\mathbf{e}_1$ and $\mathbf{E}_1$, measured counterclockwise from $\mathbf{e}_1$ to $\mathbf{E}_1$. Thus, if $\dot \phi > 0$, the disk rolls in the $-\mathbf{e}_1$ direction and $\bOm$ has the same direction as $-\mathbf{e}_2$, and  if $\dot \phi < 0$, the disk rolls in the $\mathbf{e}_1$ direction and $\bOm$ has the same direction as $\mathbf{e}_2$. 

As shown in \cite{Putkaradze2018dynamicsP,putkaradze2018normalpub}, the uncontrolled equation of motion for this rolling disk is
\begin{equation} \label{eqmo_chap_disk_4}
\ddot \phi = \frac{ -r F_{\mathrm{e},1}+ \sum_{i=0}^n m_i K_i }{d_2+\sum_{i=0}^n m_i \left[\left( r \sin \phi + \zeta_{i,1} \right)^2+\left( r \cos \phi+ \zeta_{i,3} \right)^2 \right]} \equiv \kappa\left(t,\btheta,\dot \btheta,\phi,\dot \phi,\ddot \btheta\right),
\end{equation}
where 
\begin{equation} \label{eq_K_i}
\begin{split}
K_i &\equiv  \left(g+ r {\dot \phi}^2 \right) \left(\zeta_{i,3} \sin \phi - \zeta_{i,1} \cos \phi  \right)+
\left(r \cos \phi + \zeta_{i,3} \right) \left(- 2 \dot \phi {\dot \theta}_i \zeta_{i,3}^{\prime} + {\dot \theta}_i^2 \zeta_{i,1}^{\dprime} + {\ddot \theta}_i \zeta_{i,1}^{\prime} \right)\\
&\hphantom{\equiv}   - \left(r \sin \phi + \zeta_{i,1} \right) \left( 2 \dot \phi {\dot \theta}_i \zeta_{i,1}^{\prime}+ {\dot \theta}_i^2 \zeta_{i,3}^{\dprime} + {\ddot \theta}_i \zeta_{i,3}^{\prime} \right). 
\end{split}
\end{equation}
In \eqref{eqmo_chap_disk_4}, $\kappa$ is a function that depends on time ($t$) through the possibly time-varying external force $F_{\mathrm{e},1}(t)$, on the  point mass parameterized positions ($\btheta$), velocities ($\dot \btheta$), and accelerations ($\ddot \btheta$), and on the disk's orientation angle ($\phi$) and its time-derivative ($\dot \phi$). The spatial $\mathbf{e}_1$ position $z$ of the disk's GC is given by 
\begin{equation} \label{eq_disk_GC}
z = z_a-r \left(\phi-\phi_a \right),
\end{equation}
where $z_a$ is the spatial $\mathbf{e}_1$ position of the disk's GC at time $t=a$ and $\phi_a$ is the disk's angle at time $t=a$.

\section{Controlled Equations of Motion}
\label{sec_ball_controlled}

This section uses the results of the previous section to derive the controlled equations of motion for the rolling disk and ball.
Subsection~\ref{ssec_disk_controlled} exploits the uncontrolled equation of motion \eqref{eqmo_chap_disk_4} to derive the controlled equations of motion for a rolling disk actuated by internal point masses that move along arbitrarily-shaped rails fixed within the disk. Subsection~\ref{ssec_ball_controlled} exploits the uncontrolled equations of motion \eqref{uncon_ball_eqns_explicit_1d} to derive the controlled equations of motion for a rolling ball actuated by internal point masses that move along arbitrarily-shaped rails fixed within the ball.

\subsection{Controlled Equations of Motion for the Rolling Disk} \label{ssec_disk_controlled}

\paragraph{Notation, assumptions, and the motion of the center of mass.}
Before  treating the problem of  the rolling ball, the controlled equations of motion are first developed for the rolling disk, which is a simpler mechanical system. Let $z$ and $\dot z$ denote the spatial $\mathbf{e}_1$ position and velocity, respectively, of the disk's GC, and recall that $\btheta = \begin{bmatrix} \theta_1 & \theta_2 & \dots & \theta_n \end{bmatrix}^\mathsf{T}$ denotes the vector of the control mass parameterizations. If the disk's GC is at initial spatial $\mathbf{e}_1$ position $z_a$ and if the disk's initial orientation is $\phi_a$ at initial time $t=a$, note that the spatial $\mathbf{e}_1$ position and velocity of the disk's GC are $z = z_a - r \left(\phi-\phi_a \right)$ and $\dot z = -r \dot \phi$, respectively, due to the sign convention adopted for $\phi$ in Subsection~\ref{ssec_disk_uncontrolled}. The rotation matrix that maps the body to spatial frame at time $t$ is a function of $\phi(t)$ and is given by
\begin{equation}
\tilde \Lambda (\phi(t)) = \Lambda(t) = \begin{bmatrix} \cos \phi(t) & 0 & - \sin \phi(t) \\ 0 & 1 & 0 \\ \sin \phi(t) & 0 & \cos \phi(t) \end{bmatrix}.
\end{equation}

\revision{R1Q2}{It is desired to roll the disk from some initial configuration at a prescribed or free initial time $a$ to some final configuration at a prescribed or free final time $b$, without moving the control masses too rapidly along their control rails, while the disk's GC tracks a prescribed or minimum energy path or while minimizing the maneuver's duration $b-a$. } 
\revision{R1Q11}{Concretely, at the prescribed or free initial time $a$,  the positions of the control mass parameterizations are prescribed to be $\btheta(a) = \btheta_a$, the velocities of the control mass parameterizations are prescribed to be ${\dot \btheta}(a) = {\dot \btheta}_a$, and the spatial $\mathbf{e}_1$ velocity of the disk's GC is prescribed to be ${\dot z}(a) = -r \dot \phi(a) = {\dot z}_a$ (which is equivalent to prescribing the rate of change of the disk's orientation to be $\dot \phi(a)=-\frac{{\dot z}_a}{r}$).}

Furthermore, at the prescribed or free final time $b$, some components (determined by the projection operator $\bPi$) \revision{R1Q12}{of the position} of the disk's center of mass expressed in the spatial frame translated to the GC are prescribed to be 
\begin{equation}
\bPi \left(\tilde \Lambda\left(\phi(b)\right) \left[ \frac{1}{M} \sum_{i=0}^n m_i \bzeta_i\left(\theta_i(b)\right) \right] \right) = \bDelta_b, 
\end{equation}
the velocities of the control mass parameterizations are prescribed to be ${\dot \btheta}(b) = {\dot \btheta}_b$, the spatial $\mathbf{e}_1$ position of the disk's GC is prescribed to be $z(b) =  z_a - r \left(\phi(b)-\phi_a\right) = z_b$ (which is equivalent to prescribing the disk's orientation to be $\phi(b)=\phi_a-\frac{z_b-z_a}{r}$), and the spatial $\mathbf{e}_1$ velocity of the disk's GC is prescribed to be ${\dot z}(b) = -r \dot \phi(b) = {\dot z}_b$ (which is equivalent to prescribing the rate of change of the disk's orientation to be $\dot \phi(b)=-\frac{{\dot z}_b}{r}$).

For example, if it is desired to start and stop the disk at rest, then $\bPi$ is \revision{R1Q13}{the} projection onto the first component, $\bDelta_b=\Delta_b = 0$, $\btheta_a$ and $\phi_a$ are such that 
\begin{equation} \label{eq_disk_stab_init}
\bPi \left(\tilde \Lambda\left(\phi_a\right) \left[ \frac{1}{M} \sum_{i=0}^n m_i \bzeta_i\left(\theta_{a,i}\right) \right] \right) = 0, 
\end{equation} 
$ {\dot \btheta}_a = \mathbf{0}$, ${\dot z}_a = 0$, $\btheta(b)$ and $\phi(b)$ are such that 
\begin{equation} \label{eq_disk_stab_fin}
\bPi \left(\tilde \Lambda\left(\phi(b)\right) \left[ \frac{1}{M} \sum_{i=0}^n m_i \bzeta_i\left(\theta_i(b)\right) \right] \right) = 0, 
\end{equation}
${\dot \btheta}_b = \mathbf{0}$, and ${\dot z}_b = 0$. With this choice of $\bPi$, \eqref{eq_disk_stab_init} and \eqref{eq_disk_stab_fin} mean that the CM in the spatial frame translated to the GC is above or below the GC at the initial and final times. \revisionS{R1Q9}{Note that \eqref{eq_disk_stab_init} and \eqref{eq_disk_stab_fin} impose restrictions on the initial and final positions of the internal masses. } 

\paragraph{State variables, control variables, and dynamics.}
The system state $\bx$ and control $\bu$ are
\begin{equation}
\bx \equiv \begin{bmatrix} \btheta \\ \dot \btheta \\ \phi \\ \dot \phi  \end{bmatrix} \quad \mathrm{and} \quad \bu \equiv \ddot \btheta,
\end{equation}
where $\btheta, \dot \btheta, \ddot \btheta \in \mathbb{R}^n$ and $\phi, \dot \phi  \in \mathbb{R}$. The system dynamics defined for $a \le t \le b$ are
\begin{equation} \label{eq_disk_dynamics}
\dot {\bx} = \begin{bmatrix} \dot \btheta \\ \ddot \btheta \\ \dot \phi \\ \ddot \phi  \end{bmatrix}  = \mathbf{f}\left(t,\bx,\bu,\mu\right) \equiv \begin{bmatrix} \dot \btheta \\ \bu  \\ \dot \phi \\ \kappa\left(t,\bx,\bu\right)  \end{bmatrix},
\end{equation}
\revision{R1Q4}{where $\mu$ is a scalar continuation parameter and} $\kappa\left(t,\bx,\bu \right)$ is given by the right-hand side of \eqref{eqmo_chap_disk_4}. 

\paragraph{Initial and final conditions.}
The prescribed initial conditions at time $t=a$ are
\begin{equation} \label{eq_disk_initial_conds}
\bsigma\left(a,\bx(a),\mu\right) \equiv \begin{bmatrix} \btheta(a) - \btheta_a \\ \dot \btheta(a) - {\dot \btheta}_a \\ \phi(a)-\phi_a \\ -r \dot \phi(a) - {\dot z}_a  \end{bmatrix} = \mathbf{0},
\end{equation}
and the prescribed final conditions at time $t=b$ are
\begin{equation} \label{eq_disk_final_conds}
\bpsi\left(b,\bx(b),\mu\right) \equiv \begin{bmatrix} \bPi \left(\tilde \Lambda\left(\phi(b)\right) \left[ \frac{1}{M} \sum_{i=0}^n m_i \bzeta_i\left(\theta_i(b)\right) \right] \right) - \bDelta_b \\ \dot \btheta(b) - {\dot \btheta}_b  \\ z_a - r \left(\phi(b)-\phi_a\right) - z_b \\ -r \dot \phi(b) - {\dot z}_b  \end{bmatrix} = \mathbf{0}.
\end{equation} 

\paragraph{Cost functions and performance index.}
Consider the endpoint and integrand cost functions
\begin{equation} \label{eq_disk_endpoint_cost}
p\left(a,\bx(a),b,\bx(b),\mu\right) \equiv \frac{\upsilon_a}{2}\left(a-a_\mathrm{e} \right)^2+\frac{\upsilon_b}{2}\left(b-b_\mathrm{e} \right)^2
\end{equation}
and
\begin{equation} \label{eq_disk_integrand_cost}
L\left(t,\bx,\bu,\mu\right) \equiv \frac{\alpha}{2} \left(z_a - r \left(\phi-\phi_a\right) - z_\mathrm{d} \right)^2+ \frac{\beta}{2} \left(- r \dot \phi \right)^2 + \sum_{i=1}^n \frac{\gamma_i}{2} {\ddot \theta}_i^2+\delta,
\end{equation}
for constants $a_\mathrm{e}$ and $b_\mathrm{e}$ and for fixed nonnegative constants $\upsilon_a$, $\upsilon_b$, $\alpha$, $\beta$, $\gamma_i$, $1 \le i \le n$, and $\delta$ so that the performance index is
\begin{equation} \label{eq_disk_J}
\begin{split}
J &\equiv p\left(a,\bx(a),b,\bx(b),\mu\right)+ \int_a^b L\left(t,\bx,\bu,\mu\right) \dt \\
&= \frac{\upsilon_a}{2}\left(a-a_\mathrm{e} \right)^2+\frac{\upsilon_b}{2}\left(b-b_\mathrm{e} \right)^2+\int_a^b \left[ \frac{\alpha}{2} \left(z_a - r \left(\phi-\phi_a \right) - z_\mathrm{d} \right)^2+ \frac{\beta}{2} \left(- r \dot \phi \right)^2  + \sum_{i=1}^n \frac{\gamma_i}{2} {\ddot \theta}_i^2 + \delta \right] \dt.
\end{split}
\end{equation}

\begin{figure}[h]
	\centering
	\includegraphics[width=0.5\linewidth]{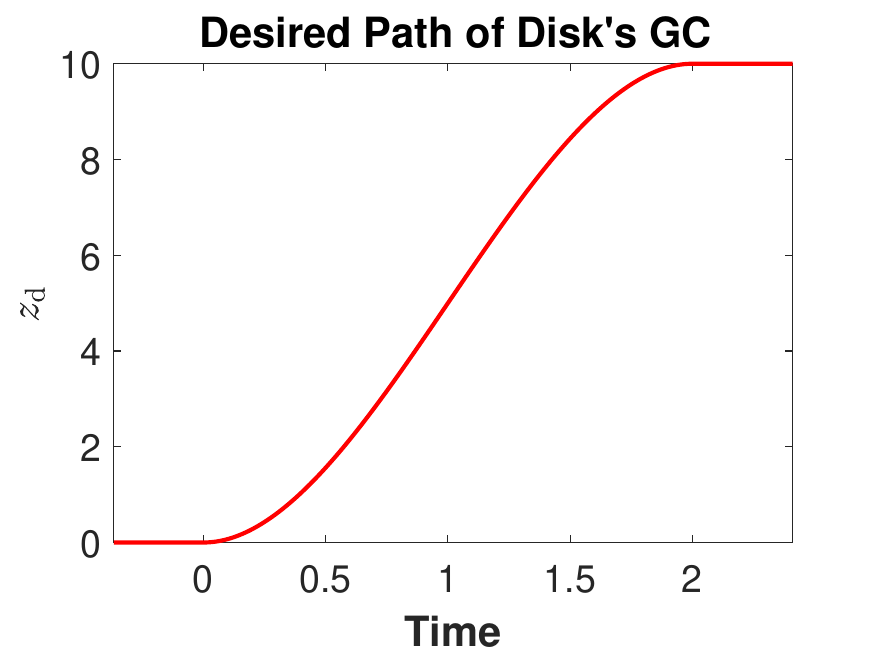}
	\caption{Plot of the desired path of the disk's GC. The disk's GC starts from rest at $z=0$ at time $t=0$, moves to the right for $0<t<2$, and stops at rest at $z=10$ at time $t=2$.}
	\label{fig:plot_rolling_disk_z_d}
\end{figure}

\paragraph{Discussion of different terms in the performance index.}
\revisionS{R1Q10}{The first summand $\frac{\upsilon_a}{2}\left(a-a_\mathrm{e} \right)^2$ in $p$ encourages the initial time $a$ to be near $a_\mathrm{e}$ if the initial time is free, while the second summand $\frac{\upsilon_b}{2}\left(b-b_\mathrm{e} \right)^2$ in $p$ encourages the final time $b$ to be near $b_\mathrm{e}$ if the final time is free. Note that these are soft constraints that encourage, but do not force, the values of $a$ and $b$ to be approximately $a_\mathrm{e}$ and $b_\mathrm{e}$, respectively. Similarly,} the first summand $\frac{\alpha}{2} \left(z_a - r (\phi-\phi_a) - z_\mathrm{d} \right)^2$ in $L$ encourages the disk's GC to track the desired spatial $\mathbf{e}_1$ path $z_\mathrm{d}$, the second summand $\frac{\beta}{2} \left(- r \dot \phi \right)^2 $ in $L$ encourages the disk's GC to track a minimum energy path, the next $n$ summands $\frac{\gamma_i}{2} {\ddot \theta}_i^2$, $1 \le i \le n$, in $L$ limit the magnitude of the acceleration of the $i^\mathrm{th}$ control mass parameterization, and the final summand $\delta$ in $L$ encourages a minimum time maneuver.
For example, the desired spatial $\mathbf{e}_1$ path might have the form
\begin{equation} \label{eq_z_d}
z_\mathrm{d}(t) \equiv \left[z_a w(t)+\tilde{z}_\mathrm{d}(t) \left(1-w(t)\right) \right] \left(1-y(t)\right)+z_b y(t),
\end{equation}
where
\begin{equation} \label{eq_sigmoid}
S(t) \equiv \frac{1}{2} \left[1+\tanh{\left( \frac{-t}{\epsilon}\right)} \right],
\end{equation}
\begin{equation} \label{eq_tran_sigmoidl}
w(t) \equiv S\left(t-a_\mathrm{e}\right),
\end{equation}
and
\begin{equation} \label{eq_tran_sigmoidr}
y(t) \equiv S\left(-t+b_\mathrm{e}\right).
\end{equation}
$z_\mathrm{d}$ \eqref{eq_z_d} holds steady at $z_a$ for $t < a_\mathrm{e}$, smoothly transitions between $z_a$ and $\tilde{z}_\mathrm{d}$ at $t = a_\mathrm{e}$, follows $\tilde{z}_\mathrm{d}$ for $a_\mathrm{e} < t < b_\mathrm{e}$, smoothly transitions between $\tilde{z}_\mathrm{d}$ and $z_b$ at $t = b_\mathrm{e}$, and holds steady at $z_b$ for $b_\mathrm{e} < t$.  $S$ \eqref{eq_sigmoid} is a time-reversed sigmoid function, i.e. a smooth approximation of the time-reversed unit step function; $\epsilon$ in \eqref{eq_sigmoid}  is a parameter such as $.01$ that determines how rapidly $S$ \eqref{eq_sigmoid} transitions from $1$ to $0$ at time 0. $w$ \eqref{eq_tran_sigmoidl} is the time-translation of $S$ \eqref{eq_sigmoid} to time $a_\mathrm{e}$ and $y$ \eqref{eq_tran_sigmoidr} is the time-translation of the time-reversal of $S$ \eqref{eq_sigmoid} to time $b_\mathrm{e}$. $w$ \eqref{eq_tran_sigmoidl} enables $z_\mathrm{d}$ \eqref{eq_z_d} to smoothly transition between $z_a$ and $\tilde{z}_\mathrm{d}$ at $t = a_\mathrm{e}$, while $y$ \eqref{eq_tran_sigmoidr} enables $z_\mathrm{d}$ \eqref{eq_z_d} to smoothly transition between $\tilde{z}_\mathrm{d}$ and $z_b$ at $t = b_\mathrm{e}$. $\tilde{z}_\mathrm{d}$, which appears in \eqref{eq_z_d}, might be the cubic polynomial
\begin{equation} \label{eq_tz_d}
\tilde{z}_\mathrm{d}(t) \equiv k_1 \left[-\frac{1}{3}t^3+\frac{1}{2}\left(a_\mathrm{e}+b_\mathrm{e} \right)t-a_\mathrm{e} b_\mathrm{e} t+k_2 \right]= k_1 \left[q(t)+k_2 \right],
\end{equation}
where 
\begin{equation} q(t) \equiv -\frac{1}{3}t^3+\frac{1}{2}\left(a_\mathrm{e}+b_\mathrm{e} \right)t-a_\mathrm{e} b_\mathrm{e} t, \quad k_1 \equiv \frac{z_b-z_a}{q(b_\mathrm{e})-q(a_\mathrm{e})}, \quad \mathrm{and} \quad k_2 \equiv \frac{z_b}{k_1}-q(b_\mathrm{e}).  \end{equation}
$\tilde{z}_\mathrm{d}$ has the special properties $\tilde{z}_\mathrm{d}(a_\mathrm{e})=z_a$, $\dot {\tilde{z}}_\mathrm{d}(a_\mathrm{e})=0$, $\tilde{z}_\mathrm{d}(b_\mathrm{e})=z_b$, and $\dot {\tilde{z}}_\mathrm{d}(b_\mathrm{e})=0$, so that the disk's GC is encouraged to start with zero velocity at $\mathbf{e}_1$-coordinate $z_a$ at $t=a_\mathrm{e}$ and to stop with zero velocity at $\mathbf{e}_1$-coordinate $z_b$ at $t=b_\mathrm{e}$. Figure~\ref{fig:plot_rolling_disk_z_d} illustrates \eqref{eq_z_d} using \eqref{eq_tz_d} with $a_\mathrm{e}=0$, $z_a=0$, $b_\mathrm{e}=2$, $z_b=10$, and $\epsilon=.01$.

\paragraph{Optimal control problem.}
The optimal control problem for the rolling disk is
\begin{equation}
\min_{a,b,\bu} J 
\mbox{\, s.t. \,}
\left\{
\begin{array}{ll}
\dot {\bx} = \mathbf{f}\left(t,\bx,\bu,\mu\right), \\
\bsigma\left(a,\bx(a),\mu\right) = \mathbf{0},\\
\bpsi\left(b,\bx(b),\mu\right) = \mathbf{0}.
\end{array}
\right.
\label{dyn_opt_problem_disk}
\end{equation}
Observe that the optimal control problem encapsulated by \eqref{dyn_opt_problem_disk} ignores path inequality constraints such as $\mathbf{D}\left(t,\bx,\bu,\mu\right) \le \mathbf{0}$, where $\mathbf{D}$ is an $r \times 1$ vector-valued function. Path inequality constraints can be incorporated in \eqref{dyn_opt_problem_disk} as soft constraints through penalty functions in the integrand cost function $L$ or the endpoint cost function $p$.

The indirect method \cite{hull2013optimal,RBOCnumerics} is applied now to \eqref{dyn_opt_problem_disk} to construct the endpoint function and regular Hamiltonian needed to formulate an ordinary differential equation two-point boundary value problem (ODE TPBVP), which renders the controlled equations of motion for the rolling disk.

\paragraph{Endpoint function, Pontryagin Hamiltonian, and regular Hamiltonian.}
The endpoint function is 
\begin{equation} \label{eq_disk_endpoint_fcn}
\begin{split}
G\left(a,\bx(a),\bxi,b,\bx(b),\bnu,\mu\right) &\equiv p\left(a,\bx(a),b,\bx(b),\mu\right)+\bxi^\mathsf{T} \bsigma\left(a,\bx(a),\mu\right)+\bnu^\mathsf{T} \bpsi\left(b,\bx(b),\mu\right) \\
&=\frac{\upsilon_a}{2}\left(a-a_\mathrm{e} \right)^2+\frac{\upsilon_b}{2}\left(b-b_\mathrm{e} \right)^2+ \bxi^\mathsf{T}  \begin{bmatrix} \btheta(a) - \btheta_a \\ \dot \btheta(a) - {\dot \btheta}_a \\ \phi(a)-\phi_a \\ -r \dot \phi(a) - {\dot z}_a  \end{bmatrix}\\
&\hphantom{=}+\bnu^\mathsf{T} \begin{bmatrix} \bPi \left(\tilde \Lambda\left(\phi(b)\right) \left[ \frac{1}{M} \sum_{i=0}^n m_i \bzeta_i\left(\theta_i(b)\right) \right] \right) - \bDelta_b \\ \dot \btheta(b) - {\dot \btheta}_b  \\ z_a - r \left(\phi(b)-\phi_a\right) - z_b \\ -r \dot \phi(b) - {\dot z}_b  \end{bmatrix},
\end{split}
\end{equation}
\revision{R1Q4}{where $\bxi \in \mathbb{R}^{2n+2}$ and $\bnu \in \mathbb{R}^{n+2+c}$, $c\in \left\{0,1,2\right\}$, are constant Lagrange multiplier vectors enforcing the initial and final conditions, \eqref{eq_disk_initial_conds} and \eqref{eq_disk_final_conds}, respectively.}  The control (Pontryagin) Hamiltonian is 
\begin{equation} \label{eq_disk_ham}
\begin{split}
H\left(t,\bx,\blam,\bu,\mu\right) &\equiv L\left(t,\bx,\bu,\mu\right) + \blam^\mathsf{T} \mathbf{f}\left(t,\bx,\bu,\mu\right) \\
&= \frac{\alpha}{2} \left(z_a - r \left(\phi-\phi_a \right) - z_\mathrm{d} \right)^2+ \frac{\beta}{2} \left(- r \dot \phi \right)^2  + \sum_{i=1}^n \frac{\gamma_i}{2} {\ddot \theta}_i^2+\delta + \blam^\mathsf{T} \begin{bmatrix} \dot \btheta \\ \bu  \\ \dot \phi \\ \kappa\left(t,\bx,\bu \right)  \end{bmatrix},
\end{split}
\end{equation}
\revision{R1Q4}{where $\blam \in \mathbb{R}^{2n+2}$ is a time-varying Lagrange multiplier vector enforcing the dynamics \eqref{eq_disk_dynamics}.} Recall  from \eqref{eqmo_chap_disk_4} that  the function describing the  uncontrolled equation of motion  for the rolling disk is 
\begin{equation} \label{eq_disk_kappa}
\kappa\left(t,\bx,\bu\right) \equiv \frac{ -r F_{\mathrm{e},1}+ \sum_{i=0}^n m_i K_i }{d_2+\sum_{i=0}^n m_i \left[\left( r \sin \phi + \zeta_{i,1} \right)^2+\left( r \cos \phi+ \zeta_{i,3} \right)^2 \right]},
\end{equation}
where  $K_i$ for $0 \le i \le n$  is defined  in \eqref{eq_K_i} as  
\begin{equation}
\begin{split}
K_i 
&\equiv \left(g+ r {\dot \phi}^2 \right) \left(\zeta_{i,3} \sin \phi - \zeta_{i,1} \cos \phi  \right)+
\left(r \cos \phi + \zeta_{i,3} \right) \left(- 2 \dot \phi {\dot \theta}_i \zeta_{i,3}^{\prime} + {\dot \theta}_i^2 \zeta_{i,1}^{\dprime} \right)\\
&\hphantom{\equiv}  - \left(r \sin \phi + \zeta_{i,1} \right) \left( 2 \dot \phi {\dot \theta}_i \zeta_{i,1}^{\prime}+ {\dot \theta}_i^2 \zeta_{i,3}^{\dprime} \right) 
+\left[ \left(r \cos \phi + \zeta_{i,3} \right)\zeta_{i,1}^{\prime}- \left(r \sin \phi + \zeta_{i,1} \right) \zeta_{i,3}^{\prime} \right] {\ddot \theta}_i \, . 
\end{split}
\end{equation}
Differentiating the Hamiltonian \eqref{eq_disk_ham} with respect to the components of the control $\bu$ gives
\begin{equation} \label{eq_H_u_i}
H_{u_i} = H_{\ddot{\theta}_i} = \gamma_i {\ddot \theta}_i + \lambda_{n+i} + \lambda_{2n+2} \frac{m_i \left[ \left(r \cos \phi + \zeta_{i,3} \right)\zeta_{i,1}^{\prime}- \left(r \sin \phi + \zeta_{i,1} \right) \zeta_{i,3}^{\prime} \right]}{d_2+\sum_{k=0}^n m_k \left[\left( r \sin \phi + \zeta_{k,1} \right)^2+\left( r \cos \phi+ \zeta_{k,3} \right)^2 \right]},
\end{equation} 
\begin{equation} 
H_{u_i u_j} = H_{\ddot{\theta}_i \ddot{\theta}_j} = \gamma_i \delta_{ij},
\end{equation}
and
\begin{equation} \label{eq_Huu_disk}
H_{\bu \bu} = \diag \begin{bmatrix} \gamma_1 & \gamma_2 & \dots & \gamma_n \end{bmatrix}.
\end{equation}
By \eqref{eq_Huu_disk},  $H_{\bu \bu}>0$ iff $\gamma_i >0$ for all $1 \le i \le n$. Consequently, the optimal control problem is regular iff $\gamma_i >0$ for all $1 \le i \le n$. \revision{R1Q5}{Let  us choose $\gamma_i >0$ for all $1 \le i \le n$, so that the optimal control problem is regular.} $H_{\bu}=0$ iff $H_{u_i}=0$ for all $1 \le i \le n$. From \eqref{eq_H_u_i},
\begin{equation}  \label{eq_ddtheta_exp}
\begin{split}
H_{u_i} &= 0 \iff \\  & {\ddot \theta}_i = \pi_i \left(t,\bx,\blam,\mu\right) \equiv -\gamma_i^{-1} \left\{ \lambda_{n+i} + \lambda_{2n+2} \frac{m_i \left[ \left(r \cos \phi + \zeta_{i,3} \right)\zeta_{i,1}^{\prime}- \left(r \sin \phi + \zeta_{i,1} \right) \zeta_{i,3}^{\prime} \right]}{d_2+\sum_{k=0}^n m_k \left[\left( r \sin \phi + \zeta_{k,1} \right)^2+\left( r \cos \phi+ \zeta_{k,3} \right)^2 \right]} \right\}.
\end{split}
\end{equation}
\eqref{eq_ddtheta_exp} shows that the control $\ddot \btheta$ may be expressed as a function $\bpi$ of $\bx$, $\blam$, and $\mu$; for generality, $\bpi$ will also depend on $t$ even though in this particular example it does not. The regular Hamiltonian is
\begin{equation} \label{eq_disk_regular_Hamiltonian}
\begin{split}
\hat H\left(t,\bx,\blam,\mu\right) &\equiv H\left(t,\bx,\blam,\bpi \left(t,\bx,\blam,\mu\right),\mu\right) \\
&= \frac{\alpha}{2} \left(z_a - r \left(\phi-\phi_a \right) - z_\mathrm{d} \right)^2+ \frac{\beta}{2} \left(- r \dot \phi \right)^2  + \sum_{i=1}^n \frac{\gamma_i}{2} {\pi}_i^2 \left(t,\bx,\blam,\mu\right)+\delta \\
&\hphantom{=} + \blam^\mathsf{T} \begin{bmatrix} \dot \btheta \\ \bpi \left(t,\bx,\blam,\mu\right)  \\ \dot \phi \\ \kappa\left(t,\bx,\bpi \left(t,\bx,\blam,\mu\right) \right)  \end{bmatrix}.
\end{split}
\end{equation}

\paragraph{Controlled equations of motion.}
As explained in \cite{caillau2012differential,RBOCnumerics}, one way to solve the optimal control problem \eqref{dyn_opt_problem_disk} for the rolling disk is to solve the ODE TPBVP:
\begin{equation} \label{eq_pmp_bvp_disk}
\begin{split}
\dot {\bx} &= \hat{H}_{\blam}^\mathsf{T} \left(t,\bx,\blam,\mu\right) = \hat{\mathbf{f}}\left(t,\bx,\blam,\mu\right)\equiv \mathbf{f}\left(t,\bx,\blam,\bpi\left(t,\bx,\blam,\mu\right),\mu\right), \\
\dot {\blam} &= - \hat{H}_{\bx}^\mathsf{T} \left(t,\bx,\blam,\mu\right)=-H_{\bx}^\mathsf{T}\left(t,\bx,\blam,\bpi\left(t,\bx,\blam,\mu\right),\mu\right), \\
\left. \hat{H} \right|_{t=a} = G_{a}, \quad \left. \blam \right|_{t=a} &= -G_{\bx(a)}^\mathsf{T}, \quad G_{\bxi}^\mathsf{T} = \bsigma\left(a,\bx(a),\mu\right) = \mathbf{0}, \\
\left. \hat{H} \right|_{t=b} = -G_{b}, \quad \left. \blam \right|_{t=b} &= G_{\bx(b)}^\mathsf{T}, \quad G_{\bnu}^\mathsf{T} = \bpsi\left(b,\bx(b),\mu\right) = \mathbf{0}.
\end{split}
\end{equation}
Subappendix~\ref{app_ccem_rdisk} derives the formulas for constructing $H_{\bx}^\mathsf{T}$.

\rem{One way to solve the optimal control problem \eqref{dyn_opt_problem_disk} for the rolling disk is to solve the ODE TPBVP given by \eqref{eq_pmp_bvp}, \eqref{eq_pmp_lbc2}, and \eqref{eq_pmp_rbc2} using the endpoint function \eqref{eq_disk_endpoint_fcn} and the regular Hamiltonian \eqref{eq_disk_regular_Hamiltonian}. In order to realize the ODE velocity function in \eqref{eq_pmp_bvp}, $\hat{\mathbf{f}}$ and $\hat{H}_{\bx}$ must be constructed. From its definition \eqref{eq_fhat_def}, $\hat{\mathbf{f}}$ is readily constructed from $\mathbf{f}$ and $\bpi$. From \eqref{eq_Hhat_x}, $\hat{H}_{\bx}$ may be obtained from $H_{\bx}$ and $\bpi$; Subappendix~\ref{app_ccem_rdisk} gives the formulas for constructing $H_{\bx}$, so that $\hat{H}_{\bx}$ may be realized without resorting to computational tools like symbolic or automatic differentiation.}

\subsection{Controlled Equations of Motion for the Rolling Ball} \label{ssec_ball_controlled}
\rem{Before proceeding, some useful terminology is defined or recalled. Given a  vector 
	\begin{equation}
	\bv = \begin{bmatrix} v_1 \\ v_2 \\ v_3 \end{bmatrix} = \begin{bmatrix} v_1 & v_2 & v_3 \end{bmatrix}^\mathsf{T} \in \mathbb{R}^3,
	\end{equation}
	the projected vector consisting of the first two components of $\bv$ is
	\begin{equation}
	\bv_{12} = \begin{bmatrix}v_1 \\ v_2 \end{bmatrix} = \begin{bmatrix}v_1 & v_2 \end{bmatrix}^\mathsf{T} \in \mathbb{R}^2.
	\end{equation}
	Since a versor is used to parameterize the rolling ball's orientation matrix, quaternions and versors are briefly reviewed here; see Appendix~\ref{app_quaternions} and the references mentioned there for a more complete review. 
	$\mathbb{H}$ denotes the set of quaternions, which is isomorphic to $\mathbb{R}^4$. A quaternion $\mathfrak{p} \in \mathbb{H}$ can be expressed as the column vector 
	\begin{equation}
	\mathfrak{p} =
	\begin{bmatrix} p_0 \\ p_1 \\ p_2 \\ p_3 \end{bmatrix}=
	\begin{bmatrix} p_0 & p_1 & p_2 & p_3 \end{bmatrix}^\mathsf{T}=\begin{bmatrix} p_0\,;\,p_1\,;\,p_2\,;\,p_3 \end{bmatrix}.
	\end{equation}
	Given a column vector $\bv \in \mathbb{R}^3$, $\bv^\sharp$ is the quaternion $\begin{bmatrix}0\,;\,\bv\end{bmatrix} \in \mathbb{H}$; that is, 
	\begin{equation}
	\bv^\sharp = \begin{bmatrix}0 \\ \bv \end{bmatrix} = \begin{bmatrix}0\,;\,\bv\end{bmatrix}.
	\end{equation}
	Given a quaternion $\mathfrak{p} \in \mathbb{H}$, $\mathfrak{p}^\flat \in \mathbb{R}^3$ is the column vector such that 
	\begin{equation}
	\mathfrak{p} = \begin{bmatrix} p_0 \\ \mathfrak{p}^\flat \end{bmatrix}= \begin{bmatrix} p_0\,;\,\mathfrak{p}^\flat \end{bmatrix}.
	\end{equation}
	Given a column vector $\bv \in \mathbb{R}^3$, note that 
	\begin{equation}
	\left(\bv^\sharp\right)^\flat=\bv.
	\end{equation}
	However, given a quaternion $\mathfrak{p} \in \mathbb{H}$, 
	\begin{equation}
	\left( \mathfrak{p}^\flat \right)^\sharp = \mathfrak{p} \quad \mathrm{iff} \quad \mathfrak{p} = \begin{bmatrix} 0 \\ \mathfrak{p}^\flat \end{bmatrix}= \begin{bmatrix}0\,;\,\mathfrak{p}^\flat \end{bmatrix}. 
	\end{equation}
	Given quaternions $\mathfrak{p},\mathfrak{q} \in \mathbb{H}$, their product is
	\begin{equation}
	\mathfrak{p} \mathfrak{q}=\begin{bmatrix} p_0\,;\,\mathfrak{p}^\flat \end{bmatrix} \begin{bmatrix}q_0\,;\,\mathfrak{q}^\flat \end{bmatrix}=\begin{bmatrix}p_0q_0-\mathfrak{p}^\flat \cdot \mathfrak{q}^\flat\,;\,p_0\mathfrak{q}^\flat + q_0\mathfrak{p}^\flat+\mathfrak{p}^\flat \times \mathfrak{q}^\flat \end{bmatrix}
	\end{equation}
	and their dot product is 
	\begin{equation}
	\mathfrak{p} \cdot \mathfrak{q}=\begin{bmatrix} p_0\,;\,\mathfrak{p}^\flat \end{bmatrix} \cdot \begin{bmatrix}q_0\,;\,\mathfrak{q}^\flat \end{bmatrix}=\begin{bmatrix}p_0\,;\,p_1\,;\,p_2\,;\,p_3\end{bmatrix} \cdot \begin{bmatrix}q_0\,;\,q_1\,;\,q_2\,;\,q_3\end{bmatrix}= p_0 q_0 + \mathfrak{p}^\flat \cdot \mathfrak{q}^\flat=p_0 q_0 +p_1 q_1 + p_2 q_2 + p_3 q_3.
	\end{equation}
	$\mathscr{S} \subset \mathbb{H}$ denotes the set of unit quaternions, also called versors, which is isomorphic to $\mathbb{S}^3 \subset \mathbb{R}^4$. A versor $\mathfrak{q} \in \mathscr{S}$ can be expressed as the column vector 
	\begin{equation}
	\mathfrak{q} = \begin{bmatrix} q_0 \\ q_1 \\ q_2 \\ q_3 \end{bmatrix}=
	\begin{bmatrix} q_0 & q_1 & q_2 & q_3 \end{bmatrix}^\mathsf{T}=\begin{bmatrix} q_0\,;\,q_1\,;\,q_2\,;\,q_3 \end{bmatrix} \quad \mathrm{such \; that} \quad \mathfrak{q} \cdot \mathfrak{q} = q_0^2+q_1^2+q_2^2+q_3^2=1.
	\end{equation}
	The rolling ball's orientation matrix $\Lambda \in SO(3)$ is parameterized by the versor $\mathfrak{q} \in \mathscr{S}$. If $\bOm \in \mathbb{R}^3$ is the rolling ball's body angular velocity, then the time derivative of $\mathfrak{q}$ is 
	\begin{equation}
	\dot {\mathfrak{q}} = \frac{1}{2} \mathfrak{q} \bOm^\sharp.
	\end{equation}
	If $\bY \in \mathbb{R}^3$ is a body frame vector, then the rotation of $\bY$ by $\Lambda$ is 
	\begin{equation}
	\Lambda \bY = \left[\mathfrak{q} \bY^\sharp \mathfrak{q}^{-1} \right]^\flat. 
	\end{equation}
	If $\mathbf{y} \in \mathbb{R}^3$ is a spatial frame vector, then the rotation of $\mathbf{y}$ by $\Lambda^{-1}$ is 
	\begin{equation}
	\Lambda^{-1} \mathbf{y} = \left[\mathfrak{q}^{-1} \mathbf{y}^\sharp \mathfrak{q} \right]^\flat.
	\end{equation}
}

\paragraph{Notation, versors, assumptions, and the motion of the center of mass.}
Having derived the controlled equations of motion for the rolling disk, the controlled equations of motion are now developed for the rolling ball. The rolling ball's orientation matrix $\Lambda \in SO(3)$ is parameterized by $\mathfrak{q} \in \mathscr{S}$, where $\mathscr{S}$ denotes the set of versors (i.e. unit quaternions) \cite{Ho2011_pII,graf2008quaternions,stevens2015aircraft,baraff2001physically}. The properties of versors and the notation used to manipulate versors are explained in Appendix D of \cite{Putkaradze2018dynamicsP}. \revision{R1Q14}{Recall from \cite{Putkaradze2018dynamicsP} that given a column vector $\bv \in \mathbb{R}^3$, $\bv^\sharp$ is the quaternion 
	\begin{equation}
	\bv^\sharp = \begin{bmatrix}0 \\ \bv \end{bmatrix},
	\end{equation}
	and given a quaternion $\mathfrak{p} \in \mathbb{H}$, $\mathfrak{p}^\flat \in \mathbb{R}^3$ is the column vector such that 
	\begin{equation}
	\mathfrak{p} = \begin{bmatrix} p_0 \\ \mathfrak{p}^\flat \end{bmatrix}.
	\end{equation} }

\revision{R1Q2}{It is desired to roll the ball from some initial configuration at a prescribed or free initial time $a$ to some final configuration at a prescribed or free final time $b$, without moving the control masses too rapidly along their control rails, while the ball's GC tracks a prescribed or minimum energy path, while the ball avoids some obstacles, or while minimizing the maneuver's duration $b-a$.} Concretely, at the prescribed or free initial time $a$, the positions of the control mass parameterizations are prescribed to be $\btheta(a)=\btheta_a$, the velocities of the control mass parameterizations are prescribed to be $\dot \btheta(a) = {\dot \btheta}_a$, the orientation of the ball is prescribed to be $\mathfrak{q}(a)=\mathfrak{q}_a$, the body angular velocity of the ball is prescribed to be $\bOm(a) = \bOm_a$, and the spatial $\mathbf{e}_1$-$\mathbf{e}_2$ position of the ball's GC is prescribed to be $\bz(a)=\bz_a$. 

\revisionS{R1Q8}{ In the body frame, the position of the center of mass relative to the geometric center at the final time is given by 
$\frac{1}{M} \sum_{i=0}^n m_i \bzeta_i\left(\theta_i(b)\right)$. In the spatial frame, the same position is given by $\Lambda(b) \left[\frac{1}{M} \sum_{i=0}^n m_i \bzeta_i\left(\theta_i(b)\right)\right]$.   Then,  using the versor parameterization of $\Lambda$, at the prescribed or free final time $b$,  the components determined by the projection operator $\bPi$ of the position of the ball's center of mass expressed in the spatial frame translated to the GC are prescribed to be 
}
\begin{equation}
\label{eqCM_spatial} 
\bPi \left(\left[\mathfrak{q}(b) \left[ \frac{1}{M} \sum_{i=0}^n m_i \bzeta_i\left(\theta_i(b)\right) \right]^\sharp \mathfrak{q}(b)^{-1} \right]^\flat \right) = \bDelta_b.
\end{equation} 
Moreover, at the final time, the velocities of the control mass parameterizations are prescribed to be $\dot \btheta(b) = {\dot \btheta}_b$, the body angular velocity of the ball is prescribed to be $\bOm(b) = \bOm_b$, and the spatial $\mathbf{e}_1$-$\mathbf{e}_2$ position of the ball's GC is prescribed to be $\bz(b)=\bz_b$. 

For example, if it is desired to start and stop the ball at rest, then $\bPi$ is the projection onto the first two components, $\bDelta_b=\begin{bmatrix} 0 \\ 0 \end{bmatrix}$, $\btheta_a$ and $\mathfrak{q}_a$ are such that 
\begin{equation} \label{eq_ball_stab_init}
\bPi \left(\left[\mathfrak{q}_a \left[ \frac{1}{M} \sum_{i=0}^n m_i \bzeta_i\left(\theta_{a,i}\right) \right]^\sharp \mathfrak{q}_a^{-1} \right]^\flat \right) = \begin{bmatrix} 0 \\ 0 \end{bmatrix},
\end{equation} 
${\dot \btheta}_a = \mathbf{0}$, $\bOm_a = \mathbf{0}$, $\btheta(b)$ and $\mathfrak{q}(b)$ are such that  
\begin{equation} \label{eq_ball_stab_fin}
\bPi \left(\left[\mathfrak{q}(b) \left[ \frac{1}{M} \sum_{i=0}^n m_i \bzeta_i\left(\theta_i(b)\right) \right]^\sharp \mathfrak{q}(b)^{-1} \right]^\flat \right) = \begin{bmatrix} 0 \\ 0 \end{bmatrix},
\end{equation} 
${\dot \btheta}_b = \mathbf{0}$, and ${\bOm}_b = \mathbf{0}$. With this choice of $\bPi$, \eqref{eq_ball_stab_init} and \eqref{eq_ball_stab_fin} mean that the CM in the spatial frame translated to the GC is above or below the GC at the initial and final times.

\paragraph{State variables, control variables, and dynamics.}
The system state $\bx$ and control $\bu$ are
\begin{equation}
\bx \equiv \begin{bmatrix} \btheta \\ \dot \btheta \\ \mathfrak{q} \\ \bOm \\ \bz  \end{bmatrix} \quad \mathrm{and} \quad \bu \equiv \ddot \btheta,
\end{equation}
where $\btheta, \, \dot \btheta, \, \ddot \btheta \in \mathbb{R}^n$, $\mathfrak{q} \in \mathscr{S} \cong \mathbb{S}^3 \subset \mathbb{R}^4$, $\bOm \in \mathbb{R}^3$, and $\bz \in \mathbb{R}^2$. 
The system dynamics defined for $a \le t \le b$ are
\begin{equation} \label{rolling_ball_opt_con_f}
\dot {\bx} = \begin{bmatrix} \dot \btheta \\ \ddot \btheta \\ \dot {\mathfrak{q}} \\ \dot \bOm \\ \dot \bz  \end{bmatrix}  = \mathbf{f}\left(t,\bx,\bu,\mu\right) \equiv \begin{bmatrix} \dot \btheta \\ \bu  \\ \frac{1}{2} \mathfrak{q} \bOm^\sharp \\ \bkappa\left(t,\bx,\bu \right) \\ \left( \left[\mathfrak{q} \bOm^\sharp \mathfrak{q}^{-1} \right]^\flat \times r \mathbf{e}_3  \right)_{12}  \end{bmatrix},
\end{equation}
\revision{R1Q4}{where $\mu$ is a scalar continuation parameter and} $\bkappa\left(t,\bx,\bu \right)$ is given by the right-hand side of the formula for $\dot \bOm$ in \eqref{uncon_ball_eqns_explicit_1d}:
\begin{equation} \label{eq_ball_kappa}
\begin{split}
\bkappa\left(t,\bx,\bu \right) \equiv
\left[\sum_{i=0}^n m_i \widehat{\mathbf{s}_i}^2  -\inertia \right]^{-1}  \Bigg[&\bOm \times \inertia \bOm+r \tilde \bGamma \times \bGamma\\
&+ \sum_{i=0}^n m_i \mathbf{s}_i \times  \left\{ g \bGamma+ \bOm \times \left(\bOm \times \bzeta_i +2 \dot \theta_i \bzeta_i^{\prime} \right) + \dot \theta_i^2 \bzeta_i^{\dprime} + \ddot \theta_i  \bzeta_i^{\prime} \right\}  \Bigg].
\end{split}
\end{equation}
Note that in order to construct $\bkappa\left(t,\bx,\bu \right)$,  $\bGamma \equiv \Lambda^{-1} \mathbf{e}_3$ and $\tilde \bGamma \equiv \Lambda^{-1} \mathbf{F}_\mathrm{e}$ must be constructed. Given $\mathfrak{q}$, this can be accomplished by first constructing $\Lambda$ from $\mathfrak{q}$ or directly from $\mathfrak{q}$ by using the formulas $\bGamma \equiv \Lambda^{-1} \mathbf{e}_3=\left[\mathfrak{q}^{-1} \mathbf{e}_3^\sharp \mathfrak{q} \right]^\flat$ and $\tilde \bGamma \equiv \Lambda^{-1} \mathbf{F}_\mathrm{e} = \left[\mathfrak{q}^{-1} \mathbf{F}_\mathrm{e}^\sharp \mathfrak{q} \right]^\flat$. Likewise, the final formula in \eqref{rolling_ball_opt_con_f} is $\dot \bz=\left( \bom \times r \mathbf{e}_3  \right)_{12}$, where $\bom \equiv \left[ \dot{\Lambda} \Lambda^{-1} \right]^\vee =  \Lambda \bOm=\left[\mathfrak{q} \bOm^\sharp \mathfrak{q}^{-1} \right]^\flat$. Thus, given $\mathfrak{q}$, $\bom$ can be constructed by first constructing $\Lambda$ from $\mathfrak{q}$ or directly from $\mathfrak{q}$ via $\bom=\left[\mathfrak{q} \bOm^\sharp \mathfrak{q}^{-1} \right]^\flat$. The most computationally efficient method to construct $\bGamma \equiv \Lambda^{-1} \mathbf{e}_3=\left[\mathfrak{q}^{-1} \mathbf{e}_3^\sharp \mathfrak{q} \right]^\flat$, $\tilde \bGamma \equiv \Lambda^{-1} \mathbf{F}_\mathrm{e} = \left[\mathfrak{q}^{-1} \mathbf{F}_\mathrm{e}^\sharp \mathfrak{q} \right]^\flat$, and $\bom = \Lambda \bOm=\left[\mathfrak{q} \bOm^\sharp \mathfrak{q}^{-1} \right]^\flat$ is to construct $\Lambda$ from $\mathfrak{q}$ and to then multiply $\Lambda^{-1} = \Lambda^\mathsf{T}$ against $\mathbf{e}_3$ and $\mathbf{F}_\mathrm{e}$ and to multiply $\Lambda$ against $\bOm$.

\paragraph{Initial and final conditions.}
The prescribed initial conditions at time $t=a$ are
\begin{equation} \label{eq_ball_initial_conds}
\bsigma\left(a,\bx(a),\mu\right) \equiv \begin{bmatrix} \btheta(a) - \btheta_a \\ \dot \btheta(a) - {\dot \btheta}_a \\ \mathfrak{q}(a)-\mathfrak{q}_a \\ \bOm(a) - \bOm_a \\ \bz(a)- \bz_a \end{bmatrix} = \mathbf{0},
\end{equation}
and the prescribed final conditions at time $t=b$ are
\begin{equation} \label{eq_ball_final_conds}
\bpsi\left(b,\bx(b),\mu\right) \equiv \begin{bmatrix} \bPi \left(\left[\mathfrak{q}(b) \left[ \frac{1}{M} \sum_{i=0}^n m_i \bzeta_i\left(\theta_i(b)\right) \right]^\sharp \mathfrak{q}(b)^{-1} \right]^\flat \right) - \bDelta_b \\ \dot \btheta(b) - {\dot \btheta}_b \\ \bOm(b) - \bOm_b \\ \bz(b)- \bz_b  \end{bmatrix} = \mathbf{0}.
\end{equation}

\paragraph{Cost functions and performance index.}
Consider the endpoint and integrand cost functions
\begin{equation} \label{eq_ball_endpoint_cost}
p\left(a,\bx(a),b,\bx(b),\mu\right) \equiv \frac{\upsilon_a}{2}\left(a-a_\mathrm{e} \right)^2+\frac{\upsilon_b}{2}\left(b-b_\mathrm{e} \right)^2 \\
\end{equation}
and
\begin{equation} \label{eq_ball_integrand_cost}
L\left(t,\bx,\bu,\mu\right) \equiv \frac{\alpha}{2} \left| \bz - \bz_\mathrm{d} \right|^2+ \frac{\beta}{2} \left| \left( \left[\mathfrak{q} \bOm^\sharp \mathfrak{q}^{-1} \right]^\flat \times r \mathbf{e}_3  \right)_{12}\right|^2 + \sum_{i=1}^n \frac{\gamma_i}{2} {\ddot \theta}_i^2+\sum_{j=1}^K V_j \left(\bz,\mu \right)+\delta,
\end{equation}
for constants $a_\mathrm{e}$ and $b_\mathrm{e}$ and for fixed nonnegative constants $\upsilon_a$, $\upsilon_b$, $\alpha$, $\beta$, $\gamma_i$, $1 \le i \le n$, and $\delta$ so that the performance index is
\begin{equation} \label{eq_ball_J}
\begin{split}
J &\equiv p\left(a,\bx(a),b,\bx(b),\mu\right)+ \int_a^b L\left(t,\bx,\bu,\mu\right) \dt \\
&= \frac{\upsilon_a}{2}\left(a-a_\mathrm{e} \right)^2+\frac{\upsilon_b}{2}\left(b-b_\mathrm{e} \right)^2\\
&\hphantom{=}+\int_a^b \left[ \frac{\alpha}{2} \left| \bz - \bz_\mathrm{d} \right|^2+ \frac{\beta}{2} \left| \left( \left[\mathfrak{q} \bOm^\sharp \mathfrak{q}^{-1} \right]^\flat \times r \mathbf{e}_3  \right)_{12}\right|^2 + \sum_{i=1}^n \frac{\gamma_i}{2} {\ddot \theta}_i^2+\sum_{j=1}^K V_j \left(\bz,\mu \right)+\delta \right] \dt.
\end{split}
\end{equation}

\paragraph{Discussion of different terms in the performance index.}
The first summand $\frac{\upsilon_a}{2}\left(a-a_\mathrm{e} \right)^2$ in $p$ encourages the initial time $a$ to be near $a_\mathrm{e}$ if the initial time is free, while the second summand $\frac{\upsilon_b}{2}\left(b-b_\mathrm{e} \right)^2$ in $p$ encourages the final time $b$ to be near $b_\mathrm{e}$ if the final time is free. The first summand $\frac{\alpha}{2} \left| \bz - \bz_\mathrm{d} \right|^2$ in $L$ encourages the ball's GC to track the desired spatial $\mathbf{e}_1$-$\mathbf{e}_2$ path $\bz_\mathrm{d}$, the second summand $\frac{\beta}{2} \left| \left( \left[\mathfrak{q} \bOm^\sharp \mathfrak{q}^{-1} \right]^\flat \times r \mathbf{e}_3  \right)_{12}\right|^2$ in $L$ encourages the ball's GC to track a minimum energy path, the next $n$ summands $\frac{\gamma_i}{2} {\ddot \theta}_i^2$, $1 \le i \le n$, in $L$ limit the magnitude of the acceleration of the $i^\mathrm{th}$ control mass parameterization,  the next $K$ summands $V_j\left(\bz,\mu\right)$, $1 \le j \le K$, in $L$ represent obstacles to be avoided, and the final summand $\delta$ in $L$ encourages a minimum time maneuver. It is useful to note that one can also consider more general terms in $L$ to minimize the magnitudes of the control mass accelerations. For example, to minimize the magnitudes of the control mass accelerations in the body frame, $\sum_{i=1}^n \frac{\gamma_i}{2} \left| \dot \theta_i^2 \bzeta_i^{\dprime} + \ddot \theta_i  \bzeta_i^{\prime} \right|^2$ should appear in $L$. Such expressions lead to quite cumbersome expressions in the controlled equations of motion and strongly depend on the shapes of the 1-d rails. We thus use the expression $\sum_{i=1}^n \frac{\gamma_i}{2} \ddot \theta_i^2$ in $L$ as an acceptable compromise.

Note that a solution  obtained by the optimal control procedure, that minimizes \eqref{eq_disk_J} for the rolling disk or \eqref{eq_ball_J} for the rolling ball, is a ``compromise" between several, often conflicting, components, where some components of the performance index can be made more prominent by making their coefficients appropriately larger. The minimization of the  performance index  does not guarantee the minimization of each  component individually.

\begin{figure}[h]
	\centering
	\includegraphics[width=0.5\linewidth]{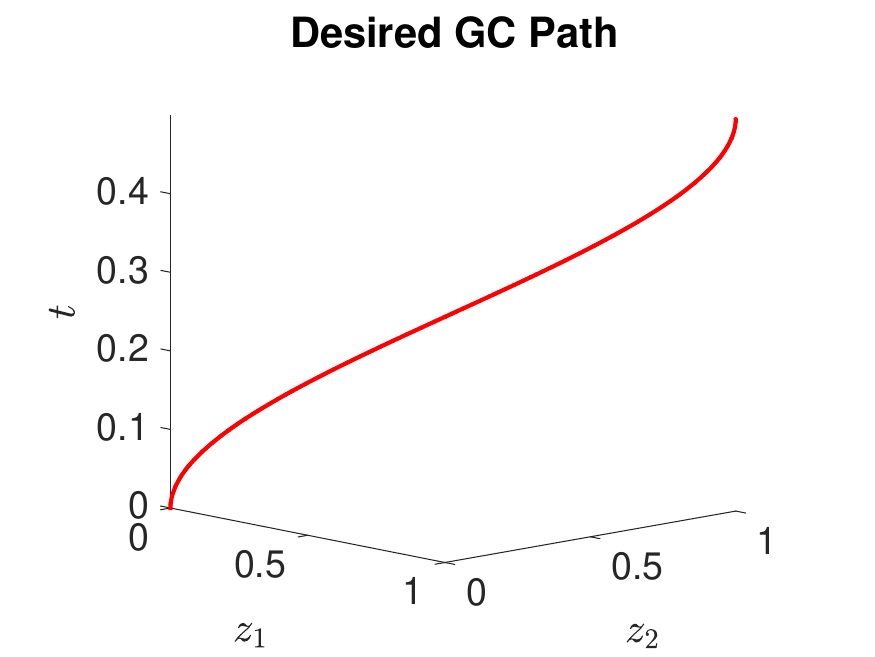}
	\caption{Plot of the desired path of the ball's GC. The ball's GC starts from rest at $\bz=\protect\begin{bmatrix}  0 & 0 \protect\end{bmatrix}^\mathsf{T}$ at time $t=0$ and stops at rest at $\bz=\protect\begin{bmatrix}  1 & 1 \protect\end{bmatrix}^\mathsf{T}$ at time $t=.5$.}
	\label{fig_rolling_ball_z_d}
\end{figure}

For example, the desired spatial $\mathbf{e}_1$-$\mathbf{e}_2$ path might have the form
\begin{equation} \label{eq_ball_z_d}
\bz_\mathrm{d}(t) \equiv \left[\bz_a w(t)+\tilde{\bz}_\mathrm{d}(t) \left(1-w(t)\right) \right] \left(1-y(t)\right)+\bz_b y(t),
\end{equation}
where $w$, $y$, and $S$ are given by \eqref{eq_tran_sigmoidl}, \eqref{eq_tran_sigmoidr}, and \eqref{eq_sigmoid}, respectively.
$\bz_\mathrm{d}$ \eqref{eq_ball_z_d} holds steady at $\bz_a$ for $t < a_\mathrm{e}$, smoothly transitions between $\bz_a$ and $\tilde{\bz}_\mathrm{d}$ at $t = a_\mathrm{e}$, follows $\tilde{\bz}_\mathrm{d}$ for $a_\mathrm{e} < t < b_\mathrm{e}$, smoothly transitions between $\tilde{\bz}_\mathrm{d}$ and $\bz_b$ at $t = b_\mathrm{e}$, and holds steady at $\bz_b$ for $b_\mathrm{e} < t$.  $\tilde{\bz}_\mathrm{d}$, which appears in \eqref{eq_ball_z_d}, might be 
\begin{equation} \label{eq_ball_tz_d}
\tilde{\bz}_\mathrm{d}(t) \equiv \bk_1 \left[ \left\{ -\frac{1}{3}t^3+\frac{1}{2}\left(a_\mathrm{e}+b_\mathrm{e} \right)t-a_\mathrm{e} b_\mathrm{e} t \right\}\begin{bmatrix}  1 \\ 1 \end{bmatrix}+\bk_2 \right]= \bk_1 \left[q(t) \begin{bmatrix}  1 \\ 1 \end{bmatrix}+\bk_2 \right],
\end{equation}
where 
\begin{equation} \label{eq_ball_tz_details}
q(t) \equiv -\frac{1}{3}t^3+\frac{1}{2}\left(a_\mathrm{e}+b_\mathrm{e} \right)t-a_\mathrm{e} b_\mathrm{e} t, \quad \bk_1 \equiv \frac{\bz_b-\bz_a}{q(b_\mathrm{e})-q(a_\mathrm{e})}, \quad \mathrm{and} \quad \bk_2 \equiv \frac{\bz_b}{\bk_1}-q(b_\mathrm{e}) \begin{bmatrix}  1 \\ 1 \end{bmatrix}.  \end{equation}
The multiplication between the vectors $\bk_1$ and $q(t) \begin{bmatrix}  1 \\ 1 \end{bmatrix}+\bk_2$ in \eqref{eq_ball_tz_d} is meant to be performed componentwise.
The division $\frac{\bz_b}{\bk_1}$ used in the construction of $\bk_2$ in \eqref{eq_ball_tz_details} is meant to be performed componentwise; to avoid division by zero, if a component of $\bk_1$ is zero, then the corresponding component of $\bk_2$ is set to zero. $\tilde{\bz}_\mathrm{d}$ has the special properties $\tilde{\bz}_\mathrm{d}(a_\mathrm{e})=\bz_a$, $\dot {\tilde{\bz}}_\mathrm{d}(a_\mathrm{e})=\begin{bmatrix}  0 & 0 \end{bmatrix}^\mathsf{T}$, $\tilde{\bz}_\mathrm{d}(b_\mathrm{e})=\bz_b$, and $\dot {\tilde{\bz}}_\mathrm{d}(b_\mathrm{e})=\begin{bmatrix}  0 & 0 \end{bmatrix}^\mathsf{T}$, so that the ball's GC is encouraged to start with zero velocity at $\bz_a$ at $t=a_\mathrm{e}$ and to stop with zero velocity at $\bz_b$ at $t=b_\mathrm{e}$. Figure~\ref{fig_rolling_ball_z_d} illustrates \eqref{eq_ball_z_d} using \eqref{eq_ball_tz_d} with $a_\mathrm{e}=0$, $\bz_a=\begin{bmatrix}  0 & 0 \end{bmatrix}^\mathsf{T}$, $b_\mathrm{e}=.5$, $\bz_b=\begin{bmatrix}  1 & 1 \end{bmatrix}^\mathsf{T}$, and $\epsilon=.01$.

Illustratively, for $1 \le j \le K$, the $j^\mathrm{th}$ obstacle of height $h_j$ and radius $\rho_j$ with center at spatial $\mathbf{e}_1$-$\mathbf{e}_2$  position $\boldsymbol{v}_j = \begin{bmatrix}  v_{j,1} & v_{j,2} \end{bmatrix}^\mathsf{T}$ might be modeled via the function
\begin{equation} \label{eq_obstacle_pot_func}
V_j\left(\bz,\mu\right) \equiv h_j \, S \left(\left| \bz - \bv_j \right| - \rho_j \right),
\end{equation}
where
\begin{equation} \label{eq_obst_sigmoid}
S(y)\equiv \frac{1}{2}+\frac{1}{2} \tanh \left(\frac{-y}{\epsilon} \right)
\end{equation} 
is a time-reversed sigmoid function or
\begin{equation} \label{eq_obst_cutoff}
S(y)\equiv \left[\max \left\{0,-y \right\} \right]^4  =\ReLU^4\left(-y\right) 
\end{equation}
is a $C^2$ cutoff function.  In \eqref{eq_obst_cutoff}, $\ReLU$ is the rectified linear unit  function frequently used in the machine learning literature \cite{lecun2015deep}. 
\rem{where $S$ is given by \eqref{eq_sigmoid}.}
As indicated by \eqref{eq_obstacle_pot_func}, the radial distance from the ball's GC to the obstacle center should exceed the obstacle radius $\rho_j$ for successful obstacle avoidance. In order to encourage the entire ball to avoid the obstacle, the obstacle radius $\rho_j$ must include the ball's radius $r$. That is, if the physical radius of the obstacle is $\epsilon_j$, then set $\rho_j=r+\epsilon_j$ to encourage the entire ball to stay away from the obstacle; if $\rho_j=\epsilon_j$, then only the ball's GC is encouraged to stay away from the obstacle. \revisionS{R1Q11}{Note that the obstacle avoidance function is independent of the control variables $\bu=\ddot \btheta$, which is important for the regularity of the control (Pontryagin) Hamiltonian. } 

\paragraph{Optimal control problem.} 
\revisionS{R2Q2}{ The ODE formulation of the optimal control problem for the rolling ball  is the following constrained minimization of the performance index:  }
\begin{equation}
\min_{a,b,\bu} J 
\mbox{\, s.t. \,}
\left\{
\begin{array}{ll}
\dot {\bx} = \mathbf{f}\left(t,\bx,\bu,\mu\right), \\
\bsigma\left(a,\bx(a),\mu\right) = \mathbf{0},\\
\bpsi\left(b,\bx(b),\mu\right) = \mathbf{0}.
\end{array}
\right.
\label{dyn_opt_problem_ball}
\end{equation}
\revisionS{R2Q2}{ One method for deriving the controlled equations of motion is to compute the Euler-Lagrange equations corresponding to an extremum of \eqref{dyn_opt_problem_ball}, provided the appropriate differentiability conditions on $J$ are satisfied, leading to a differential-algebraic equation two-point boundary value problem (DAE TPBVP).  This is known as the indirect method. The direct method, in contrast, computes a minimum of \eqref{dyn_opt_problem_ball} directly  by a finite-dimensional parameterization of $\bu$ and by discretizing the uncontrolled equations of motion in the  set of constraints. This paper focuses on  the indirect method, which is preferable for a theoretical treatment. }
 
There are also two DAE formulations of the optimal control problem for the rolling ball which explicitly enforce the algebraic versor constraint on $\mathfrak{q}$ and which are mathematically equivalent to \eqref{dyn_opt_problem_ball}. In the first DAE formulation an additional control, $\dot q_0$, is added to the control $\bu$. The first DAE formulation is
\begin{equation}
\min_{a,b,\bu_1} J 
\mbox{\, s.t. \,}
\left\{
\begin{array}{ll}
\dot {\bx} = \mathbf{f}_1\left(t,\bx,\bu_1,\mu\right), \\
h_1\left(\bx\right)=1, \\
\bsigma\left(a,\bx(a),\mu\right) = \mathbf{0},\\
\bpsi\left(b,\bx(b),\mu\right) = \mathbf{0},
\end{array}
\right.
\label{dyn_opt_problem_dae1_ball}
\end{equation}
where
\revision{R1Q15}{
	\begin{equation}
	\bu_1 \equiv \begin{bmatrix} \ddot \btheta \\ \dot q_0 \end{bmatrix}, \quad \mathbf{f}_1\left(t,\bx,\bu_1,\mu\right) \equiv \begin{bmatrix} \dot \btheta \\ \bu_1  \\ \frac{1}{2} \left( \mathfrak{q} \bOm^\sharp \right)^\flat \\ \bkappa\left(t,\bx,\ddot \btheta \right) \\ \left( \left[\mathfrak{q} \bOm^\sharp \mathfrak{q}^{-1} \right]^\flat \times r \mathbf{e}_3  \right)_{12}  \end{bmatrix}, \quad \mathrm{and} \quad h_1\left(\bx\right) \equiv  \left| \mathfrak{q} \right|^2.
	\end{equation}}
In the second DAE formulation, the first component, $q_0$, of the versor $\mathfrak{q}$ is moved from the state $\bx$ to the control $\bu$ and an imitator state, $\tilde q_0$, is used to replace $q_0$ in $\bx$. $\tilde q_{a,0} = q_{a,0}$, so that with perfect integration (i.e. no numerical integration errors), $\tilde q_0(t) = q_0(t)$ for $a \le t \le b$. Defining
\begin{equation} 
\tilde {\mathfrak{q}} \equiv \begin{bmatrix} \tilde q_0 \\  \mathfrak{q}^\flat \end{bmatrix},
\end{equation}
then with perfect integration, 
\begin{equation}
\tilde {\mathfrak{q}}(t) \equiv \begin{bmatrix} \tilde q_0(t) \\  \mathfrak{q}^\flat(t) \end{bmatrix}=\begin{bmatrix} q_0(t) \\  \mathfrak{q}^\flat(t) \end{bmatrix}=\mathfrak{q}(t)
\end{equation}
for $a \le t \le b$. $\tilde q_0$ is added to the state since the final conditions require knowledge of $q_0$, which is unavailable if it has been moved to the control since the final conditions are not a function of the control. The second DAE formulation is
\begin{equation}
\min_{a,b,\bu_2} J 
\mbox{\, s.t. \,}
\left\{
\begin{array}{ll}
\dot {\bx}_2 = \mathbf{f}_2\left(t,\bx_2,\bu_2,\mu\right), \\
h_2\left(\bx_2,\bu_2\right)=1, \\
\bsigma_2\left(a,\bx_2(a),\mu\right) = \mathbf{0},\\
\bpsi_2\left(b,\bx_2(b),\mu\right) = \mathbf{0},
\end{array}
\right.
\label{dyn_opt_problem_dae2_ball}
\end{equation}
where
\revision{R1Q15}{
	\begin{equation}
	\bx_2 \equiv \begin{bmatrix} \btheta \\ \dot \btheta \\ \tilde {\mathfrak{q}} \\ \bOm \\ \bz  \end{bmatrix}, \quad
	\bu_2 \equiv \begin{bmatrix} \ddot \btheta \\ q_0 \end{bmatrix}, \quad \mathbf{f}_2\left(t,\bx_2,\bu_2,\mu\right) \equiv \begin{bmatrix} \dot \btheta \\ \ddot \btheta  \\ \frac{1}{2} \mathfrak{q} \bOm^\sharp \\ \bkappa\left(t,\bx,\ddot \btheta \right) \\ \left( \left[\mathfrak{q} \bOm^\sharp \mathfrak{q}^{-1} \right]^\flat \times r \mathbf{e}_3  \right)_{12}  \end{bmatrix}, \quad h_2\left(\bx_2,\bu_2\right) \equiv \left| \mathfrak{q} \right|^2,
	\end{equation}}
\begin{equation}
\bsigma_2\left(a,\bx_2(a),\mu\right) \equiv \begin{bmatrix} \btheta(a) - \btheta_a \\ \dot \btheta(a) - {\dot \btheta}_a \\ \tilde {\mathfrak{q}}(a)-\mathfrak{q}_a \\ \bOm(a) - \bOm_a \\ \bz(a)- \bz_a \end{bmatrix},
\end{equation}
and
\begin{equation}
\bpsi_2\left(b,\bx_2(b),\mu\right) \equiv \begin{bmatrix} \bPi \left(\left[\tilde {\mathfrak{q}}(b) \left[ \frac{1}{M} \sum_{i=0}^n m_i \bzeta_i\left(\theta_i(b)\right) \right]^\sharp \tilde {\mathfrak{q}}(b)^{-1} \right]^\flat \right) - \bDelta_b \\ \dot \btheta(b) - {\dot \btheta}_b \\ \bOm(b) - \bOm_b \\ \bz(b)- \bz_b  \end{bmatrix}.
\end{equation}
Even though both DAE formulations \eqref{dyn_opt_problem_dae1_ball} and \eqref{dyn_opt_problem_dae2_ball} are mathematically equivalent to the ODE formulation \eqref{dyn_opt_problem_ball}, the DAE formulations \eqref{dyn_opt_problem_dae1_ball} and \eqref{dyn_opt_problem_dae2_ball} tend to be numerically more stable to solve than the ODE formulation \eqref{dyn_opt_problem_ball}, as explained in Example 6.12 ``Reorientation of an Asymmetric Rigid Body'' of \cite{betts2010practical}. While the second DAE formulation \eqref{dyn_opt_problem_dae2_ball} is computationally more efficient (i.e. faster) than the first \eqref{dyn_opt_problem_dae1_ball} because it explicitly constructs the control $q_0$ rather than $\dot q_0$, the second DAE formulation \eqref{dyn_opt_problem_dae2_ball} is not as accurate as the first \eqref{dyn_opt_problem_dae1_ball}, because it only constructs an approximation, $\tilde q_0(b)$, of $q_0(b)$, which is needed for the final conditions. 

Observe that the optimal control problem encapsulated by \eqref{dyn_opt_problem_ball} ignores path inequality constraints such as $\mathbf{D}\left(t,\bx,\bu,\mu\right) \le \mathbf{0}$, where $\mathbf{D}$ is an $r \times 1$ vector-valued function. Path inequality constraints can be incorporated in \eqref{dyn_opt_problem_ball} and additional path inequality constraints can be incorporated in \eqref{dyn_opt_problem_dae1_ball} and \eqref{dyn_opt_problem_dae2_ball} as soft constraints through penalty functions in the integrand cost function $L$ or the endpoint cost function $p$.

The indirect method \cite{hull2013optimal,RBOCnumerics} is applied now to \eqref{dyn_opt_problem_ball} to construct the endpoint function and regular Hamiltonian needed to formulate an ODE TPBVP, which renders the controlled equations of motion for the rolling ball.

\paragraph{Endpoint function, Pontryagin Hamiltonian, and regular Hamiltonian.}
The endpoint function is 
\begin{equation} \label{eq_ball_endpoint_fcn}
\begin{split}
G\left(a,\bx(a),\bxi,b,\bx(b),\bnu,\mu\right) &\equiv p\left(a,\bx(a),b,\bx(b),\mu\right)+\bxi^\mathsf{T} \bsigma\left(a,\bx(a),\mu\right)+\bnu^\mathsf{T} \bpsi\left(b,\bx(b),\mu\right) \\
&=\frac{\upsilon_a}{2}\left(a-a_\mathrm{e} \right)^2+\frac{\upsilon_b}{2}\left(b-b_\mathrm{e} \right)^2+ \bxi^\mathsf{T} \begin{bmatrix} \btheta(a) - \btheta_a \\ \dot \btheta(a) - {\dot \btheta}_a \\ \mathfrak{q}(a)-\mathfrak{q}_a \\ \bOm(a) - \bOm_a \\ \bz(a)- \bz_a \end{bmatrix}\\
&\hphantom{=}+\bnu^\mathsf{T} \begin{bmatrix} \bPi \left(\left[\mathfrak{q}(b) \left[ \frac{1}{M} \sum_{i=0}^n m_i \bzeta_i\left(\theta_i(b)\right) \right]^\sharp \mathfrak{q}(b)^{-1} \right]^\flat \right) - \bDelta_b \\ \dot \btheta(b) - {\dot \btheta}_b \\ \bOm(b) - \bOm_b \\ \bz(b)- \bz_b  \end{bmatrix},
\end{split}
\end{equation}
\revision{R1Q4}{where $\bxi \in \mathbb{R}^{2n+9}$ and $\bnu \in \mathbb{R}^{n+5+c}$, $c \in \left\{0,1,2,3\right\}$, are constant Lagrange multiplier vectors enforcing the initial and final conditions, \eqref{eq_ball_initial_conds} and \eqref{eq_ball_final_conds}, respectively.} The control (Pontryagin) Hamiltonian is
\begin{equation} \label{eq_ball_ham}
\begin{split}
H\left(t,\bx,\blam,\bu,\mu\right) &\equiv L\left(t,\bx,\bu,\mu\right) + \blam^\mathsf{T} \mathbf{f}\left(t,\bx,\bu,\mu\right) \\
&= \frac{\alpha}{2} \left| \bz - \bz_\mathrm{d} \right|^2+ \frac{\beta}{2} \left| \left( \left[\mathfrak{q} \bOm^\sharp \mathfrak{q}^{-1} \right]^\flat \times r \mathbf{e}_3  \right)_{12}\right|^2  + \sum_{i=1}^n \frac{\gamma_i}{2} {\ddot \theta}_i^2+\sum_{j=1}^K V_j \left(\bz,\mu \right)+\delta\\
&\hphantom{=} + \blam^\mathsf{T} \begin{bmatrix} \dot \btheta \\ \bu  \\ \frac{1}{2} \mathfrak{q} \bOm^\sharp \\ \bkappa\left(t,\bx,\bu \right) \\ \left( \left[\mathfrak{q} \bOm^\sharp \mathfrak{q}^{-1} \right]^\flat \times r \mathbf{e}_3  \right)_{12}  \end{bmatrix},
\end{split}
\end{equation}
\revision{R1Q4}{where $\blam \in \mathbb{R}^{2n+9}$ is a time-varying Lagrange multiplier vector enforcing the dynamics \eqref{rolling_ball_opt_con_f}.} Let $\blam_{\bOm}=\begin{bmatrix} \lambda_{2n+5} & \lambda_{2n+6} & \lambda_{2n+7} \end{bmatrix}^\mathsf{T}$. Recall that
\begin{equation} 
\begin{split}
\bkappa\left(t,\bx,\bu \right) \equiv
\left[\sum_{i=0}^n m_i \widehat{\mathbf{s}_i}^2  -\inertia \right]^{-1}  \Bigg[&\bOm \times \inertia \bOm+r \tilde \bGamma \times \bGamma\\
&+ \sum_{i=0}^n m_i \mathbf{s}_i \times  \left\{ g \bGamma+ \bOm \times \left(\bOm \times \bzeta_i +2 \dot \theta_i \bzeta_i^{\prime} \right) + \dot \theta_i^2 \bzeta_i^{\dprime} + \ddot \theta_i  \bzeta_i^{\prime} \right\}  \Bigg].
\end{split}
\end{equation}
Differentiating the Hamiltonian \eqref{eq_ball_ham} with respect to the components of the control $\bu$ gives
\begin{equation} \label{eq_H_u_i_rball}
H_{u_i} = H_{\ddot{\theta}_i} = \gamma_i {\ddot \theta}_i + \lambda_{n+i} + \blam_{\bOm}^\mathsf{T} \left[\sum_{k=0}^n m_k \widehat{\mathbf{s}_k}^2  -\inertia \right]^{-1}  \left[ m_i \mathbf{s}_i \times  \bzeta_i^{\prime}  \right],
\end{equation}
\begin{equation} 
H_{u_i u_j} = H_{\ddot{\theta}_i \ddot{\theta}_j} = \gamma_i \delta_{ij},
\end{equation}
and
\begin{equation} \label{eq_Huu_ball}
H_{\bu \bu} = \diag \begin{bmatrix} \gamma_1 & \gamma_2 & \dots & \gamma_n \end{bmatrix}.
\end{equation}
By \eqref{eq_Huu_ball},  $H_{\bu \bu}>0$ iff $\gamma_i >0$ for all $1 \le i \le n$. Consequently, the optimal control problem is regular iff $\gamma_i >0$ for all $1 \le i \le n$. \revision{R1Q5}{Let us choose $\gamma_i >0$ for all $1 \le i \le n$, so that the optimal control problem is regular.} $H_{\bu}=0$ iff $H_{u_i}=0$ for all $1 \le i \le n$. From \eqref{eq_H_u_i_rball},
\begin{equation}  \label{eq_ddtheta_exp_rball}
H_{u_i} = 0 \iff  {\ddot \theta}_i = \pi_i \left(t,\bx,\blam,\mu\right) \equiv -\gamma_i^{-1} \left\{ \lambda_{n+i} + \blam_{\bOm}^\mathsf{T} \left[\sum_{k=0}^n m_k \widehat{\mathbf{s}_k}^2  -\inertia \right]^{-1}  \left[ m_i \mathbf{s}_i \times  \bzeta_i^{\prime}  \right]  \right\}.
\end{equation}
\eqref{eq_ddtheta_exp_rball} shows that $\ddot \btheta$ may be expressed as a function $\bpi$ of $\bx$, $\blam$, and $\mu$; for generality, $\bpi$ will also depend on $t$ even though in this particular example it does not. The regular Hamiltonian is
\begin{equation} \label{eq_ball_regular_Hamiltonian}
\begin{split}
\hat H\left(t,\bx,\blam,\mu\right) &\equiv H\left(t,\bx,\blam,\bpi \left(t,\bx,\blam,\mu\right),\mu\right) \\
&=  \frac{\alpha}{2} \left| \bz - \bz_\mathrm{d} \right|^2+ \frac{\beta}{2} \left| \left( \left[\mathfrak{q} \bOm^\sharp \mathfrak{q}^{-1} \right]^\flat \times r \mathbf{e}_3  \right)_{12}\right|^2 + \sum_{i=1}^n \frac{\gamma_i}{2} {\pi}_i^2 \left(t,\bx,\blam,\mu\right)+\sum_{j=1}^K V_j \left(\bz,\mu \right)+\delta\\
&\hphantom{=} + \blam^\mathsf{T} \begin{bmatrix} \dot \btheta \\ \bpi \left(t,\bx,\blam,\mu\right)  \\ \frac{1}{2} \mathfrak{q} \bOm^\sharp \\ \bkappa\left(t,\bx,\bpi \left(t,\bx,\blam,\mu\right) \right) \\ \left( \left[\mathfrak{q} \bOm^\sharp \mathfrak{q}^{-1} \right]^\flat \times r \mathbf{e}_3  \right)_{12}  \end{bmatrix}.
\end{split}
\end{equation}

\paragraph{Controlled equations of motion.}
As explained in \cite{caillau2012differential,RBOCnumerics}, one way to solve the optimal control problem \eqref{dyn_opt_problem_ball} for the rolling ball is to solve the ODE TPBVP:
\begin{equation} \label{eq_pmp_bvp_ball}
\begin{split}
\dot {\bx} &= \hat{H}_{\blam}^\mathsf{T} \left(t,\bx,\blam,\mu\right) = \hat{\mathbf{f}}\left(t,\bx,\blam,\mu\right) \equiv \mathbf{f}\left(t,\bx,\blam,\bpi\left(t,\bx,\blam,\mu\right),\mu\right), \\
\dot {\blam} &= - \hat{H}_{\bx}^\mathsf{T} \left(t,\bx,\blam,\mu\right) =-H_{\bx}^\mathsf{T}\left(t,\bx,\blam,\bpi\left(t,\bx,\blam,\mu\right),\mu\right), \\
\left. \hat{H} \right|_{t=a} = G_{a}, \quad \left. \blam \right|_{t=a} &= -G_{\bx(a)}^\mathsf{T}, \quad G_{\bxi}^\mathsf{T} = \bsigma\left(a,\bx(a),\mu\right) = \mathbf{0}, \\
\left. \hat{H} \right|_{t=b} = -G_{b}, \quad \left. \blam \right|_{t=b} &= G_{\bx(b)}^\mathsf{T}, \quad G_{\bnu}^\mathsf{T} = \bpsi\left(b,\bx(b),\mu\right) = \mathbf{0}.
\end{split}
\end{equation}
Subappendix~\ref{app_ccem_rball} derives the formulas for constructing $H_{\bx}^\mathsf{T}$.

\rem{One way to solve the optimal control problem \eqref{dyn_opt_problem_ball} for the rolling ball is to solve the ODE TPBVP given by \eqref{eq_pmp_bvp}, \eqref{eq_pmp_lbc2}, and \eqref{eq_pmp_rbc2} using the endpoint function \eqref{eq_ball_endpoint_fcn} and the regular Hamiltonian \eqref{eq_ball_regular_Hamiltonian}. In order to realize the ODE velocity function in \eqref{eq_pmp_bvp}, $\hat{\mathbf{f}}$ and $\hat{H}_{\bx}$ must be constructed. From its definition \eqref{eq_fhat_def}, $\hat{\mathbf{f}}$ is readily constructed from $\mathbf{f}$ and $\bpi$. From \eqref{eq_Hhat_x}, $\hat{H}_{\bx}$ may be obtained from $H_{\bx}$ and $\bpi$; Subappendix~\ref{app_ccem_rball} gives the formulas for constructing $H_{\bx}$, so that $\hat{H}_{\bx}$ may be realized without resorting to computational tools like symbolic or automatic differentiation.}

These controlled equations of motion for the rolling ball actuated by internal point masses that move along arbitrarily-shaped rails fixed within the ball are new and have not appeared previously in the literature, as far as we know. Note that these controlled equations of motion track the ball's orientation (i.e. the mapping from body to spatial frames) and involve many specialized parameters and terms that enable the ball to execute a wide variety of interesting and useful maneuvers. These controlled equations of motion constitute a novel contribution of this paper.

\begin{remark}[On the validity of the controlled equations of motion]{The uncontrolled equations of motion \eqref{uncon_ball_eqns_explicit_1d} are only valid if the rolling constraint is in effect, which requires that the magnitude $N$ of the normal force be positive. The formula for $N$ is derived in \cite{putkaradze2018normalpub}.  Since the controlled equations of motion \eqref{eq_pmp_bvp_ball} do not enforce $N>0$, a solution of \eqref{eq_pmp_bvp_ball} may not satisfy $N>0$, in which case the solution will be physically inconsistent with the uncontrolled equations of motion \eqref{uncon_ball_eqns_explicit_1d}. }\end{remark}

\begin{figure}[!ht] 
	\centering
       \includegraphics[scale=.35]{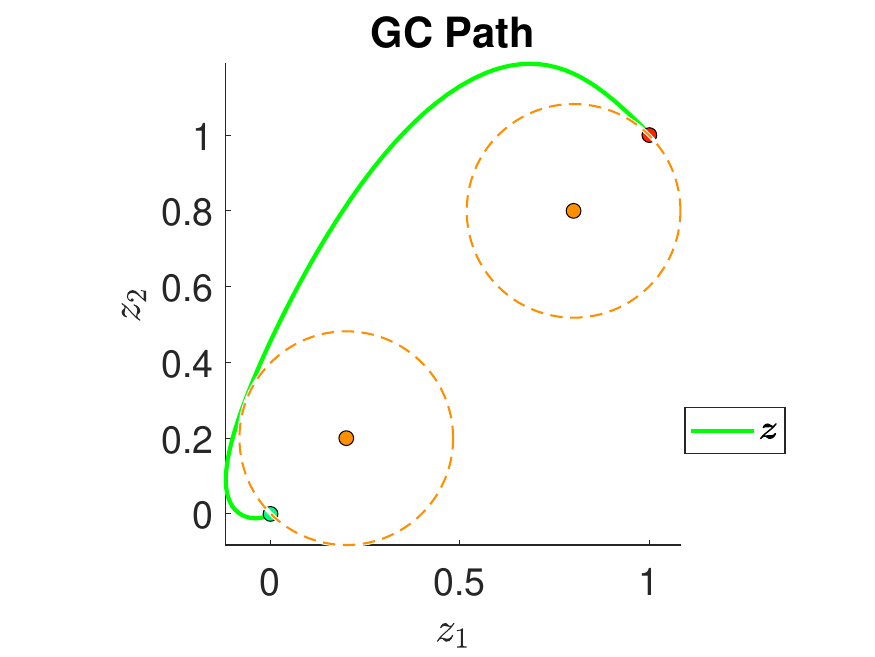} 
	\includegraphics[scale=.35]{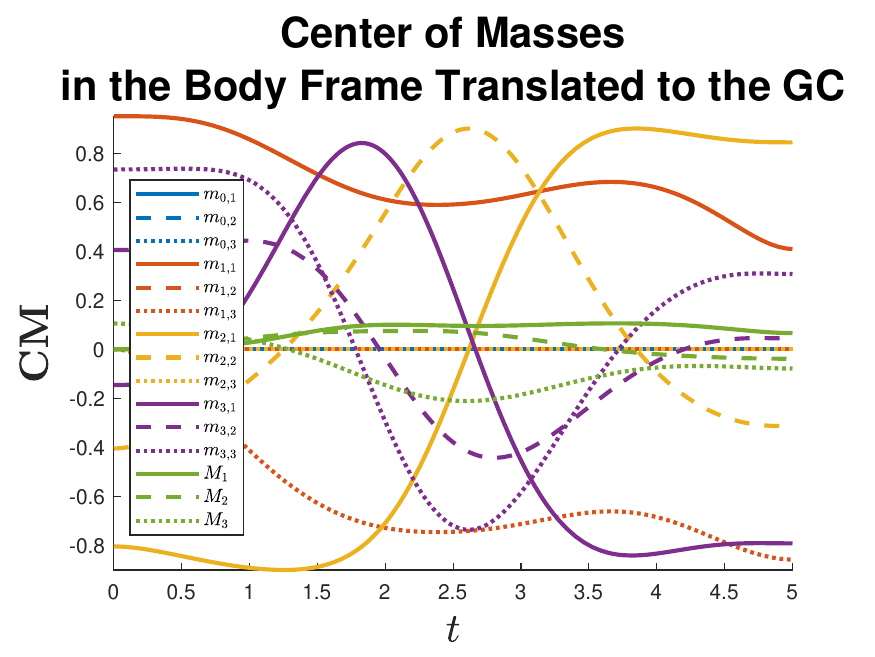} 
	\includegraphics[scale=.35]{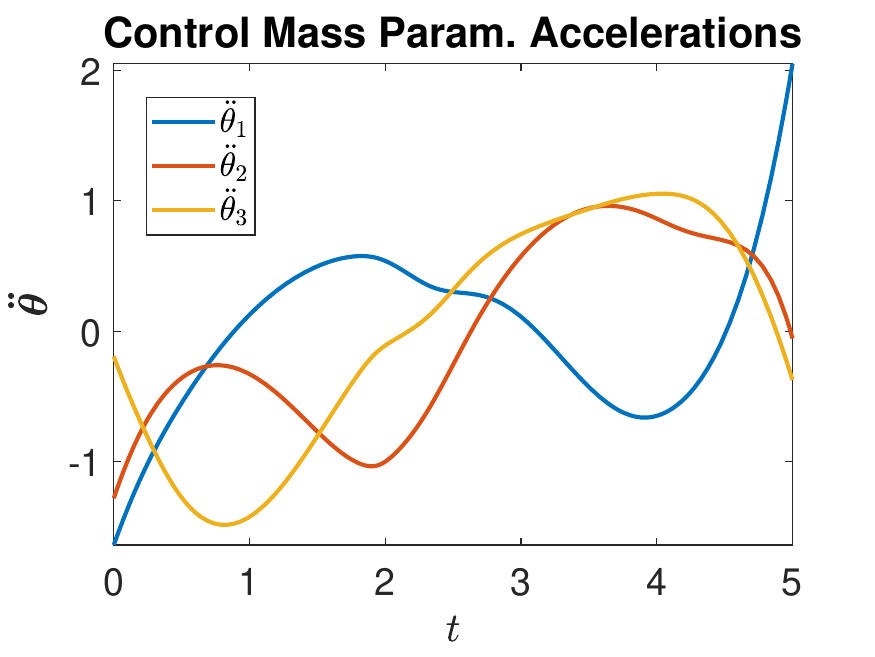}
	\caption{An example of the motion of the projection of the geometric center on the plane during obstacle avoidance (left panel), the motion of the control masses (center panel), and the motion of the controls (right panel).  }
	\label{fig_bsim10}
\end{figure}

\revisionS{R2Q3}{
\paragraph{Example numerical solution.} 
A subsequent paper \cite{RBOCnumerics} numerically solves particular cases of these controlled equations of motion.
The numerical solution poses additional technical problems which go beyond the theoretical description presented in this paper. In this paper, as an illustration of the method, we present an example in Figure~\ref{fig_bsim10} of an obstacle avoidance maneuver, together with the motion of the control masses and of the controls, i.e. the variables $\bu=\ddot \btheta$. We refer the reader to the paper \cite{RBOCnumerics} for details of the algorithm used to obtain this numerical solution. 
}

\section{Summary and Future Work} \label{sec_conclusions}

The use of rolling  ball robots in practical, real-life applications poses additional challenges. Some of the most interesting challenges that can be tackled by our method are, in our opinion,  the navigation of a  rolling ball robot over complex terrains and  imposing the  associated constraints on the  ball's motion.  Complex terrains will affect both the  uncontrolled equations of motion due to gravity and the performance index, for example, through penalty functions that discourage the ball from  ascending steep slopes. Constraints on the ball's motion may  arise during operation on a slippery surface, in which case the force at the contact  point must be bounded to avoid slippage. These constraints can be  imposed as soft constraints through penalty functions in the performance index, although it would   also  be interesting to consider the possibility of  imposing hard constraints, which represents a much more complex problem. These interesting questions will be considered in future work.

\section*{Acknowledgements}

We are indebted to our colleagues A.M. Bloch, D.M. de Diego, F. Gay-Balmaz, D.D. Holm, M. Leok, A. Lewis, T. Ohsawa, and D.V. Zenkov for useful and fruitful discussions. This research was partially  supported by  the NSERC Discovery Grant, the University of Alberta Centennial Fund, and the Alberta Innovates Technology Funding (AITF) which came through the Alberta Centre for Earth Observation Sciences (CEOS). S.M. Rogers also received support from the University of Alberta Doctoral Recruitment Scholarship, the FGSR Graduate Travel Award, the IGR Travel Award, the GSA Academic Travel Award, the AMS Fall Sectional Graduate Student Travel Grant, Target Corporation, and the Institute for Mathematics and its Applications at the University of Minnesota, Twin Cities.

\phantomsection
\addcontentsline{toc}{section}{References}
\printbibliography
\hypertarget{References}{}

\appendix

\section{Calculations for the Controlled Equations of Motion} \label{app_ccem}

This appendix derives formulas for constructing $H_{\bx}^\mathsf{T}$, which is needed to realize the controlled equations of motion \eqref{eq_pmp_bvp_disk} and \eqref{eq_pmp_bvp_ball} for a rolling disk and ball actuated by internal point masses that move along arbitrarily-shaped rails fixed within the disk and ball, respectively.

\subsection{Rolling Disk} \label{app_ccem_rdisk}

The system state $\bx$, costate $\blam$, and control $\bu$ are
\begin{equation}
\bx \equiv \begin{bmatrix} \btheta \\ \dot \btheta \\ \phi \\ \dot \phi  \end{bmatrix}, \quad \blam \equiv \begin{bmatrix} \blam_{\btheta} \\ \blam_{\dot \btheta} \\ \lambda_{\phi} \\ \lambda_{\dot \phi}  \end{bmatrix}, \quad \mathrm{and} \quad \bu \equiv \ddot \btheta,
\end{equation}
where
\begin{equation}
\btheta, \dot \btheta, \ddot \btheta \in \mathbb{R}^n \quad \mathrm{and} \quad \phi, \dot \phi  \in \mathbb{R},
\end{equation}
and
\begin{equation}
\blam_{\btheta}, \blam_{\dot \btheta}  \in \mathbb{R}^n \quad \mathrm{and} \quad \lambda_{\phi}, \lambda_{\dot \phi}  \in \mathbb{R}.
\end{equation}
The system dynamics defined for $a \le t \le b$ are
\begin{equation} \label{eq_rdisk_dynamics}
\dot {\bx} = \begin{bmatrix} \dot \btheta \\ \ddot \btheta \\ \dot \phi \\ \ddot \phi  \end{bmatrix}  = \mathbf{f}\left(t,\bx,\bu,\mu\right) \equiv \begin{bmatrix} \dot \btheta \\ \bu  \\ \dot \phi \\ \kappa\left(t,\bx,\bu\right)  \end{bmatrix},
\end{equation}
where 
\begin{equation} \label{eq_rdisk_kappa}
\kappa\left(t,\bx,\bu\right) \equiv \frac{N}{D} = \frac{ -r F_{\mathrm{e},1}+ \sum_{i=0}^n m_i K_i }{d_2+\sum_{i=0}^n m_i P_i},
\end{equation}
\begin{equation} \label{eq_rdisk_N}
N \equiv -r F_{\mathrm{e},1}+ \sum_{i=0}^n m_i K_i,
\end{equation}
\begin{equation} \label{eq_rdisk_D}
D \equiv d_2+\sum_{i=0}^n m_i P_i,
\end{equation}
\begin{equation} \label{eq_rdisk_Ki}
\begin{split}
K_i &\equiv  \left(g+ r {\dot \phi}^2 \right) \left(\zeta_{i,3} \sin \phi - \zeta_{i,1} \cos \phi  \right)+
\left(r \cos \phi + \zeta_{i,3} \right) \left(- 2 \dot \phi {\dot \theta}_i \zeta_{i,3}^{\prime} + {\dot \theta}_i^2 \zeta_{i,1}^{\dprime} + {\ddot \theta}_i \zeta_{i,1}^{\prime} \right)\\
&\hphantom{\equiv}   - \left(r \sin \phi + \zeta_{i,1} \right) \left( 2 \dot \phi {\dot \theta}_i \zeta_{i,1}^{\prime}+ {\dot \theta}_i^2 \zeta_{i,3}^{\dprime} + {\ddot \theta}_i \zeta_{i,3}^{\prime} \right), 
\end{split}
\end{equation}
and
\begin{equation} \label{eq_rdisk_Pi}
P_i \equiv \left( r \sin \phi + \zeta_{i,1} \right)^2+\left( r \cos \phi+ \zeta_{i,3} \right)^2.
\end{equation}
The Hamiltonian is
\begin{equation} \label{eq_rdisk_H}
\begin{split}
H\left(t,\bx,\blam,\bu,\mu\right) &\equiv L\left(t,\bx,\bu,\mu\right) + \blam^\mathsf{T} \mathbf{f}\left(t,\bx,\bu,\mu\right) \\
&= \frac{\alpha}{2} \left(z_a - r \left(\phi-\phi_a \right) - z_\mathrm{d} \right)^2+ \frac{\beta}{2} \left(- r \dot \phi \right)^2  + \sum_{i=1}^n \frac{\gamma_i}{2} {\ddot \theta}_i^2+\delta + \blam^\mathsf{T} \begin{bmatrix} \dot \btheta \\ \bu  \\ \dot \phi \\ \kappa\left(t,\bx,\bu \right)  \end{bmatrix} \\
&= \frac{\alpha}{2} \left(z_a - r \left(\phi-\phi_a \right) - z_\mathrm{d} \right)^2+ \frac{\beta r^2}{2} {\dot \phi}^2  + \sum_{i=1}^n \frac{\gamma_i}{2} {\ddot \theta}_i^2+\delta + \blam^\mathsf{T}_{\btheta} \dot \btheta + \blam^\mathsf{T}_{\dot \btheta} \bu \\
&\hphantom{=} + \lambda_{\phi} \dot \phi + \lambda_{\dot \phi} \kappa\left(t,\bx,\bu \right).
\end{split}
\end{equation}

In order to realize the ODE velocity function in the controlled equations of motion \eqref{eq_pmp_bvp_disk} for the rolling disk, $\hat{H}_{\bx}^\mathsf{T}$ must be constructed. $\hat{H}_{\bx}^\mathsf{T}$ may be obtained from $H_{\bx}^\mathsf{T}$ and $\bpi$, where for the rolling disk 
\begin{equation}
H_{\bx}^\mathsf{T} = \begin{bmatrix} H_{\btheta} & H_{\dot \btheta} & H_{\phi} & H_{\dot \phi} \end{bmatrix}^\mathsf{T} = \begin{bmatrix} H_{\btheta}^\mathsf{T} \\ H_{\dot \btheta}^\mathsf{T} \\ H_{\phi} \\ H_{\dot \phi} \end{bmatrix}.
\end{equation}
Formulas for constructing the components of $H_{\bx}^\mathsf{T}$, namely $H_{\btheta}^\mathsf{T}$, $H_{\dot \btheta}^\mathsf{T}$, $H_{\phi}$, and $H_{\dot \phi}$, are derived below.

\begin{equation}
\pp{H}{\btheta} = \lambda_{\dot \phi} \pp{\kappa}{\btheta},
\end{equation}

\begin{equation}
\left[\pp{H}{\btheta}\right]^\mathsf{T} = \lambda_{\dot \phi} \left[ \pp{\kappa}{\btheta}\right]^\mathsf{T},
\end{equation}

\begin{equation}
\pp{\kappa}{\theta_i} = \frac{D\pp{N}{\theta_i}-N\pp{D}{\theta_i}}{D^2} = \frac{D m_i \pp{K_i}{\theta_i}-N m_i \pp{P_i}{\theta_i}}{D^2} = \frac{m_i}{D^2} \left[ D \pp{K_i}{\theta_i}-N \pp{P_i}{\theta_i} \right],
\end{equation}

\begin{equation}
\begin{split}
\pp{K_i}{\theta_i} &= \left(g+ r {\dot \phi}^2 \right) \left(\zeta^{\prime}_{i,3} \sin \phi - \zeta^{\prime}_{i,1} \cos \phi  \right)+
\zeta^{\prime}_{i,3} \left(- 2 \dot \phi {\dot \theta}_i \zeta_{i,3}^{\prime} + {\dot \theta}_i^2 \zeta_{i,1}^{\dprime} + {\ddot \theta}_i \zeta_{i,1}^{\prime} \right) \\
&\hphantom{=}+\left(r \cos \phi + \zeta_{i,3} \right) \left(- 2 \dot \phi {\dot \theta}_i \zeta_{i,3}^{\dprime} + {\dot \theta}_i^2 \zeta_{i,1}^{\thprime} + {\ddot \theta}_i \zeta_{i,1}^{\dprime} \right) - \zeta^{\prime}_{i,1} \left( 2 \dot \phi {\dot \theta}_i \zeta_{i,1}^{\prime}+ {\dot \theta}_i^2 \zeta_{i,3}^{\dprime} + {\ddot \theta}_i \zeta_{i,3}^{\prime} \right) \\
&\hphantom{=}   - \left(r \sin \phi + \zeta_{i,1} \right) \left( 2 \dot \phi {\dot \theta}_i \zeta_{i,1}^{\dprime}+ {\dot \theta}_i^2 \zeta_{i,3}^{\thprime} + {\ddot \theta}_i \zeta_{i,3}^{\dprime} \right),
\end{split}
\end{equation}

\begin{equation}
\pp{P_i}{\theta_i} = 2 \left( r \sin \phi + \zeta_{i,1} \right)\zeta^{'}_{i,1}+2 \left( r \cos \phi+ \zeta_{i,3} \right)\zeta^{'}_{i,3} = 2 \left[\left( r \sin \phi + \zeta_{i,1} \right)\zeta^{'}_{i,1}+ \left( r \cos \phi+ \zeta_{i,3} \right)\zeta^{'}_{i,3} \right],
\end{equation}

\begin{equation}
\pp{H}{\dot \btheta} = \blam^\mathsf{T}_{\btheta} + \lambda_{\dot \phi} \pp{\kappa}{\dot \btheta},
\end{equation}

\begin{equation}
\left[\pp{H}{\dot \btheta}\right]^\mathsf{T} = \blam_{\btheta} + \lambda_{\dot \phi} \left[ \pp{\kappa}{\dot \btheta} \right]^\mathsf{T},
\end{equation}

\begin{equation}
\begin{split}
\pp{\kappa}{\dot \theta_i} = \frac{m_i}{D} \pp{K_i}{\dot \theta_i} &= \frac{m_i}{D} \left[\left(r \cos \phi + \zeta_{i,3} \right) \left(- 2 \dot \phi \zeta_{i,3}^{\prime} + 2 {\dot \theta}_i \zeta_{i,1}^{\dprime} \right) - \left(r \sin \phi + \zeta_{i,1} \right) \left(2 \dot \phi \zeta_{i,1}^{\prime} + 2 {\dot \theta}_i \zeta_{i,3}^{\dprime} \right) \right] \\
&= \frac{2 m_i}{D} \left[\left(r \cos \phi + \zeta_{i,3} \right) \left(- \dot \phi \zeta_{i,3}^{\prime} + {\dot \theta}_i \zeta_{i,1}^{\dprime} \right) - \left(r \sin \phi + \zeta_{i,1} \right) \left(\dot \phi \zeta_{i,1}^{\prime} + {\dot \theta}_i \zeta_{i,3}^{\dprime} \right) \right],
\end{split}
\end{equation}

\begin{equation}
\pp{H}{\phi} = - \alpha r \left(z_a - r \left(\phi-\phi_a \right) - z_\mathrm{d} \right)+ \lambda_{\dot \phi} \pp{\kappa}{\phi},
\end{equation}

\begin{equation}
\pp{\kappa}{\phi} = \frac{D\pp{N}{\phi}-N\pp{D}{\phi}}{D^2} = \frac{D \left( \sum_{i=0}^n m_i \pp{K_i}{\phi} \right)-N \left( \sum_{i=0}^n m_i \pp{P_i}{\phi} \right) }{D^2},
\end{equation}

\begin{equation}
\begin{split}
\pp{K_i}{\phi} &=  \left(g+ r {\dot \phi}^2 \right) \left(\zeta_{i,3} \cos \phi + \zeta_{i,1} \sin \phi  \right)+
\left(-r \sin \phi \right) \left(- 2 \dot \phi {\dot \theta}_i \zeta_{i,3}^{\prime} + {\dot \theta}_i^2 \zeta_{i,1}^{\dprime} + {\ddot \theta}_i \zeta_{i,1}^{\prime} \right) \\
&\hphantom{=} - \left(r \cos \phi \right) \left( 2 \dot \phi {\dot \theta}_i \zeta_{i,1}^{\prime}+ {\dot \theta}_i^2 \zeta_{i,3}^{\dprime} + {\ddot \theta}_i \zeta_{i,3}^{\prime} \right) \\
&= \left(g+ r {\dot \phi}^2 \right) \left(\zeta_{i,3} \cos \phi + \zeta_{i,1} \sin \phi  \right)-r \Big[
\sin \phi \left(- 2 \dot \phi {\dot \theta}_i \zeta_{i,3}^{\prime} + {\dot \theta}_i^2 \zeta_{i,1}^{\dprime} + {\ddot \theta}_i \zeta_{i,1}^{\prime} \right)  \\
&\hphantom{=\left(g+ r {\dot \phi}^2 \right) \left(\zeta_{i,3} \cos \phi + \zeta_{i,1} \sin \phi  \right)-r \Big[} 
+ \cos \phi \left( 2 \dot \phi {\dot \theta}_i \zeta_{i,1}^{\prime}+ {\dot \theta}_i^2 \zeta_{i,3}^{\dprime} + {\ddot \theta}_i \zeta_{i,3}^{\prime} \right) \Big]  , 
\end{split}
\end{equation}

\begin{equation}
\pp{P_i}{\phi} = 2 \left( r \sin \phi + \zeta_{i,1} \right) \left( r \cos \phi \right)+2 \left( r \cos \phi+ \zeta_{i,3} \right) \left(-r \sin \phi \right)=2r \left( \zeta_{i,1} \cos \phi - \zeta_{i,3} \sin \phi \right),
\end{equation}

\begin{equation}
\pp{H}{\dot \phi} = \beta r^2 \dot \phi+\lambda_{\phi} + \lambda_{\dot \phi} \pp{\kappa}{\dot \phi},
\end{equation}

\begin{equation} 
\begin{split}
\pp{\kappa}{\dot \phi} &= \frac{1}{D}\sum_{i=0}^n m_i \pp{K_i}{\dot \phi} \\
&= \frac{1}{D}\sum_{i=0}^n m_i \left[ \left(2 r {\dot \phi} \right) \left(\zeta_{i,3} \sin \phi - \zeta_{i,1} \cos \phi  \right)-2 \dot \theta_i \left\{
\left(r \cos \phi + \zeta_{i,3} \right) \zeta_{i,3}^{\prime}  + \left(r \sin \phi + \zeta_{i,1} \right) \zeta_{i,1}^{\prime} \right\} \right] \\
&= \frac{2}{D}\sum_{i=0}^n m_i \left[ r {\dot \phi} \left(\zeta_{i,3} \sin \phi - \zeta_{i,1} \cos \phi  \right)- \dot \theta_i \left\{
\left(r \cos \phi + \zeta_{i,3} \right) \zeta_{i,3}^{\prime}  + \left(r \sin \phi + \zeta_{i,1} \right) \zeta_{i,1}^{\prime} \right\} \right].
\end{split}
\end{equation}

\subsection{Rolling Ball} \label{app_ccem_rball}
The system state $\bx$, costate $\blam$, and control $\bu$ are
\begin{equation}
\bx \equiv \begin{bmatrix} \btheta \\ \dot \btheta \\ \mathfrak{q} \\ \bOm \\ \bz  \end{bmatrix}, \quad \blam \equiv \begin{bmatrix} \blam_{\btheta} \\ \blam_{\dot \btheta} \\ \blam_{\mathfrak{q}} \\ \blam_{\bOm} \\ \blam_{\bz}  \end{bmatrix}, \quad \mathrm{and} \quad \bu \equiv \ddot \btheta,
\end{equation}
where
\begin{equation}
\btheta, \dot \btheta, \ddot \btheta \in \mathbb{R}^n, \quad \mathfrak{q} = \begin{bmatrix} q_0 \\ q_1 \\ q_2 \\ q_3 \end{bmatrix} \in \mathscr{S} \cong \mathbb{S}^3 \subset \mathbb{R}^4, \quad \bOm = \begin{bmatrix} \Omega_1 \\ \Omega_2 \\ \Omega_3 \end{bmatrix} \in \mathbb{R}^3, \quad \mathrm{and} \quad \bz = \begin{bmatrix} z_1 \\ z_2 \end{bmatrix} \in \mathbb{R}^2,
\end{equation}
and
\begin{equation}
\blam_{\btheta}, \blam_{\dot \btheta} \in \mathbb{R}^n, \quad \blam_{\mathfrak{q}}  \in \mathbb{R}^4, \quad \blam_{\bOm} \in \mathbb{R}^3, \quad \mathrm{and} \quad \blam_{\bz} \in \mathbb{R}^2.
\end{equation}
Recall that given a versor 
\begin{equation}
\mathfrak{q} = \begin{bmatrix}q_0 \\ q_1 \\ q_2 \\ q_3 \end{bmatrix} \in \mathscr{S},
\end{equation}
the corresponding rotation matrix $\Lambda \in SO(3)$ is
\begin{equation} \label{eq_versor_rot}
\Lambda = \begin{bmatrix} 1-2\left(q_2^2+q_3^2\right) & 2\left(q_1q_2-q_0q_3 \right) & 2\left(q_1 q_3+q_0q_2 \right) \\ 2\left(q_1 q_2+q_0q_3 \right) & 1-2\left(q_1^2+q_3^2\right) &  2\left(q_2 q_3-q_0q_1 \right) \\ 2\left(q_1 q_3-q_0q_2 \right) & 2\left(q_2 q_3+q_0q_1 \right) & 1-2\left(q_1^2+q_2^2\right) \end{bmatrix}.
\end{equation}
The system dynamics defined for $a \le t \le b$ are
\begin{equation} \label{eq_rball_dynamics}
\dot {\bx} = \begin{bmatrix} \dot \btheta \\ \ddot \btheta \\ \dot {\mathfrak{q}} \\ \dot \bOm \\ \dot \bz  \end{bmatrix}  = \mathbf{f}\left(t,\bx,\bu,\mu\right) \equiv \begin{bmatrix} \dot \btheta \\ \bu  \\ \frac{1}{2} \mathfrak{q} \bOm^\sharp \\ \bkappa\left(t,\bx,\bu \right) \\ \left( \left[\mathfrak{q} \bOm^\sharp \mathfrak{q}^{-1} \right]^\flat \times r \mathbf{e}_3  \right)_{12}  \end{bmatrix},
\end{equation}
where
\begin{equation} \label{eq_rball_kappa}
\begin{split}
\bkappa\left(t,\bx,\bu \right) \equiv
\left[\sum_{i=0}^n m_i \widehat{\mathbf{s}_i}^2  -\inertia \right]^{-1}  \Bigg[&\bOm \times \inertia \bOm+r \tilde \bGamma \times \bGamma\\
&+ \sum_{i=0}^n m_i \mathbf{s}_i \times  \left\{ g \bGamma+ \bOm \times \left(\bOm \times \bzeta_i +2 \dot \theta_i \bzeta_i^{\prime} \right) + \dot \theta_i^2 \bzeta_i^{\dprime} + \ddot \theta_i  \bzeta_i^{\prime} \right\}  \Bigg],
\end{split}
\end{equation}
\begin{equation} \label{eq_rball_Gamma}
\bGam \equiv \Lambda^{-1} \mathbf{e}_3 = \Lambda^\mathsf{T}\mathbf{e}_3 = \begin{bmatrix}  2\left(q_1 q_3-q_0q_2 \right) \\ 2\left(q_2 q_3+q_0q_1 \right) \\ 1-2\left(q_1^2+q_2^2\right) \end{bmatrix}, 
\end{equation}
and
\begin{equation} \label{eq_rball_tGamma}
\tilde \bGamma \equiv \Lambda^{-1} \mathbf{F}_\mathrm{e} = \Lambda^\mathsf{T} \mathbf{F}_\mathrm{e} = \begin{bmatrix} \left[1-2\left(q_2^2+q_3^2\right)\right]\mathbf{F}_\mathrm{e,1} + 2\left(q_1q_2+q_0q_3 \right)\mathbf{F}_\mathrm{e,2} + 2\left(q_1 q_3-q_0q_2 \right)\mathbf{F}_\mathrm{e,3} \\ 2\left(q_1 q_2-q_0q_3 \right)\mathbf{F}_\mathrm{e,1} + \left[1-2\left(q_1^2+q_3^2\right)\right]\mathbf{F}_\mathrm{e,2} +  2\left(q_2 q_3
+q_0q_1 \right)\mathbf{F}_\mathrm{e,3} \\ 2\left(q_1 q_3+q_0q_2 \right)\mathbf{F}_\mathrm{e,1} + 2\left(q_2 q_3-q_0q_1 \right)\mathbf{F}_\mathrm{e,2} + \left[1-2\left(q_1^2+q_2^2\right)\right] \mathbf{F}_\mathrm{e,3} \end{bmatrix}.
\end{equation}

For $1 \le j \le K$, let $\bv_j \in \mathbb{R}^2$ denote the $j^\mathrm{th}$ obstacle center, $h_j \in \mathbb{R}^{+}$ denote the $j^\mathrm{th}$ obstacle height, $\rho_j \in \mathbb{R}^{+}$ denote the $j^\mathrm{th}$ obstacle radius, and 
\begin{equation}
V_j \left(\bz,\mu \right) \equiv h_j \, S \left( \left| \bz - \bv_j \right| - \rho_j \right)
\end{equation}
denote the $j^\mathrm{th}$ obstacle avoidance potential function where
\begin{equation}
S(y)\equiv \frac{1}{2}+\frac{1}{2} \tanh \left(\frac{-y}{\epsilon} \right)
\end{equation}
is a time-reversed sigmoid function or
\begin{equation}
S(y)\equiv \left[\max \left\{0,-y \right\} \right]^4  =\ReLU^4\left(-y\right)
\end{equation}
is a $C^2$ cutoff function.

The Hamiltonian is 
\begin{equation} \label{eq_rball_H}
\begin{split}
H\left(t,\bx,\blam,\bu,\mu\right) &\equiv L\left(t,\bx,\bu,\mu\right) + \blam^\mathsf{T} \mathbf{f}\left(t,\bx,\bu,\mu\right) \\
&= \frac{\alpha}{2} \left| \bz - \bz_\mathrm{d} \right|^2+ \frac{\beta}{2} \left| \left( \left[\mathfrak{q} \bOm^\sharp \mathfrak{q}^{-1} \right]^\flat \times r \mathbf{e}_3  \right)_{12}\right|^2  + \sum_{i=1}^n \frac{\gamma_i}{2} {\ddot \theta}_i^2+\sum_{j=1}^K V_j \left(\bz,\mu \right)+\delta\\
&\hphantom{=} + \blam^\mathsf{T} \begin{bmatrix} \dot \btheta \\ \bu  \\ \frac{1}{2} \mathfrak{q} \bOm^\sharp \\ \bkappa\left(t,\bx,\bu \right) \\ \left( \left[\mathfrak{q} \bOm^\sharp \mathfrak{q}^{-1} \right]^\flat \times r \mathbf{e}_3  \right)_{12}  \end{bmatrix} \\
&= \frac{\alpha}{2} \left| \bz - \bz_\mathrm{d} \right|^2+ \frac{\beta}{2} \left| \left( \left[\mathfrak{q} \bOm^\sharp \mathfrak{q}^{-1} \right]^\flat \times r \mathbf{e}_3  \right)_{12}\right|^2  + \sum_{i=1}^n \frac{\gamma_i}{2} {\ddot \theta}_i^2+\sum_{j=1}^K V_j \left(\bz,\mu \right)+\delta\\
&\hphantom{=} + \blam^\mathsf{T}_{\btheta} \dot \btheta +\blam^\mathsf{T}_{\dot \btheta} \bu + \blam^\mathsf{T}_{\mathfrak{q}} \left( \frac{1}{2} \mathfrak{q} \bOm^\sharp \right) + \blam^\mathsf{T}_{\bOm} \bkappa\left(t,\bx,\bu \right) + \blam^\mathsf{T}_{\bz} \left( \left[\mathfrak{q} \bOm^\sharp \mathfrak{q}^{-1} \right]^\flat \times r \mathbf{e}_3  \right)_{12} \\
&= \frac{\alpha}{2} \left| \bz - \bz_\mathrm{d} \right|^2+ \frac{\beta r^2}{2} \left(\tilde{Q}\bOm\right)^\mathsf{T} \left(\tilde{Q}\bOm\right)  + \sum_{i=1}^n \frac{\gamma_i}{2} {\ddot \theta}_i^2+\sum_{j=1}^K V_j \left(\bz,\mu \right)+\delta\\
&\hphantom{=} + \blam^\mathsf{T}_{\btheta} \dot \btheta +\blam^\mathsf{T}_{\dot \btheta} \bu + \frac{1}{2} \blam^\mathsf{T}_{\mathfrak{q}} \left(Q \bOm\right) + \blam^\mathsf{T}_{\bOm} \bkappa\left(t,\bx,\bu \right) + r\blam^\mathsf{T}_{\bz} \left(\tilde{Q}\bOm\right), 
\end{split}
\end{equation}
since
\begin{equation}
\frac{1}{2} \mathfrak{q} \bOm^\sharp = \frac{1}{2} \begin{bmatrix} -\Omega_1 q_1-\Omega_2 q_2 - \Omega_3 q_3 \\ \Omega_1 q_0+\Omega_3 q_2 - \Omega_2 q_3 \\ \Omega_2 q_0+\Omega_1 q_3 - \Omega_3 q_1 \\ \Omega_3 q_0+\Omega_2 q_1 - \Omega_1 q_2 \end{bmatrix} = \frac{1}{2} \begin{bmatrix} -q_1 & -q_2 & -q_3 \\ q_0 & -q_3  & q_2 \\ q_3 & q_0 & - q_1 \\ -q_2 & q_1 & q_0 \end{bmatrix}  \bOm = \frac{1}{2} Q \bOm
\end{equation}
and
\begin{equation}
\begin{split}
\left( \left[\mathfrak{q} \bOm^\sharp \mathfrak{q}^{-1} \right]^\flat \times r \mathbf{e}_3  \right)_{12} &= \left( \Lambda \bOm \times r \mathbf{e}_3  \right)_{12} = \left( \bomega \times r \mathbf{e}_3  \right)_{12} = r \begin{bmatrix} \omega_2 \\ -\omega_1 \end{bmatrix}  \\
&=  r \begin{bmatrix} 2\left(q_1 q_2+q_0q_3 \right)\Omega_1 + \left[1-2\left(q_1^2+q_3^2\right) \right] \Omega_2 +  2\left(q_2 q_3-q_0q_1 \right)\Omega_3 \\ -\left[1-2\left(q_2^2+q_3^2\right)\right]\Omega_1  -2\left(q_1q_2-q_0q_3 \right)\Omega_2 -2\left(q_1 q_3+q_0q_2 \right)\Omega_3 \end{bmatrix} \\
&= r \begin{bmatrix} 2\left(q_1 q_2+q_0q_3 \right) & 1-2\left(q_1^2+q_3^2\right) &  2\left(q_2 q_3-q_0q_1 \right) \\ -\left[1-2\left(q_2^2+q_3^2\right)\right] & -2\left(q_1q_2-q_0q_3 \right) & -2\left(q_1 q_3+q_0q_2 \right) \end{bmatrix} \begin{bmatrix} \Omega_1 \\ \Omega_2 \\ \Omega_3 \end{bmatrix} = r \tilde{Q} \bOm,
\end{split}
\end{equation}
where
\begin{equation}
Q \equiv \begin{bmatrix} -q_1 & -q_2 & -q_3 \\ q_0 & -q_3  & q_2 \\ q_3 & q_0 & - q_1 \\ -q_2 & q_1 & q_0 \end{bmatrix}
\end{equation}
and
\begin{equation}
\tilde{Q} \equiv \begin{bmatrix} 2\left(q_1 q_2+q_0q_3 \right) & 1-2\left(q_1^2+q_3^2\right) &  2\left(q_2 q_3-q_0q_1 \right) \\ -\left[1-2\left(q_2^2+q_3^2\right)\right] & -2\left(q_1q_2-q_0q_3 \right) & -2\left(q_1 q_3+q_0q_2 \right) \end{bmatrix}.
\end{equation}

In order to realize the ODE velocity function in the controlled equations of motion \eqref{eq_pmp_bvp_ball} for the rolling ball, $\hat{H}_{\bx}^\mathsf{T}$ must be constructed. $\hat{H}_{\bx}^\mathsf{T}$ may be obtained from $H_{\bx}^\mathsf{T}$ and $\bpi$, where for the rolling ball 
\begin{equation}
H_{\bx}^\mathsf{T} = \begin{bmatrix} H_{\btheta} & H_{\dot \btheta} & H_{\mathfrak{q}} & H_{\bOm} & H_{\bz} \end{bmatrix}^\mathsf{T} = \begin{bmatrix} H_{\btheta}^\mathsf{T} \\ H_{\dot \btheta}^\mathsf{T} \\ H_{\mathfrak{q}}^\mathsf{T} \\ H_{\bOm}^\mathsf{T} \\ H_{\bz}^\mathsf{T} \end{bmatrix} .
\end{equation}
Formulas for constructing the components of $H_{\bx}^\mathsf{T}$, namely $H_{\btheta}^\mathsf{T}$, $H_{\dot \btheta}^\mathsf{T}$, $H_{\mathfrak{q}}^\mathsf{T}$, $H_{\bOm}^\mathsf{T}$, and $H_{\bz}^\mathsf{T}$, are derived below.

\begin{equation}
\pp{H}{\btheta} = \blam_{ 
	\bOm}^\mathsf{T} \pp{\bkappa}{\btheta},
\end{equation}
so that
\begin{equation}
\left[\pp{H}{\btheta}\right]^\mathsf{T} =\left[\pp{\bkappa}{\btheta}\right]^\mathsf{T} \blam_{\bOm}.
\end{equation}

\begin{equation}
\pp{H}{\dot \btheta} = \blam^\mathsf{T}_{\btheta} + \blam_{ 
	\bOm}^\mathsf{T} \pp{\bkappa}{\dot \btheta},
\end{equation}
so that
\begin{equation}
\left[\pp{H}{\dot \btheta}\right]^\mathsf{T} = \blam_{\btheta} +  \left[\pp{\bkappa}{\dot \btheta}\right]^\mathsf{T} \blam_{ 
	\bOm}.
\end{equation}

If 
\begin{equation}
S(y)\equiv \frac{1}{2}+\frac{1}{2} \tanh \left(\frac{-y}{\epsilon} \right),
\end{equation}
then
\begin{equation}
\begin{split}
\pp{H}{\bz} &= \alpha \left( \bz - \bz_\mathrm{d} \right)^\mathsf{T}+\sum_{j=1}^K \pp {V_j \left(\bz,\mu \right)}{\bz} \\
&= \alpha \left( \bz - \bz_\mathrm{d} \right)^\mathsf{T}-\frac{1}{2 \epsilon}\sum_{j=1}^K h_j \left[1-\tanh^2 \left(\frac{\rho_j-\left| \bz - \bv_j \right|}{\epsilon} \right) \right] \frac{\left( \bz - \bv_j \right)^\mathsf{T}}{\left| \bz - \bv_j \right|},
\end{split}
\end{equation}
since
\begin{equation}
\pp {V_j \left(\bz,\mu \right)}{\bz} = h_j \, S^\prime \left(\left| \bz - \bv_j \right|-\rho_j \right) \frac{\left( \bz - \bv_j \right)^\mathsf{T}}{\left| \bz - \bv_j \right|}  = -\frac{h_j}{2 \epsilon}  \left[1-\tanh^2 \left(\frac{\rho_j-\left| \bz - \bv_j \right|}{\epsilon} \right) \right] \frac{\left( \bz - \bv_j \right)^\mathsf{T}}{\left| \bz - \bv_j \right|}
\end{equation}
and
\begin{equation}
S^\prime(y)= -\frac{1}{2 \epsilon} \left[1-\tanh^2 \left(\frac{-y}{\epsilon} \right) \right].
\end{equation}
Therefore,
\begin{equation}
\left[ \pp{H}{\bz} \right]^\mathsf{T} = \alpha \left( \bz - \bz_\mathrm{d} \right)-\frac{1}{2 \epsilon}\sum_{j=1}^K h_j \left[1-\tanh^2 \left(\frac{\rho_j-\left| \bz - \bv_j \right|}{\epsilon} \right) \right] \frac{\left( \bz - \bv_j \right)}{\left| \bz - \bv_j \right|}.
\end{equation}

If 
\begin{equation}
S(y)\equiv \left[\max \left\{0,-y \right\} \right]^4=\ReLU^4\left(-y\right),
\end{equation}
then
\begin{equation}
\begin{split}
\pp{H}{\bz} &= \alpha \left( \bz - \bz_\mathrm{d} \right)^\mathsf{T}+\sum_{j=1}^K \pp {V_j \left(\bz,\mu \right)}{\bz} \\
&= \alpha \left( \bz - \bz_\mathrm{d} \right)^\mathsf{T}-4\sum_{j=1}^K h_j \left[\max \left\{0,\rho_j-\left| \bz - \bv_j \right| \right\} \right]^3 \frac{\left( \bz - \bv_j \right)^\mathsf{T}}{\left| \bz - \bv_j \right|},
\end{split}
\end{equation}
since
\begin{equation}
\pp {V_j \left(\bz,\mu \right)}{\bz} = h_j \, S^\prime \left(\left| \bz - \bv_j \right|-\rho_j \right) \frac{\left( \bz - \bv_j \right)^\mathsf{T}}{\left| \bz - \bv_j \right|}  = -4 h_j \left[\max \left\{0,\rho_j-\left| \bz - \bv_j \right| \right\} \right]^3 \frac{\left( \bz - \bv_j \right)^\mathsf{T}}{\left| \bz - \bv_j \right|}
\end{equation}
and
\begin{equation}
S^\prime(y)= -4 \left[\max \left\{0,-y\right\} \right]^3.
\end{equation}
Therefore,
\begin{equation}
\left[ \pp{H}{\bz} \right]^\mathsf{T} = \alpha \left( \bz - \bz_\mathrm{d} \right)-4 \sum_{j=1}^K h_j \left[\max \left\{0,\rho_j-\left| \bz - \bv_j \right| \right\} \right]^3 \frac{\left( \bz - \bv_j \right)}{\left| \bz - \bv_j \right|}.
\end{equation}

\begin{equation}
\pp{H}{\bOm} = \beta r^2 \left(\tilde{Q} \bOm \right)^\mathsf{T} \tilde{Q}+\frac{1}{2} \blam^\mathsf{T}_{\mathfrak{q}} Q+\blam^\mathsf{T}_{\bOm} \pp{\bkappa}{\bOm}+r \blam^\mathsf{T}_{\bz} \tilde{Q} = r \left[\beta r \left(\tilde{Q} \bOm \right)^\mathsf{T} +\blam^\mathsf{T}_{\bz}\right] \tilde{Q}+\frac{1}{2} \blam^\mathsf{T}_{\mathfrak{q}} Q+\blam^\mathsf{T}_{\bOm} \pp{\bkappa}{\bOm},
\end{equation}
so that
\begin{equation}
\begin{split}
\left[ \pp{H}{\bOm} \right]^\mathsf{T} &= r \tilde{Q}^\mathsf{T} \left[\beta r \left(\tilde{Q} \bOm \right) +\blam_{\bz}\right] +\frac{1}{2} Q^\mathsf{T} \blam_{\mathfrak{q}} + \left[\pp{\bkappa}{\bOm} \right]^\mathsf{T} \blam_{\bOm} \\
&= r \begin{bmatrix} 2\left(q_1 q_2+q_0q_3 \right) & -\left[1-2\left(q_2^2+q_3^2\right)\right] \\ 1-2\left(q_1^2+q_3^2\right) & -2\left(q_1q_2-q_0q_3 \right) \\ 2\left(q_2 q_3-q_0q_1 \right) & -2\left(q_1 q_3+q_0q_2 \right)   \end{bmatrix} \left[\beta r \begin{bmatrix} 2\left(q_1 q_2+q_0q_3 \right)\Omega_1 + \left[1-2\left(q_1^2+q_3^2\right) \right] \Omega_2 +  2\left(q_2 q_3-q_0q_1 \right)\Omega_3 \\ -\left[1-2\left(q_2^2+q_3^2\right)\right]\Omega_1  -2\left(q_1q_2-q_0q_3 \right)\Omega_2 -2\left(q_1 q_3+q_0q_2 \right)\Omega_3 \end{bmatrix} +\blam_{\bz}\right] \\
&\hphantom{=} +\frac{1}{2} \begin{bmatrix} -q_1 & q_0 & q_3 & -q_2 \\ -q_2 & -q_3  & q_0 & q_1 \\ -q_3 & q_2 & - q_1 & q_0 \end{bmatrix} \blam_{\mathfrak{q}} + \left[\pp{\bkappa}{\bOm} \right]^\mathsf{T} \blam_{\bOm}.
\end{split}
\end{equation}

\begin{equation}
\begin{split}
\pp{H}{\mathfrak{q}} &= \frac{\beta r^2}{2} \pp{}{\mathfrak{q}} \left[ \left(\tilde{Q}\bOm\right)^\mathsf{T} \left(\tilde{Q}\bOm\right) \right]  
+ \frac{1}{2} \blam^\mathsf{T}_{\mathfrak{q}} \pp{}{\mathfrak{q}} \left(Q \bOm\right) + \blam^\mathsf{T}_{\bOm} \pp{\bkappa}{\mathfrak{q}} + r\blam^\mathsf{T}_{\bz} \pp{}{\mathfrak{q}} \left(\tilde{Q}\bOm\right) \\
&= \beta r^2  \left[ \left(\tilde{Q}\bOm\right)^\mathsf{T} \pp{}{\mathfrak{q}} \left(\tilde{Q}\bOm\right) \right]  
+ \frac{1}{2} \blam^\mathsf{T}_{\mathfrak{q}} \pp{}{\mathfrak{q}} \left(Q \bOm\right) + \blam^\mathsf{T}_{\bOm} \pp{\bkappa}{\mathfrak{q}} + r\blam^\mathsf{T}_{\bz} \pp{}{\mathfrak{q}} \left(\tilde{Q}\bOm\right) \\
&=  \left[ \beta r^2 \left(\tilde{Q}\bOm\right)^\mathsf{T}+r\blam^\mathsf{T}_{\bz} \right] \pp{}{\mathfrak{q}} \left(\tilde{Q}\bOm\right)  
+ \frac{1}{2} \blam^\mathsf{T}_{\mathfrak{q}} \pp{}{\mathfrak{q}} \left(Q \bOm\right) + \blam^\mathsf{T}_{\bOm} \pp{\bkappa}{\mathfrak{q}} \\
&= r \left[ \beta r \left(\tilde{Q}\bOm\right)^\mathsf{T}+\blam^\mathsf{T}_{\bz} \right] \pp{}{\mathfrak{q}} \left(\tilde{Q}\bOm\right)  
+ \frac{1}{2} \blam^\mathsf{T}_{\mathfrak{q}} \pp{}{\mathfrak{q}} \left(Q \bOm\right) + \blam^\mathsf{T}_{\bOm} \pp{\bkappa}{\mathfrak{q}},
\end{split}
\end{equation}
so that
\begin{equation}
\begin{split}
\left[ \pp{H}{\mathfrak{q}} \right]^\mathsf{T} &=r \left[ \pp{}{\mathfrak{q}} \left(\tilde{Q}\bOm\right) \right]^\mathsf{T} \left[ \beta r \left(\tilde{Q}\bOm\right)+\blam_{\bz} \right]   
+ \frac{1}{2} \left[ \pp{}{\mathfrak{q}} \left(Q \bOm\right) \right]^\mathsf{T} \blam_{\mathfrak{q}} + \left[ \pp{\bkappa}{\mathfrak{q}} \right]^\mathsf{T} \blam_{\bOm} \\
&=2 r \begin{bmatrix} q_3 \Omega_1-q_1 \Omega_3 & q_3 \Omega_2-q_2\Omega_3 \\ q_2 \Omega_1-2q_1 \Omega_2-q_0 \Omega_3 & -q_2 \Omega_2 -q_3 \Omega_3 \\ q_1 \Omega_1+q_3 \Omega_3 & 2q_2 \Omega_1-q_1\Omega_2-q_0 \Omega_3 \\ q_0 \Omega_1-2q_3\Omega_2+q_2\Omega_3 & 2q_3\Omega_1+q_0\Omega_2-q_1\Omega_3 \end{bmatrix} \cdot \\
&\hphantom{=} \left[ \beta r \begin{bmatrix} 2\left(q_1 q_2+q_0q_3 \right)\Omega_1 + \left[1-2\left(q_1^2+q_3^2\right) \right] \Omega_2 +  2\left(q_2 q_3-q_0q_1 \right)\Omega_3 \\ -\left[1-2\left(q_2^2+q_3^2\right)\right]\Omega_1  -2\left(q_1q_2-q_0q_3 \right)\Omega_2 -2\left(q_1 q_3+q_0q_2 \right)\Omega_3 \end{bmatrix}+\blam_{\bz} \right] \\
&\hphantom{=}   
+ \frac{1}{2} \begin{bmatrix} 0 & \Omega_1 & \Omega_2 & \Omega_3 \\ -\Omega_1 & 0 & -\Omega_3 & \Omega_2 \\ -\Omega_2 & \Omega_3 & 0 & -\Omega_1 \\ -\Omega_3 & -\Omega_2 & \Omega_1 & 0 \end{bmatrix} \blam_{\mathfrak{q}} + \left[ \pp{\bkappa}{\mathfrak{q}} \right]^\mathsf{T} \blam_{\bOm}, 
\end{split}
\end{equation}
since
\begin{equation}
\pp{}{\mathfrak{q}} \left(Q \bOm\right) = \begin{bmatrix} 0 & -\Omega_1 & -\Omega_2 & -\Omega_3 \\ \Omega_1 & 0 & \Omega_3 & -\Omega_2 \\ \Omega_2 & -\Omega_3 & 0 & \Omega_1 \\ \Omega_3 & \Omega_2 & -\Omega_1 & 0 \end{bmatrix}
\end{equation}
and
\begin{equation}
\begin{split}
\pp{}{\mathfrak{q}} \left(\tilde{Q}\bOm\right) &= 
\begin{bmatrix} 2q_3 \Omega_1-2q_1 \Omega_3 & 2q_2 \Omega_1-4q_1 \Omega_2-2q_0 \Omega_3 & 2q_1 \Omega_1+2q_3 \Omega_3 & 2q_0 \Omega_1-4q_3\Omega_2+2q_2\Omega_3 \\ 2q_3 \Omega_2-2q_2\Omega_3 & -2q_2 \Omega_2 -2q_3 \Omega_3 & 4q_2 \Omega_1-2q_1\Omega_2-2q_0 \Omega_3 & 4 q_3\Omega_1+2q_0\Omega_2-2q_1\Omega_3 \end{bmatrix} \\
&= 2 \begin{bmatrix} q_3 \Omega_1-q_1 \Omega_3 & q_2 \Omega_1-2q_1 \Omega_2-q_0 \Omega_3 & q_1 \Omega_1+q_3 \Omega_3 & q_0 \Omega_1-2q_3\Omega_2+q_2\Omega_3 \\ q_3 \Omega_2-q_2\Omega_3 & -q_2 \Omega_2 -q_3 \Omega_3 & 2q_2 \Omega_1-q_1\Omega_2-q_0 \Omega_3 & 2q_3\Omega_1+q_0\Omega_2-q_1\Omega_3 \end{bmatrix}.
\end{split}
\end{equation}

\begin{equation}
\begin{split}
\pp{\bkappa}{\theta_j} &= -\left[\sum_{i=0}^n m_i \widehat{\mathbf{s}_i}^2  -\inertia \right]^{-1} \left[ m_j \left(\widehat{\mathbf{s}_j} \widehat{{\bzeta}^{\prime}_j} + \widehat{{\bzeta}^{\prime}_j} \widehat{\mathbf{s}_j} \right) \right] \left[\sum_{i=0}^n m_i \widehat{\mathbf{s}_i}^2  -\inertia \right]^{-1} \left[\sum_{i=0}^n m_i \widehat{\mathbf{s}_i}^2  -\inertia \right] \bkappa\left(t,\bx,\bu \right) \\
&\hphantom{=} + \left[\sum_{i=0}^n m_i \widehat{\mathbf{s}_i}^2  -\inertia \right]^{-1}  \Bigg[ m_j \bzeta^{\prime}_j \times  \left\{ g \bGamma+ \bOm \times \left(\bOm \times \bzeta_j +2 \dot \theta_j \bzeta_j^{\prime} \right) + \dot \theta_j^2 \bzeta_j^{\dprime} + \ddot \theta_j  \bzeta_j^{\prime} \right\} \\
&\hphantom{=\left[\sum_{i=0}^n m_i \widehat{\mathbf{s}_i}^2  -\inertia \right]^{-1}  \Bigg[} + m_j \mathbf{s}_j \times  \left\{ \bOm \times \left(\bOm \times \bzeta_j^{\prime} +2 \dot \theta_j \bzeta_j^{\dprime} \right) + \dot \theta_j^2 \bzeta_j^{\thprime} + \ddot \theta_j  \bzeta_j^{\dprime} \right\}  \Bigg] \\
&= \left[\sum_{i=0}^n m_i \widehat{\mathbf{s}_i}^2  -\inertia \right]^{-1} m_j   \Bigg[ \bzeta^{\prime}_j \times  \left\{ g \bGamma+ \bOm \times \left(\bOm \times \bzeta_j +2 \dot \theta_j \bzeta_j^{\prime} \right) + \dot \theta_j^2 \bzeta_j^{\dprime} + \ddot \theta_j  \bzeta_j^{\prime} \right\} \\
&\hphantom{=\left[\sum_{i=0}^n m_i \widehat{\mathbf{s}_i}^2  -\inertia \right]^{-1} m_j \Bigg[} + \mathbf{s}_j \times  \left\{ \bOm \times \left(\bOm \times \bzeta_j^{\prime} +2 \dot \theta_j \bzeta_j^{\dprime} \right) + \dot \theta_j^2 \bzeta_j^{\thprime} + \ddot \theta_j  \bzeta_j^{\dprime} \right\}  \\
&\hphantom{=\left[\sum_{i=0}^n m_i \widehat{\mathbf{s}_i}^2  -\inertia \right]^{-1} m_j \Bigg[} - 2 \, \Sym\left( \widehat{\mathbf{s}_j} \widehat{{\bzeta}^{\prime}_j} \right) \bkappa\left(t,\bx,\bu \right) 
\Bigg],
\end{split}
\end{equation}
where in the second equality the following result and definition are used:
\begin{equation} \label{eq_hat_rev}
\widehat{\mathbf{a}} \widehat{\mathbf{b}} = \left[ \widehat{\mathbf{b}} \widehat{\mathbf{a}} \right]^\mathsf{T} \quad \forall \, \mathbf{a}, \mathbf{b} \in \mathbf{R}^3
\end{equation}
and 
\begin{equation} \label{eq_Sym}
\Sym\left(A \right) \equiv \frac{1}{2} \left(A+A^\mathsf{T} \right) \quad \forall \, A \in \mathbf{R}^{3 \times 3} . 
\end{equation}

\begin{equation}
\begin{split}
\pp{\bkappa}{\dot \theta_j} &=
\left[\sum_{i=0}^n m_i \widehat{\mathbf{s}_i}^2  -\inertia \right]^{-1}  \left[ m_j \mathbf{s}_j \times  \left\{2 \bOm \times \bzeta_j^{\prime} + 2 \dot \theta_i \bzeta_i^{\dprime}  \right\}  \right] \\
&=2 \left[\sum_{i=0}^n m_i \widehat{\mathbf{s}_i}^2  -\inertia \right]^{-1}  \left[ m_j \mathbf{s}_j \times  \left\{ \bOm \times \bzeta_j^{\prime} + \dot \theta_i \bzeta_i^{\dprime}  \right\}  \right],
\end{split}
\end{equation}

\begin{equation}
\begin{split}
\pp{\bkappa}{\Omega_j} =
\left[\sum_{i=0}^n m_i \widehat{\mathbf{s}_i}^2  -\inertia \right]^{-1}  \left[\mathbf{e}_j \times \inertia \bOm+\bOm \times \inertia \mathbf{e}_j
+ \sum_{i=0}^n m_i \mathbf{s}_i \times  \left\{\mathbf{e}_j \times \left(\bOm \times \bzeta_i +2 \dot \theta_i \bzeta_i^{\prime} \right) + \bOm \times \left(\mathbf{e}_j \times \bzeta_i \right)  \right\}  \right],
\end{split}
\end{equation}

\begin{equation}
\begin{split}
\pp{\bkappa}{q_j} &= -\left[\sum_{i=0}^n m_i \widehat{\mathbf{s}_i}^2  -\inertia \right]^{-1} \left[ \sum_{i=0}^n m_i \left( \widehat{\mathbf{s}_i} \left(r \widehat{ \pp{\bGamma}{q_j}} \right) + \left(r \widehat{ \pp{\bGamma}{q_j}} \right)  \widehat{\mathbf{s}_i}  \right)  \right] \left[\sum_{i=0}^n m_i \widehat{\mathbf{s}_i}^2  -\inertia \right]^{-1} \left[\sum_{i=0}^n m_i \widehat{\mathbf{s}_i}^2  -\inertia \right] \bkappa\left(t,\bx,\bu \right) \\
&\hphantom{=} + \left[\sum_{i=0}^n m_i \widehat{\mathbf{s}_i}^2  -\inertia \right]^{-1}  \Bigg[r \pp{\tilde \bGamma}{q_j} \times \bGamma+r \tilde \bGamma \times \pp{\bGamma}{q_j}\\
&\hphantom{= + \left[\sum_{i=0}^n m_i \widehat{\mathbf{s}_i}^2  -\inertia \right]^{-1}  \Bigg[}+ \sum_{i=0}^n m_i r \pp{\bGamma}{q_j} \times  \left\{ g \bGamma+ \bOm \times \left(\bOm \times \bzeta_i +2 \dot \theta_i \bzeta_i^{\prime} \right) + \dot \theta_i^2 \bzeta_i^{\dprime} + \ddot \theta_i  \bzeta_i^{\prime} \right\} \\
&\hphantom{= + \left[\sum_{i=0}^n m_i \widehat{\mathbf{s}_i}^2  -\inertia \right]^{-1}  \Bigg[}+ \sum_{i=0}^n m_i \mathbf{s}_i \times  g \pp{\bGamma}{q_j} \Bigg] \\
&= \left[\sum_{i=0}^n m_i \widehat{\mathbf{s}_i}^2  -\inertia \right]^{-1}  \Bigg[r \pp{\tilde \bGamma}{q_j} \times \bGamma \\
&\hphantom{= + \left[\sum_{i=0}^n m_i \widehat{\mathbf{s}_i}^2  -\inertia \right]^{-1}  \Bigg[}+\pp{\bGamma}{q_j} \times \left\{\sum_{i=0}^n m_i \left[ -g \bzeta_i + r\left(\bOm \times \left(\bOm \times \bzeta_i +2 \dot \theta_i \bzeta_i^{\prime} \right) + \dot \theta_i^2 \bzeta_i^{\dprime} + \ddot \theta_i  \bzeta_i^{\prime} \right) \right]-r \tilde \bGamma \right\} \\
&\hphantom{= + \left[\sum_{i=0}^n m_i \widehat{\mathbf{s}_i}^2  -\inertia \right]^{-1}  \Bigg[}-2r \sum_{i=0}^n m_i \, \Sym\left( \widehat{\mathbf{s}_i} \widehat{ \pp{\bGamma}{q_j}} \right) \bkappa\left(t,\bx,\bu \right) \Bigg] \\
&= \left[\sum_{i=0}^n m_i \widehat{\mathbf{s}_i}^2  -\inertia \right]^{-1}  \Bigg[ \pp{\bGamma}{q_j} \times \left\{\sum_{i=0}^n m_i \left[ -g \bzeta_i + r\left(\bOm \times \left(\bOm \times \bzeta_i +2 \dot \theta_i \bzeta_i^{\prime} \right) + \dot \theta_i^2 \bzeta_i^{\dprime} + \ddot \theta_i  \bzeta_i^{\prime} \right) \right]-r \tilde \bGamma \right\} \\
&\hphantom{= + \left[\sum_{i=0}^n m_i \widehat{\mathbf{s}_i}^2  -\inertia \right]^{-1}  \Bigg[}+r \left(\pp{\tilde \bGamma}{q_j} \times \bGamma -2 \sum_{i=0}^n m_i \, \Sym\left( \widehat{\mathbf{s}_i} \widehat{ \pp{\bGamma}{q_j}} \right) \bkappa\left(t,\bx,\bu \right) \right) \Bigg],
\end{split}
\end{equation}
where in the second equality \eqref{eq_hat_rev} and \eqref{eq_Sym} are used.

Differentiating $\bGam$, given in \eqref{eq_rball_Gamma}, with respect to $\mathfrak{q}$ yields
\begin{equation} 
\pp{\bGamma}{\mathfrak{q}} = \begin{bmatrix} -2q_2 & 2q_3 & -2q_0 & 2q_1 \\ 2q_1 & 2q_0 & 2q_3 & 2q_2 \\ 0 & -4q_1 & -4q_2 & 0 \end{bmatrix} = 2 \begin{bmatrix} -q_2 & q_3 & -q_0 & q_1 \\ q_1 & q_0 & q_3 & q_2 \\ 0 & -2q_1 & -2q_2 & 0 \end{bmatrix}.
\end{equation}

Differentiating $\tilde \bGamma$, given in \eqref{eq_rball_tGamma}, with respect to $\mathfrak{q}$ yields
\begin{equation} 
\begin{split}
\pp{\tilde \bGamma}{\mathfrak{q}} &= \begin{bmatrix} 2 q_3 \mathbf{F}_\mathrm{e,2}-2q_2\mathbf{F}_\mathrm{e,3} & 2 q_2 \mathbf{F}_\mathrm{e,2}+2q_3\mathbf{F}_\mathrm{e,3} & -4 q_2 \mathbf{F}_\mathrm{e,1}+2q_1\mathbf{F}_\mathrm{e,2}-2q_0\mathbf{F}_\mathrm{e,3} & -4 q_3 \mathbf{F}_\mathrm{e,1}+2q_0\mathbf{F}_\mathrm{e,2}+2q_1\mathbf{F}_\mathrm{e,3} \\
-2 q_3 \mathbf{F}_\mathrm{e,1}+2q_1\mathbf{F}_\mathrm{e,3} & 2 q_2 \mathbf{F}_\mathrm{e,1}-4q_1\mathbf{F}_\mathrm{e,2}+2q_0\mathbf{F}_\mathrm{e,3} & 2 q_1 \mathbf{F}_\mathrm{e,1}+2q_3\mathbf{F}_\mathrm{e,3} & -2 q_0 \mathbf{F}_\mathrm{e,1}-4q_3\mathbf{F}_\mathrm{e,2}+2q_2\mathbf{F}_\mathrm{e,3} \\
2 q_2 \mathbf{F}_\mathrm{e,1}-2q_1\mathbf{F}_\mathrm{e,2} & 2 q_3 \mathbf{F}_\mathrm{e,1}-2q_0\mathbf{F}_\mathrm{e,2}-4q_1\mathbf{F}_\mathrm{e,3} & 2 q_0 \mathbf{F}_\mathrm{e,1}+2q_3\mathbf{F}_\mathrm{e,2}-4q_2\mathbf{F}_\mathrm{e,3} & 2 q_1 \mathbf{F}_\mathrm{e,1}+2q_2\mathbf{F}_\mathrm{e,2} 
\end{bmatrix} \\
&= 2 \begin{bmatrix}  q_3 \mathbf{F}_\mathrm{e,2}-q_2\mathbf{F}_\mathrm{e,3} &  q_2 \mathbf{F}_\mathrm{e,2}+q_3\mathbf{F}_\mathrm{e,3} & -2 q_2 \mathbf{F}_\mathrm{e,1}+q_1\mathbf{F}_\mathrm{e,2}-q_0\mathbf{F}_\mathrm{e,3} & -2 q_3 \mathbf{F}_\mathrm{e,1}+q_0\mathbf{F}_\mathrm{e,2}+q_1\mathbf{F}_\mathrm{e,3} \\
- q_3 \mathbf{F}_\mathrm{e,1}+q_1\mathbf{F}_\mathrm{e,3} &  q_2 \mathbf{F}_\mathrm{e,1}-2q_1\mathbf{F}_\mathrm{e,2}+q_0\mathbf{F}_\mathrm{e,3} &  q_1 \mathbf{F}_\mathrm{e,1}+q_3\mathbf{F}_\mathrm{e,3} & - q_0 \mathbf{F}_\mathrm{e,1}-2q_3\mathbf{F}_\mathrm{e,2}+q_2\mathbf{F}_\mathrm{e,3} \\
q_2 \mathbf{F}_\mathrm{e,1}-q_1\mathbf{F}_\mathrm{e,2} &  q_3 \mathbf{F}_\mathrm{e,1}-q_0\mathbf{F}_\mathrm{e,2}-2q_1\mathbf{F}_\mathrm{e,3} &  q_0 \mathbf{F}_\mathrm{e,1}+q_3\mathbf{F}_\mathrm{e,2}-2q_2\mathbf{F}_\mathrm{e,3} &  q_1 \mathbf{F}_\mathrm{e,1}+q_2\mathbf{F}_\mathrm{e,2} 
\end{bmatrix}.
\end{split}
\end{equation} 

\rem{ 
	\todo{SMR: I removed the R2Q2 revision command from Section 2.2, because it is too hard to edit. I will add back the R2Q2 revision command  when the paper is ready for resubmission. I went ahead and made lots of minor grammar, English, and punctuation edits, though I think the grammar is still choppy. Other edits are marked with green textcolor. Perhaps Section~\ref{sec:intrinsic} should be moved to the appendix, since this is not what the reviewer meant by intrinsic coordinates and since you discuss both uncontrolled and controlled equations of motion, whereas Section 2 is on uncontrolled equations of motion. \\ 
		VP: I think it is close to what the reviewer wants. We can't really go all the way with abstract notation, and the use of $SO(3) \times \mathbb{R}^2$ is not really appropriate for us, but it is as close to Jurdjevic as we can get. }
} 

\revision{R2Q2}{\section{Equations of Motion in Intrinsic Coordinates} \label{app_intrinsic} } 
Let us consider the equations of motion for the rolling ball from the general geometric point of view and rederive
\eqref{uncon_ball_eqns_explicit_1d} using Euler-Poincar\'e's variational principle \cite{Ho2011_pII,MaRa2013}. The geometric approach yields an elegant and efficient way to treat optimal control problems, especially when applied to rolling bodies, e.g. see \cite{Jurdjevic1993,Jurdjevic1999a,Jurdjevic1999b,ohsawa2020geometric} and subsequent application of geometric theory to the trajectory planning of bodies rolling with and  without twisting at the contact point \cite{AlChLo2010}. In particular, the celebrated paper of Jurdjevic  \cite{Jurdjevic1993}  demonstrates the analytical solution of the equations for the kinematically controlled ball  rolling between two plates using the description of the rolling ball in terms of the Lie group $SO(3) \times \mathbb{R}^2$. The variables describing the motion are the ball's orientation and the position of the ball's geometric center above the plane. More recently, using this description, \cite{ohsawa2020geometric} found the analytical solution of the Pontryagin optimal control problem for the \emph{kinematic} control case for the ball driven by internal wheels, such as the Sphero toy. 

Our problem is a bit more complex than Jurdjevic's case in the sense that the control problem we consider is dynamic, not kinematic. The equations of motion are thus considerably more complicated than the case of kinematic control. Since the distance between the center of mass of the ball and the plane it is rolling on can vary with time, it is more convenient for us to work with the Lie group of  translations and rotations in the space $SE(3)=SO(3) \circledS \mathbb{R}^3$, where the notation $\circledS$ denotes the semidirect product of two Lie groups. The precise definition of the semidirect product Lie group is as follows \cite{Ho2011_pII}. 
Suppose a Lie group $G$ may be decomposed uniquely into a normal subgroup $K$ and a subgroup $H$ such that every group element may be written as
$g = kh$ or $g = hk$, in either order, for unique choices of $k \in K$ and $h \in H$. Then, $G$ is called a semidirect product
of $H$ and $K$ and denoted as $G=H \circledS K$. When $K$ is a vector space, the action of $G$ on itself is given by 
\begin{equation} 
\label{semidirect_prod_action} 
(h_1,k_1)(h_2,k_2)=(h_1 h_2, h_1 k_2 + k_1). 
\end{equation} 
The Lie group of rotations and translations $SE(3) = SO(3) \circledS \mathbb{R}^3$ is a particular and familiar example of a semidirect product group which is important for the dynamics and control applications considered here. We call the coordinates of that group the \emph{intrinsic coordinates}.  In this section, we shall show that the derivation of the dynamical equations can be realized quite elegantly using the group-theoretical description; however, the question of the optimal control is more complex and is difficult to accomplish without addressing particular properties of $SE(3)$. 

\rem{ 
	\todo{VP: Technically, \cite{AlChLo2010} states that they have no spinning (sticky ball), as far as I can see. Not sure how important that is, I did not find it very essential to their derivation (although, I must say, I have not used too much time on this). Maybe we should just cross out spinning above. Also, I found their IEEE paper from 2010, there seems to be nothing from 2008. Did I miss it online? Can you have a look? \\ SMR: I only found the 2010 IEEE paper. Above I see that you wrote slipping but not spinning or twisting. Above you write $SE(3)=SO(3) \times \mathbb{R}^3$, but below you write $SE(3) = SO(3) \circledS \mathbb{R}^3$.  \\ 
		VP: OK, we will keep that citation. I cleaned up the text and put the definition of semidirect product group.  }
	\todo{SMR: Should we use $\mathbf{z}_0$ instead of $\br$ below for the position of the ball's center of mass? Why do you use $\bA_i(t)$ instead of $\bxi_i(t)$, which was the notation used in the dynamics paper? 
		\\ 
		VP: In principle we could, but everyone is calling $SE(3)=(\Lambda, \br)$ so it may confuse the readers... I kept it as is for now, changed $\bA$ to $\bxi$. I did not highlight it with magenta.
} } 

The orientation and position of the center of mass of the ball with respect to a given spatial coordinate system is given by $g=(\Lambda, \mathbf{z}_0) \in SE(3)$, with an additional constraint that the ball is spherical so that the distance between the geometric center and the horizontal surface is constant. For $1 \le i \le n$, the $i^\mathrm{th}$ mass moves with respect to the ball according to a prescribed trajectory $\bxi_i(t) \in \mathbb{R}^3$. Thus, to write an equation in intrinsic coordinates, we consider a Lie group $G$, with its Lie algebra $\mathfrak{g}$, acting on a vector space $V$, with a set of variables $\xi_i \in V$, which we call $\overline{\xi}=\left(\xi_1,\ldots,\xi_n\right)$. The unreduced Lagrangian depends on the variables $\left(g,\dot g\right)$, as well as the positions and velocities of the masses 
\begin{equation} 
\label{mass_pos_vel_gen} 
\overline{\rm{a}}=g \overline{\xi} =\left(g \xi_1, \ldots g \xi_n\right), 
\quad 
\overline{\rm{v}}= \dot {\overline{\rm{a}} } = \left(\dot g \xi_1 + g \dot{\xi}_1, \ldots, \dot g \xi_n + g \dot{\xi}_n \right). 
\end{equation} 
Here, for $v_i={\dot \xi}_i \in T_{\xi_i} V$ we have defined $g v_i$ to be the corresponding action of the Lie group $G$ on the vector $v_i$. To make this discussion more concrete, for the case of $G=SE(3)$ which we will consider further in this appendix, we have $g=(\Lambda, \mathbf{z}_0)$, the orientation of the ball and the position of its center of mass in $\mathbb{R}^3$. Then, the mass positions and velocities in the body frame are $\bxi_i \in \mathbb{R}^3$ and $\bv_i = \dot{\bxi}_i \in \mathbb{R}^3$. The mass  positions and velocities in the spatial frame are given by \eqref{mass_pos_vel_gen}, which are written explicitly as $\mathbf{a}_i = g \bxi_i = \Lambda \bxi_i + \mathbf{z}_0$ (position) and $\mathbf{v}_i=\dot{\mathbf{a}}_i$ (velocity). 

The action of $g$ on the velocity is given by $g \bv_i  = \Lambda \dot{\bxi}_i$. If we call $\alpha = g^{-1} \dot g \in \mathfrak{g}$, i.e. $\alpha$ is an element of the Lie algebra of $G$, then $\dot g \xi_i$ is computed as 
$\dot g \xi_i = g  (\alpha \xi_i )$. In the case of $g=(\Lambda,\mathbf{z}_0)$, the corresponding Lie algebra element is $\alpha=(\bOm, \Lambda^{-1} \dot{\mathbf{z}}_0)$, with $\bOm = (\Lambda^{-1} \dot \Lambda)^\vee$, and  the corresponding action is computed as $\alpha \bxi_i=\bOm \times \bxi_i + \Lambda^{-1} \dot{\mathbf{z}}_0$. Thus, for the case of $g \in SE(3)$, 
\begin{equation} 
\mathbf{v}_i=\dot g \bxi_i + g {\dot \bxi}_i = \dot \Lambda \bxi_i + \Lambda {\dot \bxi}_i + \dot{\mathbf{z}}_0. 
\label{v_SE3}
\end{equation}
We shall return to the discussion of $G=SE(3)$ later in this section when considering the explicit expression for the equations written in the general form. 
\rem{ 
	\todo{SMR: Is \eqref{mass_pos_vel_gen} correct? For $g=(\Lambda, \br) \in SE(3)$, the time derivative of $gA_i = \Lambda A_i+\br$ is $\dot \Lambda A_i+\Lambda \dot A_i + \dot \br$, which does not equal $\dot g A_i + g \dot{A}_i = \dot \Lambda A_i + \dot \br+\Lambda \dot A_i + \br.$ \\ 
		VP: Actually, it is if you consider it as group actions: $\dot g (\bxi)$ is defined as $\dot \Lambda \bxi + \dot \br= \Lambda \left( \bOm \times \bxi \right) + \dot \br$. On the other hand, $g ( \dot \bxi) = \Lambda \dot \bxi$ since it transfers velocities to velocities, the transformation of velocities does not add $\br$ to velocities. So these formulas agree. I tightened up explanation above. \\ 
		VP: I also used $a$ and $A$ to distinguish between the body and spatial variables, because the capital $\xi$ (i.e. $\Xi$) looks really ugly, I think. In the text below, I kept $a$ but changed $A$ to $\xi$.  } 
} 

In what follows, we shall assume that both the Lagrangian and the constraints are invariant with respect to the Lie group $G$, as is the case of a single rolling ball. The general theory of nonholonomic systems outlined here was first derived in the Hamiltonian framework in \cite{bates1993nonholonomic} and in the Lagrangian framework in  \cite{BlKrMaMu1996}, see also \cite{Bloch2003} for a more detailed exposition. In the general case, the Lagrangian is written as $L=L\left(g, \dot g, \overline{\rm{a}}, \overline{\rm{v}}\right)$, with the symmetry-reduced Lagrangian 
\begin{equation} 
\ell=\ell\left(\Omega, \overline{V}, \overline{\xi}\right), \quad  \mbox{with} \quad 
\Omega \equiv g^{-1} \dot g \in \mathfrak{g} \quad \mbox{and} \quad \overline{V} \equiv \Omega \overline{\xi}+  \dot {\overline{\xi}} \, . 
\label{sym_reduced_lagr}
\end{equation} 

Note that in the case of $G=SE(3)$, we have $\Omega=\left(\left(\Lambda^{-1} \dot \Lambda\right)^\vee, \Lambda^{-1} \dot{\mathbf{z}}_0\right)$ as above, so $\Omega=g^{-1} \dot g$ involves both angular and linear velocities. Technically speaking, the Lagrangian of the rolling ball is only invariant with respect to rotations about the vertical axis and translations in $\mathbb{R}^3$, though it can be viewed as being a $SE(3)$-invariant quantity where gravity is the advected parameter \cite{Ho2011_pII,Bloch2003}. Let us define the variations $\eta=g^{-1} \de g \in \mathfrak{g}$. Since $\overline{\xi}$ in this case consists of prescribed functions of $t$, the components of $\overline{\xi}$ do not need to be varied to obtain the uncontrolled dynamics (they will, however, be varied to obtain the controlled dynamics). The variations of $\Omega$ and $\overline{V}$ are given by \cite{Ho2011_pII,MaRa2013}:  
\rem{\todo{VP: I think it is better to cross out $g$, technically, the variations are taken with respect to $\eta$ which is an element of the Lie algebra $\eta$. Writing 'With respect to $g$' may be confusing. \\ SMR: I thought the variations were taken with respect to $g$, after which one defined $\eta = g^{-1} \de g$. I think we argued about this when I was writing my thesis. 
		\\
		VP: The correct statement of Euler-Poincar\'e Theorem is as follows: suppose due to symmetry the Lagrangian $L(g,\dot g)$ can be written as a function which depends on $\xi=g^{-1} \dot g$, so $\ell=\ell(\xi)$. Then, the Euler-Lagrange equations 
		\[ 
		\delta \int L(g,\dot g) \mbox{d} t = \int \left( \frac{d}{dt} \pp{L}{\dot g}- \pp{L}{g} \right) \de g =0 
		\] 
		are equivalent to 
		\[ 
		\de \int \ell(\xi) \mbox{d} t = \int \left< \pp{\ell}{\xi} \, , \, \de \xi\right> \mbox{d} t
		\]
		taken on the variations $\de \xi = \dot \eta + {\rm ad}_\xi \eta$, with $\eta = g^{-1} \de g$. So if we are talking about the application of Euler-Poincar\'e theory, the variations are with respect to $\xi$ with a particular form of variations in $\xi$. If you are talking about Euler-Lagrange equations, then variations are with respect to $g$, taken of the unreduced Lagrangian $L(g, \dot g)$ with no restriction on $\de g$. 
		We can probably just say 'the variations are computed as' and not go into these details because otherwise we will need to specify a lot of other information. So I am OK with crossing out the text above. 
} }
\begin{equation} 
\de \Omega = \dot \eta +{\rm ad}_\Omega \eta  \quad \mbox{and} \quad \de \overline{V} = \de \Omega \overline{\xi} = \left(\dot \eta +{\rm ad}_\Omega \eta \right) \overline{\xi}, 
\label{var_Omega_gen} 
\end{equation} 
where ${\rm ad}$ is the adjoint action of the group $G$ acting on itself. 
Then, using the critical action principle $\de \int \ell\left(\Omega, \overline{V} , \overline{\xi}\right) \dt=0$ on variations $\de \Omega$ satisfying \eqref{var_Omega_gen} yields the Euler-Poincar\'e equations of motion: 
\begin{equation} 
\label{EP_deriv_gen_init} 
\begin{aligned} 
0 = & \int \left[ \left<  \pp{\ell}{\Omega} , \de \Omega \right>_{\mathfrak{g}} + 
\left< \sum_{i=1}^n \pp{\ell}{V_i} , \de \Omega \xi_i \right> _{V} \right] \dt 
= \int  \left<  \pp{\ell}{\Omega} -\sum_{i=1}^n \pp{\ell}{V_i} \diamond \xi_i , \de \Omega \right>_{\mathfrak{g}} \dt. 
\end{aligned} 
\end{equation} 
In \eqref{EP_deriv_gen_init}, the action of $\alpha \in \mathfrak{g}$  on $b \in V$ is expressed as $\alpha b$; for example,  if  $\alpha=g^{-1} \dot g$, then $\alpha b =g^{-1} \dot g b$. In \eqref{EP_deriv_gen_init}, $\diamond$ is the diamond operator which acts on $a \in V^*$ and $b \in V$. The diamond operator is defined, in the general case, as minus the dual of the Lie derivative action \cite{Ho2011_pII}. In our case, consider the pairing $\left< \cdot ,\cdot \right>_{\mathfrak{g}}$ between elements of $\mathfrak{g}$ and $\mathfrak{g}^*$ and the pairing $\left< \cdot ,\cdot \right>_{V}$ between elements of $V$ and $V^*$.  Then, the diamond operator is defined as 
\begin{equation} 
\label{Diamond_def}
\left< a \diamond b, \alpha \right>_{\mathfrak{g}} = - \left<a, \alpha b \right>_V \, . 
\end{equation}
To continue, we use \eqref{var_Omega_gen} and integrate \eqref{EP_deriv_gen_init} by parts to isolate the coefficient multiplying $\eta$, which yields the Euler-Poincar\'e equations of motion 
\begin{equation} 
\dd{}{t} \Pi = {\rm ad}^*_\Omega \Pi \, , \quad \Pi \equiv \pp{\ell}{\Omega} - \sum_{i=1}^n \pp{\ell}{V_i} \diamond \xi_i \, . 
\label{EP_gen} 
\end{equation} 
Here, the ${\rm ad}^*$ operator acting on elements of $\mathfrak{g}^*$ is defined as $\left<{\rm ad}^*_\Omega a, b\right>_{\mathfrak{g}} = \left<a, {\rm ad}_\Omega b\right>_{\mathfrak{g}}$ for any $\Omega, b \in \mathfrak{g}$.  

The linear nonholonomic constraint defined on the body angular velocity  can be reformulated as  
\begin{equation} 
\left< \alpha(g) , \Omega \right>_\mathfrak{g}=0 \, , 
\label{nonholonomic_gen} 
\end{equation} 
where the function $\alpha(g)$ takes its value in the dual of the Lie algebra. The case of several constraints can be considered trivially by repeating \eqref{nonholonomic_gen} for several functions $\left\{ \alpha^j(g) \right\}_{j=1}^m$. In what follows, we concentrate only on the case of a single  constraint \eqref{nonholonomic_gen}, since this is the case relevant to the rolling ball. The equations of motion are derived using Euler-Poincar\'e-Suslov's (also known as Lagrange-d'Alembert's) principle using the variational principle above and imposing an additional constraint on the variations $\eta = g^{-1} \de g$: 
\begin{equation} 
\left< \alpha(g) , \eta \right>_\mathfrak{g}=0  \, .
\label{nonholonomic_gen_var} 
\end{equation} 
In that case, the equations of motion are derived by enforcing the constraint \eqref{nonholonomic_gen_var} with a Lagrange multiplier $\lambda$  \cite{Bloch2003} and are written in the form: 
\begin{equation} 
\dd{}{t} \Pi  = {\rm ad}^*_\Omega \Pi +  \lambda \alpha (g) \, , \quad \Pi \equiv \pp{\ell}{\Omega} - \sum_{i=1}^n \pp{\ell}{V_i} \diamond \xi_i \,, 
\label{nonholonomic_eqs_gen} 
\end{equation} 
coupled with the constraint \eqref{nonholonomic_gen}, which determines the Lagrange multiplier  $\lambda$. 

The computation  preceding  \eqref{nonholonomic_eqs_gen} above is valid for an arbitrary Lie group $G$. Let us now consider a particular case when the group $G$  has the  structure of a semidirect product, such as $SE(3)$, in which an element $g =\left(\Lambda, \mathbf{z}_0\right) \in SE(3)$ acts on elements $\bx \in V = \mathbb{R}^3$ through the law $g \bx = \Lambda \bx + \mathbf{z}_0$. 

\rem{ 
	Suppose a Lie group $H$ acts from the left by linear maps on the linear space $V$. We say that a group $G$ is a semidirect product and write $G=H \circledS V$ if $G=H \times V$ and the group multiplication of elements $\left(h_1,v_1\right) \in G$ and $\left(h_2,v_2\right) \in G$ is given by 
	\begin{equation} 
	\label{semidirect_prod_action} 
	\left(h_1,v_1\right) \left(h_2,v_2\right) = \left(h_1 h_2, v_1 + h_1 v_2\right) \,. 
	\end{equation} 
} 

For the particular case of $G=SE(3) = SO(3) \circledS \mathbb{R}^3$, the adjoint ${\rm ad}$ and coadjoint ${\rm ad}^*$ actions of $G$ on itself are written explicitly as follows. The Lie algebra $\mse(3)$ of $SE(3)$ can be associated with pairs of vectors in $\mathbb{R}^3$ through the hat map as either $\left(\balpha, \boldsymbol{\beta}\right)$ or $\left(\widehat{\balpha}, \boldsymbol{\beta}\right)$, with the hat map  $\widehat{\balpha}_{ij}=-\epsilon_{ijk} \alpha^k$ identifying the mapping between $3 \times 3$ antisymmetric matrices $\widehat{\balpha}$ and vectors $\balpha \in \mathbb{R}^3$. Taking elements $\left(\balpha_1, \boldsymbol{\beta}_1\right)$ and $\left(\balpha_2, \boldsymbol{\beta}_2\right)$ in $\mse(3)$ (the Lie algebra of $SE(3)$) and $\left(\bu,\bv\right) \in \mse(3)^*$, the adjoint and coadjoint actions are \cite{Ho2011_pII}:
\begin{equation} 
\begin{aligned} 
& {\rm ad}_{\left(\balpha_1, \boldsymbol{\beta}_1\right)} \left(\balpha_2, \boldsymbol{\beta}_2\right) = \left(   \balpha_1 \times \balpha_2  , \balpha_1 \times \boldsymbol{\beta}_2- 
\balpha_2 \times \boldsymbol{\beta}_1\right) , 
\\
& \left< {\rm ad}^*_{\left(\balpha_1, \boldsymbol{\beta}_1\right)} \left(\bu,\bv\right), \left(\balpha_2, \boldsymbol{\beta}_2\right) \right> 
= 
\left( \bu\times \balpha_1 + \bv \times \boldsymbol{\beta}_1 \right) \cdot \balpha_2 + \left( \bu \times \balpha_1 \right) \cdot \boldsymbol{\beta}_2 . 
\end{aligned} 
\quad 
\label{ad_coad_actions} 
\end{equation} 
The diamond operator in \eqref{Diamond_def} can also be computed explicitly for the case of $SE(3)$. If $g=\left(\Lambda, \mathbf{z}_0\right) \in SE(3)$, then $g \bxi_i = \Lambda \bxi_i + \mathbf{z}_0$. Consider  $\Om = g^{-1} \dot g \in \mse(3)$ with 
$\Om \bxi_i = g^{-1} \dot g \bxi_i$. If we denote $\Omega=\left(\balpha, \boldsymbol{\beta}\right)$, then, physically, $\balpha$ and  $\boldsymbol{\beta}$  are the angular and linear velocities measured in the body frame. Take an arbitrary vector $\bu \in V^*=\left(\mathbb{R}^3\right)^*$ and the natural pairing between $\mathbb{R}^3$ and its dual. Then, 
\begin{equation} 
\begin{split}
\left< \bu , \Omega \bxi_i \right>_V= \bu \cdot \left( \balpha \times \bxi_i + \boldsymbol{\beta} \right) = \balpha  \cdot \left( \bxi_i \times \bu \right) + \boldsymbol{\beta} \cdot \bu =&
- 
\left< -\left( \bxi_i \times \bu, \bu\right), \left(\balpha, \boldsymbol{\beta}\right)  \right>_{\mse(3)} 
\\ = &-\left< \bu \diamond \bxi_i, \left(\balpha, \boldsymbol{\beta}\right)  \right>_{\mse(3)}. 
\label{Diamond_calc} 
\end{split}
\end{equation} 
Here, we brought the minus sign to the front of the bracket on the right-hand side of \eqref{Diamond_calc} to match the definition of the diamond operator \eqref{Diamond_def}. 

We shall note that for the case of a rolling ball $G=SE(3)$, a semidirect product group.  The single  Lagrange multiplier $\lambda$ for the nonholonomic constraint in \eqref{nonholonomic_eqs_gen} can be eliminated, since the constraint has a particular form. More specifically, $\Omega=\left(\bOm, \bY_0\right) \in \mse(3)$, with $\bOm$  and $\bY_0$ denoting the angular and linear velocities in the body frame, as discussed above. The nonholonomic constraint then has the form $\bOm \times \bs_0 = \bY_0$, where $\bs_0 \in \mathbb{R}^3$ is the body frame vector connecting the contact point with the center of mass.  The nonholonomic constraint then leads to a constraint on the variations $(\bSigma,\bPsi)  \equiv \left( \left(\Lambda^{-1} \de \Lambda\right)^\vee , \Lambda^{-1} \de \mathbf{z}_0 \right) \in \mse(3)$ of the form $\bPsi = \bSigma \times \bs_0 $. Letting $\boldsymbol{\lambda}$ denote the Lagrange multiplier enforcing this constraint on the variations, Lagrange-d'Alembert's principle for this system may be written as 
\begin{equation} 
\int \left[ \left< \mbox{Linear momentum balance} , \bPsi \right> +
\left< \mbox{Angular momentum balance}, \bSigma \right> + \left<\boldsymbol{\lambda}, \bPsi - \bSigma \times \bs_0 \right> \right] \dt =0. 
\label{LdA_gen_ball} 
\end{equation} 
In \cite{Putkaradze2018dynamicsP,putkaradze2018normalpub}, we substituted the expression $\bPsi = \bSigma \times \bs_0$ directly into \eqref{LdA_gen_ball} and collected the terms proportional to $\bSigma$ to obtain the equations of motion \eqref{uncon_ball_eqns_explicit_1d}. Alternatively, from \eqref{LdA_gen_ball}, the equations of motion are 
\begin{equation} 
\mbox{Linear momentum balance} = -\boldsymbol{\lambda} \, , \quad 
\mbox{Angular momentum balance}= \bs_0 \times \boldsymbol{\lambda} \, . 
\label{Eqs_with_lambda} 
\end{equation}
Then, the Lagrange multiplier $\boldsymbol{\lambda}$ can be eliminated by vector multiplying the linear momentum equation by $\bs_0$ and adding that result to the angular momentum equation, yielding \eqref{uncon_ball_eqns_explicit_1d}. This latter procedure is essentially the same as the standard procedure for eliminating Lagrange multipliers from nonholonomic systems of rolling ball type, see \cite{Bloch2003,Ho2011_pII}. However, elimination of the Lagrange multiplier for motion on general Lie groups and other types of constraints does not seem to be possible. 

One can also achieve further progress in the more general case by considering $G$ to be a general semidirect product Lie group. When the masses are not moving, it is not necessary to consider the additional action of this group on a vector space. One can also extend the theory developed here to the case where some of the internal masses  are not explicitly controlled, but rather exhibit additional dynamics, e.g. some of them are movable, some are connected with elastic springs and/or rods, or some are  influenced by prescribed forces and torques with their exact positions unknown. In that case, the variations of the action with respect to the  subset of dynamic variables in $\overline{\xi}$ will yield a modification of \eqref{nonholonomic_eqs_gen}. That extension of \eqref{nonholonomic_eqs_gen} can be accomplished with relative ease in the geometric setting, in contrast to generalizing \eqref{uncon_ball_eqns_explicit_1d} explicitly, which is bound to be very difficult.

Of course, the actual implementation of these equations from their general form  is still going to be cumbersome, requiring careful calculation of the adjoint actions and diamond operators. However, it seems more straightforward than the direct variational calculation and can be advantageous for further work in this area. The derivation presented here generalizes the uncontrolled dynamics of a rolling ball actuated by moving masses \cite{Putkaradze2018dynamicsP,putkaradze2018normalpub} to the case of motion on a general Lie group, and, as far as we are aware, has not been presented in this generality. 

For further discussion of the optimal control between times $t=a$ and $t=b$, we enforce the uncontrolled equations of motion \eqref{nonholonomic_eqs_gen} via a Lagrange multiplier $\kappa \in \mathfrak{g}$ and consider the  augmented performance index 
\begin{equation} 
\begin{aligned} 
\tilde J=& \int_a^b \left[ C\left( g, \Omega, \overline{\xi}, \dot{\overline{\xi}}, \ddot{\overline{\xi}} \right) - \left< \kappa, \dd{}{t}\Pi  - {\rm ad}^*_\Omega \Pi -   \lambda \alpha  (g)  \right>_{\mathfrak{g}} \right] \dt
\\ =&  
\int_a^b \left[ C\left(g, \Omega, \overline{\xi}, \dot{\overline{\xi}}  ,\ddot{\overline{\xi}} \right) + \left< \dot \kappa,   \Pi \right>_{\mathfrak{g}}  + 
\left<\kappa,  {\rm ad}^*_\Omega \Pi + \lambda  \alpha  (g)  \right>_{\mathfrak{g}} \right] \dt  - \left.  \left< \kappa,   \Pi \right>_{\mathfrak{g}} \right|_a^b  \, ,
\\ 
\quad \Pi \equiv&  \pp{\ell}{\Omega} - \sum_{i=1}^n \pp{\ell}{V_i} \diamond \xi_i.
\end{aligned} 
\label{Cost_opt_control} 
\end{equation} 
We note that the performance index may depend explicitly on $g$ if, for example, it is desired to follow a given trajectory on the plane or avoid obstacles. Alternatively, it may be desirable to have a camera installed on a rolling ball robot to point in a particular direction in space while rolling along its trajectory. The performance index will depend explicitly on the position of the ball in space (first case) and/or the ball's orientation (second case), i.e. the group element $g$. 

\rem{ 
	\todo{SMR: In the augmented performance index $\tilde J$ above, I think the integrand cost function $C$ should also depend on $\ddot {\overline{\xi}}$. \\ 
		VP: Sure. \\  
		SMR: In the line below, do you prefer to write variations w.r.t. $\overline{\xi}$ or $\xi_i$? \\ 
		VP I think we can do both, i.e. explain that we are doing variations with respect to every component. }
} 

Then, taking variations of the augmented performance index in \eqref{Cost_opt_control} with respect to $g$ or $\Omega$, if the performance index \eqref{Cost_opt_control} does not depend explicitly on $g$, and  $\overline{\xi}$ (i.e. each component $\xi_i$ for $1\le i \le n$) yields the costate  (or adjoint) equations, whereas taking variations with respect to $\kappa$ enforces the uncontrolled equations of motion \eqref{nonholonomic_eqs_gen}. One can also include the reconstruction equations for the spatial coordinate for problems involving obstacle avoidance or desired orientations. We have found out that the resulting equations, while being treatable, are quite nonintuitive (such as involving diamonds and second derivatives of the Lagrangian with respect to the variables) and involve too many technicalities for the purposes of this paper, so in the main text we  use the explicit calculations on $SE(3)$ for these derivations. The  calculations we present there are fundamentally the same in spirit as the calculations based on \eqref{Cost_opt_control}, but are more friendly towards applications without involving unnecessarily complex mathematical notation. 

\rem{ 
	\todo{VP: The biggest problem I found was that in general, unless additional assumptions are being made, such as $G=SE(3)$, and certain form of constraints, there is no way to remove 
		the Lagrange multiplier. Doing optimal control with that Lagrange multiplier is very difficult. If you assume the particular forms of contraint and $G$, you can just as well derive the equations the way we did it. Did I say it OK there? I did not want to argue strongly with Reviewer 2, but it seems to me that there is not that much value in going too abstract. \\ 
		VP: I think it will be good to keep it here, since the complete derivation through the intrinsic equations is not possible. Jurdjevic, Tomoki and others have only done kinematic control with intrinsic coordinates. So I think it is important to say why it is not possible to go all the way. } 
} 

\end{document}